%% file: Dissertation.tex
\pdfoutput=1
\documentclass[oneside, 11pt]{amsbook}
\input{LayoutPreamble.tex}

\input{MacrosPreamble.tex}

\begin{document}
   \frontmatter

   \pagestyle{prelim}
   
   %
   \fancypagestyle{plain}{%
      \fancyhf{}
      \cfoot{-\thepage-}
   }%
   \input{TitlePage.tex}
   \newpage
   
   %
   \doublespacing
   
   \tableofcontents
   \newpage
   
   \input{Abstract.tex}
   \newpage
   
   \section*{Acknowledgments}
   \input{Acknowledgments.tex}

   \mainmatter
   
   \pagestyle{maintext}
   
   %
   \fancypagestyle{plain}{%
      \renewcommand{\headrulewidth}{0pt}
      \fancyhf{}
      \rhead{\thepage}
   }%
   
   \chapter{Introduction}
   \label{ch:IntroductionLabel}
   \input{Introduction.tex}

   \chapter[%
      Idempotents in $H_0(S_N)$
   ]{%
      A Combinatorial Formula for Idempotents in the Zero-Hecke Algebra of the Symmetric Group
   }%
   \label{ch:idempotents}
   \input{chapter-idempotents.tex}

   \chapter[%
      Representation Theory of $\JJ$-Trivial Monoids
   ]{%
      Representation Theory of $\JJ$-Trivial Monoids
   }%
   \label{ch:jtrivial}
   \input{chapter-jtrivial.tex}

  \chapter[%
      $\NDPF$ and Pattern Avoidance
   ]{%
      Non-Decreasing Parking Functions and Pattern Avoidance
   }%
   \label{ch:ndpfavoid}
   \input{chapter-ndpfavoid.tex}

  \chapter[%
      Crystal Bases
   ]{%
      Some Results on Crystal Bases
   }%
   \label{ch:crystals}
   \input{chapter-crystals.tex}


       
   \backmatter
   
   \bibliographystyle{amsalpha-fi-arxlast}
   \bibliography{DissertationBibliography}
\end{document}

%% file: LayoutPreamble.tex
\usepackage[
   includehead,
   includefoot,
     left = 1.5in, 
      top = 0.5in, 
    right = 1in,
   bottom = 1in
]{geometry}
\usepackage{fancyhdr}
\usepackage{setspace}
\usepackage{calc}

\fancyheadoffset[R]{0.5in} 

\fancypagestyle{prelim}{%
   \renewcommand{\headrulewidth}{0pt} 
   \fancyhf{}           
   \pagenumbering{roman}    
   \cfoot{-\thepage-}       
}

\fancypagestyle{maintext}{%
   \renewcommand{\headrulewidth}{0.4pt}
   \pagenumbering{arabic}
   \fancyhf{}
   \fancyhead[L]{\rightmark}
   \rhead{\thepage}
}

\numberwithin{figure}{chapter} 
\numberwithin{table}{chapter}
\numberwithin{equation}{chapter}
\numberwithin{section}{chapter}

%% file: MacrosPreamble.tex
\usepackage{graphicx}
\usepackage{verbatim}
\usepackage{tikz}
\usetikzlibrary{matrix,shapes}

\usepackage{amsmath,amsthm,amssymb}

\usepackage[numbers]{natbib}
\defcitealias{bona.permutations}{\scshape [Bo04]}

\newtheorem{theorem}{Theorem}[section]
\newtheorem{lemma}[theorem]{Lemma}
\newtheorem{proposition}[theorem]{Proposition}
\newtheorem{corollary}[theorem]{Corollary}
\newtheorem{definition}[theorem]{Definition}
\newtheorem{conjecture}[theorem]{Conjecture}
\newtheorem{problem}[theorem]{Problem}

\newenvironment{remark}[1][Remark]{\begin{trivlist}
\item[\hskip \labelsep {\bfseries #1}]}{\end{trivlist}}

\theoremstyle{definition}
\newtheorem{example}[theorem]{Example}

\newcommand{\nc}{\newcommand}
\nc{\pv}{P^{\vee}}

\newcommand{\bpi}[1]{\pi_{#1}^{-}}
\newcommand{\barp}{\bar{p}}
\newcommand{\barx}{\bar{x}}
\newcommand{\bary}{\bar{y}}

\newcommand{\joins}{\operatorname{Joins}}

\newcommand{\id}{\operatorname{id}}
\newcommand{\im}{\operatorname{im}}
\newcommand{\longest}{\pi}
\newcommand{\K}{\mathbb{K}}
\newcommand{\ZZ}{\mathbb{Z}}
\newcommand{\len}{\ell}
\newcommand{\wt}{\operatorname{wt}}

\newcommand{\opi}{\overline{\pi}}
\newcommand{\rad}{\operatorname{rad}}
\newcommand{\rank}{\operatorname{rank}}
\newcommand{\sg}[1][n]{{\mathfrak{S}_{#1}}}
\newcommand{\suchthat}{\mid}
\newcommand{\tMonoid}{M}

\newcommand{\N}{\mathbb{N}}

\newcommand{\CC}{\mathbb{C}}

\newcommand{\QQ}{\mathcal{Q}}
\newcommand{\RR}{\mathcal{R}}
\newcommand{\LL}{\mathcal{L}}
\newcommand{\JJ}{\mathcal{J}}
\newcommand{\HH}{\mathcal{H}}
\newcommand{\BB}{\mathcal{B}}
\newcommand{\KK}{\mathcal{K}}

\newcommand{\lfix}[1]{\operatorname{lfix}(#1)}
\newcommand{\rfix}[1]{\operatorname{rfix}(#1)}
\newcommand{\idempMon}[1][\tMonoid]{E(#1)} 

\newcommand{\biheckemonoid}[1][]{M}

\newcommand{\OR}{\mathcal{OR}}
\newcommand{\NDPF}{{\operatorname{NDPF}}}
\newcommand{\BNDPF}{{\operatorname{BNDPF}}}
\newcommand{\ANDPF}{{\operatorname{\widetilde{NDPF}}}}
\newcommand{\ASn}{{\widetilde{S}}}
\newcommand{\unitribool}{\mathcal U}
\newcommand{\sgn}{{\operatorname{sgn}}}
\newcommand{\lex}{{\operatorname{lex}}}

\newcommand{\edge}{\!\!\rightarrow\!\!}

\newcommand{\leftexp}[2]{{\vphantom{#2}}^{#1}{#2}}
\newcommand{\ggc}{ \mathfrak{g} }
\newcommand{\ggh}{ \hat{\mathfrak{g}} }
\newcommand{\sln}{ \mathfrak{sl} }
\newcommand{\slnh}{ \hat{\mathfrak{sl}} }
\newcommand{\pr}{ \operatorname{pr} }
\newcommand{\cpr}{ \operatorname{\mathfrak{pr}} }

\usepackage{ifthen}
\newboolean{draft}
\setboolean{draft}{true}
\ifdraft
\newcommand{\TODO}[2][To do: ]{\textcolor{red}{\textbf{#1#2}}}
\else
\newcommand{\TODO}[2][]{}
\fi

%% file: TitlePage.tex
\begin{center}
   \null\vfill
   \textbf{%
      Excursions into Algebra and Combinatorics at $q=0$
   }%
   \\
   \bigskip
   By \\
   \bigskip
   TOM DENTON \\
   \bigskip
   B.S. (University of Oregon) 2006 \\
   \bigskip
   DISSERTATION \\
   \bigskip
   Submitted in partial satisfaction of the requirements for the
   degree of \\
   \bigskip
   DOCTOR OF PHILOSOPHY \\
   \bigskip
   in \\
   \bigskip
   Mathematics \\
   \bigskip
   in the \\
   \bigskip
   OFFICE OF GRADUATE STUDIES \\
   \bigskip        
   of the \\
   \bigskip
   UNIVERSITY OF CALIFORNIA \\
   \bigskip
   DAVIS \\
   \bigskip
   Approved: \\
   \bigskip
   \bigskip
   \makebox[3in]{\hrulefill} \\
   Anne Schilling \\
   \bigskip
   \bigskip
   \makebox[3in]{\hrulefill} \\
   Greg Kuperberg \\
   \bigskip

   \bigskip
   \makebox[3in]{\hrulefill} \\
   Nicolas M. Thi\'ery \\
   \bigskip

   \bigskip
   \makebox[3in]{\hrulefill} \\
   Monica Vazirani \\
   \bigskip
   Committee in Charge \\
   \bigskip
   2011 \\
   \vfill
\end{center}

%% file: Abstract.tex
{\singlespacing
   \begin{flushright}
      Samuel T. Denton \\
      September 2011 \\
      Mathematics \\
   \end{flushright}
}

\bigskip

\begin{center}
   Excursions into Algebra and Combinatorics at $q=0$ \\
\end{center}

\section*{Abstract}

We explore combinatorics associated with the degenerate Hecke algebra at $q=0$, obtaining a formula for a system of orthogonal idempotents, and also exploring various pattern avoidance results.  Generalizing constructions for the $0$-Hecke algebra, we explore the representation theory of $\JJ$-trivial monoids.  

We then discuss two-tensors of crystal bases for $U_q(\tilde{\mathfrak{sl}_2})$, establishing a complementary result to one of Bandlow, Schilling, and Thi\'ery on affine crystals arising from promotion operators.  Finally, we give a computer implementation of Stembridge's local axioms for simply-laced crystal bases.

%% file: Acknowledgments.tex
There are far too many people deserving of acknowledgement to recount here in full, but I will do my best.

First, I should thank Lucia Black for pushing me to apply for graduate school in the midst of a season of eighty-hour work weeks.  There will, after all, be plenty of time to de-lax when we're dead.

(adviser(s))
Numerous people have helped my research through discussion and input.  First and foremost are Prof. Anne Schilling, my adviser, and Prof. Nicolas M. Thi\'ery, my unofficial co-adviser, both of whom suggested interesting problems, helped me out of tight spots, and pushed me to release some of the molehill of code I've developed for my research for use by the greater community.  Prof. Schilling invited me to my first international math conference, the 2008 FPSAC in Chile, introduced me to numerous fascinating mathematical topics, and the equally fascinating mathematicians who drive them forward.  Thanks are also due for a few most excellent nights of band practice.
Prof. Thi\'ery also provided gracious hospitality in Orsay for some four months, the inspiration to finally learn to pass juggling clubs, and forgiveness (I think) for somehow still misspelling his name after all these years.
Florent Hivert was also a pleasure to collaborate and juggle with.  

(Research)
Pat Dragon was a neighbor and office-mate for four years, and often served as an audience for works in progress.  My academic siblings Qiang Wang and Steve Pon provided a great deal of guidance before I knew which way was up, modulo a choice of orientation.  Andrew Berget, Jeff Ferreira, Brant Jones, Katie O'Reilly and Rabhar Virk has also provided interesting conversation over the last couple years.  Sammy Black and David Jordan have now been sharing their beautiful ideas (and couches) for five years.   Thanks are also due to my wonderful teachers at UC Davis and the coordinators of the invaluable VIGRE Research Focus Groups I've attended over the last five years.

(Math staff)
Whether throwing up in class or arranging a last-moment retreat to France, the UC Davis math staff have been invaluable in helping me find practical ways out of difficult situations.  In particular, Celia Davis, Tina Denena and Perry Gee have been inestimably helpful.

(Pedagogy)
Prof. Ken Debevoise taught me valuable lessons in mental endurance (``College students can do anything in the world they set their mind to, so long as it doesn't take more than two weeks.'').  In spite of being a professor of political science, Prof. Debeviose's pedagogy has greatly inspired and informed my approach to teaching mathematics.  I would also like to thank Rick West for his invaluable training in Socratic methods for math instruction.  Ali Dad-Del provided two wonderful summers working with the robotics cluster in the COSMOS program for science-oriented high-school students.  Jesus de Loera, a veritable force of nature, has been inspirational in his combination of research and teaching, providing ways for students at many different levels to contribute and get involved with an active plan of research.  Finally, I should thank my students, who have taught me more than they know.

I have been blessed with a range of great teachers, from Tammy Popp (who disturbed my naps in geometry class and helped me through my first excursions into calculus) to Dmitry Fuchs (who disturbed my naps during algebraic geometry and provided some of the best homework sets I've had the pleasure to work on).  Other stand-outs along this long strange trip have included Dave Robben (who introduced me to Thoreau), Bill Thayer (who taught me to juggle, with a little calculus), Bill Flynn (who taught me more about drawing than I can write in the constrained linear medium of words alone), Martin Guterman (who took time to discuss the Putnam exam on one fateful September morning in 2001), Csaba Szab\'o (who also gave fantastic problem sets), and, again, Ken Debevoise.

(Domies)
Finally, thanks and blame are due in equal parts to the Domes at UC Davis, the utopian paradise I've called home for the last four years.  The Domes have transformed and expanded my notions of what is important and, perhaps more importantly, what is possible.  Many thanks are due to the people who originally created this space, including Ron Swenson and Clay Brandow.  Thanks also to those who have poured blood, sweat, and lost sleep into keeping it running, including (but not limited to) Ben Pearl, Chris Congleton, Jay Erker, JayLee Tuil, Veronica Pardo, Michelle Yates and the Solar Community Housing Association.  

Shout-outs are due to those who participated in the three-year run of the Domes' Gary Gygax Memorial Dungeons and Dragons Game, including Dustin Pluta, Pat Dragon, Ted Tracy, Chris Salam, Ben Miller, Jordan Thompson, Gretchen Kisler, Jay Erker, Brandon Sowers, and all those who sat in for a game (or six) here and there.

Additional individual domies deserving thanks, in no particular order, include, but are not limited to: Matt McCorkle, Kori Farrell, Francesca Claverie, Cat Callaway,  Kurt Vaughn, Kurt Kornbluth, Ina Rommeck, Jonathan Wooley, Jake Lorber, JayLee Tuil, Chuck Parker, Marguerite Wilson, Mike Gordon, Michelle Yates, Shannon Harney, Veronica Pardo, Hrubs, Liz Ernst, Chris Congleton, Isabel Call and probably three dozen others.

Thanks, finally, to Mom, Dad, and the sibs.

%% file: Introduction.tex
Many structures in mathematics have been shown to admit one-parameter deformations, which often allow a more complete understanding of the original object, and occasionally connect various objects that previously seemed quite distinct.  A very beautiful example of this phenomena is given by the $q$-binomial coefficients $\binom{n}{k}_q$, which are polynomials in $q$.  When evaluated at $q=1$, one recovers the usual binomial coefficient $\binom{n}{k}$, which counts the number of subsets of $k$ objects of a set with $n$ objects.  At $q$ a prime power, though, $\binom{n}{k}_q$ counts the number of $k$-dimensional subspaces of the $n$-dimensional vector space over the field with $q$ elements~\cite{stanley97}.  This phenomenon is symptomatic of a much larger interplay between the areas of algebra, combinatorics, and geometry.

Two of the most important examples of $q$-deformations are the \emph{Iwahori-Hecke algebra}, or \emph{Hecke algebra}\footnote{This algebra was first defined by Iwahori, who named it after Erich Hecke.  It is traditional to give credit to Iwahori at the outset and then refer to the object as the Hecke algebra forevermore.} $H_q(W)$, which is a deformation of a Coxeter group $W$, and the \emph{quantum group} $U_q(\mathfrak{g})$, deforming the enveloping algebra of a Lie algebra $\mathfrak{g}$.  In both of these cases, specialization at $q=1$ recovers the original object.  These deformations have been important in establishing canonical bases for representations of the original objects, and have also proved useful in studying representation theory over finite fields~\cite{khazdanLusztig.79, Bjorner_Brenti.2005, hongKang, klimykSchmudgen}.

The majority of this work is concerned with degenerate specializations of $q$-deformations at $q=0$.  In the case of the Hecke algebra, the $0$-Hecke algebra of a Coxeter group $W$ is no longer semi-simple, but still has a great deal of structure:  It is a monoid algebra over a monoid generated by idempotent ``anti-sorting'' operators, with a one-dimensional simple representation for each subset of a collection of simple generators of $W$.  This is an example of an algebra of a $\JJ$-trivial monoid; we will also discuss aspects of the reprsentation theory of such monoids.

We also discuss crystal bases, which arise from representations of a $q$-deformation of the enveloping algebra of a Lie algebra $\mathfrak{g}$.  This deformation $U_q(\mathfrak{g})$ is known as a quantum group, and has been studied extensively; see~\cite{hongKang, klimykSchmudgen} for background.  While one cannot set $q=0$ in this construction (there are unavoidable $q^{-1}$'s in the definition of $U_q(\mathfrak{g})$), Kashiwara demonstrated the existence of a certain lattice $\mathcal{L}$ that exists in $U_q(\mathfrak{g})$-modules, and such that the quotient $\mathcal{L}/q\mathcal{L}$ often has a convenient basis compatible with the structure of the representation~\cite{kashiwara.90}.  This basis can often be lifted to a ``global basis'' independent of $q$, and has been useful in understanding the internal structure of modules for $U_q(\mathfrak{g})$.

\section{Orthogonal Idempotents in the $0$-Hecke Algebra of the Symmetric Group.}

The $0$-Hecke algebra $\mathbb{C}H_0(S_N)$ for the symmetric group $S_N$ can be obtained as the Iwahori-Hecke algebra of the symmetric group $H_q(S_N)$ at $q=0$.  It can also be constructed as the algebra of the monoid generated by anti-sorting operators on permutations of $N$.

P. N. Norton described the full representation theory of $\mathbb{C}H_0(S_N)$ in \cite{Norton.1979}:  In brief, there is a collection of $2^{N-1}$ simple representations indexed by subsets of the usual generating set for the symmetric group, in correspondence with a collection of $2^{N-1}$ projective indecomposable modules.  Norton gave a construction for some elements generating these projective modules, however these elements were neither orthogonal nor idempotent.  While it was known that an orthogonal collection of idempotents to generate the indecomposable modules exists, there was no known formula for these elements.

Herein, we describe an explicit construction for two different families of orthogonal idempotents in $\mathbb{C}H_0(S_N)$, one for each of the two orientations of the Dynkin diagram for $S_N$.  The construction proceeds by creating a collection of $2^{N-1}$ \emph{demipotent} elements, which we call \emph{diagram demipotents}, each indexed by a copy of the Dynkin diagram with signs attached to each node.  These elements are demipotent in the sense that, for each element $X$, there exists some number $k\leq N-1$ such that $X^j$ is idempotent for all $j\geq k$.  The collection of idempotents thus obtained provides a maximal orthogonal decomposition of the identity.

An important feature of the $0$-Hecke algebra is that it is the monoid algebra of a $\mathcal{J}$-trivial monoid.  As a result, its representation theory is highly combinatorial.  This paper is part of an ongoing effort with Hivert, Schilling, and Thi\'ery \cite{dhst.2011} to characterize the representation theory of general $\mathcal{J}$-trivial monoids, continuing the work of \cite{Norton.1979, Carter.1986, Hivert.Thiery.HeckeGroup.2007}.  This effort is part of a general trend to better understand the representation theory of finite semigroups.  See, for example, \cite{Izhakian.2010.SemigroupsRepresentationsOverSemirings, Steinberg.2006.Moebius, Steinberg.2008.MoebiusII, Almeida_Margolis_Steinberg_Volkov.2009, pennell_putcha_renner.1997}, and for a general overview, \cite{Ganyushkin_Mazorchuk_Steinberg.2009}.

The diagram demipotents obey a branching rule which compares well to the situation in \cite{okounkov96} in their ``New Approach to the Representation Theory of the Symmetric Group.''  In their construction, the branching rule for $S_N$ is given primary importance, and yields a canonical basis for the irreducible modules for $S_N$ which pulls back to bases for irreducible modules for $S_{N-M}$.  

Okounkov and Vershik further make extensive use of a maximal commutative algebra generated by the Jucys-Murphy elements.  In the $0$-Hecke algebra, their construction does not directly apply, because the deformation of Jucys-Murphy elements (which span a maximal commutative subalgebra of $\mathbb{C}S_N$) to the $0$-Hecke algebra no longer commute.  Instead, the idempotents obtained from the diagram demipotents play the role of the Jucys-Murphy elements, generating a commutative subalgebra of $\mathbb{C}H_0(S_N)$ and giving a natural decomposition into indecomposable modules, while the branching diagram describes the multiplicities of the irreducible modules.

The Okounkov-Vershik construction is well-known to extend to group algebras of general finite Coxeter groups~\cite{Ram97seminormalrepresentations}.  It remains to be seen whether our construction for orthogonal idempotents generalizes beyond type $A$.  However, the existence of a process for type $A$ gives hope that the Okounkov-Vershik process might extend to more general $0$-Hecke algebras of Coxeter groups.

Following this work, Berg, Bergeron, Bhargava and Saliola described a method for constructing families of orthogonal idempotents for general $\RR$-trivial monoids~\cite{Berg_Bergeron_Bhargava_Saliola.2010}.  Their work provides an interesting middle ground between the fully combinatorial formula in this chapter and the general construction of idempotents from the semi-simple quotient.  The general method for construction of primitive idempotents is described, for example, in~\cite{curtis_reiner.1962}, and in~\cite{dhst.2011}, where very explicit algorithms are provided for $\JJ$-trivial monoids (which are a subset of $\RR$-trivial monoids).  As one might expect, these various constructions become computationally more difficult with greater generality.

The results in this chapter originally appeared in the Electronic Journal of Combinatorics~\cite{Denton.2011}.  

\section{Representation Theory of $\JJ$-Trivial Monoids}

The representation theory of the $0$-Hecke algebra (also called
\emph{degenerate Hecke algebra}) was first studied by P.-N.~Norton~\cite{Norton.1979}
in type A and expanded to other types by
Carter~\cite{Carter.1986}. Using an analogue of Young symmetrizers, they
describe the simple and indecomposable projective modules together with the
Cartan matrix. An interesting combinatorial application was then found by Krob
and Thibon~\cite{Krob_Thibon.NCSF4.1997} who explained how induction and
restriction of these modules gives an interpretation of the products and
coproducts of the Hopf algebras of noncommutative symmetric functions and
quasi-symmetric functions. Two other important steps were further made by
Duchamp--Hivert--Thibon~\cite{Duchamp_Hivert_Thibon.2002} for type $A$ and
Fayers~\cite{Fayers.2005} for other types, using the Frobenius structure to
get more results, including a description of the Ext-quiver. Through divided difference (Demazure
operator), the $0$-Hecke algebra has a central role in Schubert calculus and
also appeared has connection with $K$-theory
\cite{Demazure.1974,Lascoux.2001,Lascoux.2003,Miller.2005,
  Buch_Kresch_Shimozono.2008,Lam_Schilling_Shimozono.2010}.

Like several algebras whose representation theory was studied in
recent years in the algebraic combinatorics community (such as degenerate
left regular bands, Solomon-Tits algebras, ...), the $0$-Hecke
algebra is the algebra of a finite monoid endowed with special
properties. Yet this fact was seldom used,
despite a large body of literature on finite semigroups, including
representation theory
results~\cite{Putcha.1996,Putcha.1998,Saliola.2007,Saliola.2008,Margolis_Steinberg.2008,Schocker.2008,Steinberg.2006.Moebius,Steinberg.2008.MoebiusII,AlmeidaMargolisVolkov05,Almeida_Margolis_Steinberg_Volkov.2009,Ganyushkin_Mazorchuk_Steinberg.2009,Izhakian.2010.SemigroupsRepresentationsOverSemirings}. From these, one can see that much of the representation theory of a
semigroup algebra is combinatorial in nature (provided the
representation theory of groups is known).  One can expect, for
example, that for aperiodic semigroups (which are semigroups which
contain only trivial subgroups) most of the numerical information
(dimensions of the simple/projective indecomposable modules,
induction/restriction constants, Cartan matrix) can be
computed without using any linear algebra.
In a monoid with partial inverses, one finds (non-trivial) local groups and 
an understanding of the representation theory of these groups is necessary 
for the full representation theory of the monoid.  
In this sense, the notion of aperiodic monoids is orthogonal to that 
of groups as they contain only trivial group-like structure (there are no 
elements with partial inverses).  On the same token, their 
representation theory is orthogonal to that of groups.

The class of $\JJ$-trivial monoids is by itself an active subject of
research (see
e.g.~\cite{Straubing.Therien.1985,Henckell_Pin.2000,Vernitski.2008}),
and contains many monoids of interest, starting with the $0$-Hecke
monoid.  Another classical $\JJ$-trivial monoid is that of
nondecreasing parking functions, or monoid of order preserving
regressive functions on a chain.
Hivert and Thi\'ery~\cite{Hivert.Thiery.HeckeSg.2006,Hivert.Thiery.HeckeGroup.2007}
showed that it is a natural quotient of the $0$-Hecke monoid and used
this fact to derive its complete representation theory. It is also a
quotient of Kiselman's monoid which is studied
in~\cite{Kudryavtseva_Mazorchuk.2009} with some representation theory
results. Ganyushkin and Mazorchuk~\cite{Ganyushkin_Mazorchuk.2010}
pursued a similar line with a larger family of quotients of both the
$0$-Hecke monoid and Kiselman's monoid.

Some complications
necessarily arise in the extension of the program to larger classes of monoids, like
$\RR$-trivial or aperiodic monoids, since the simple modules are not necessarily 
one-dimensional in the latter case.
The approach taken there is to
suppress the dependence upon specific properties of orthogonal
idempotents. Following a complementary line, Berg, Bergeron, Bhargava,
and Saliola~\cite{Berg_Bergeron_Bhargava_Saliola.2010} have very
recently provided a construction for a decomposition of the identity
into orthogonal idempotents for the class of $\RR$-trivial monoids.

\section{Non-Decreasing Parking Functions and Pattern Avoidance}

In this chapter, we investigate various connections between the $0$-Hecke monoid and questions of pattern avoidance, and develop tools for approaching pattern avoidance as an algebraic problem.  

Pattern avoidance is a rich and interesting subject which has received much attention since Knuth first connected the notion of $[231]$-avoidance with stack sortability~\cite{knuth.TAOCP1}.  Pattern avoidance has also appeared in the study of smoothness of Schubert varieties~\cite{billeyLakshmibai.2000, Billey98patternavoidance}, the Temperley-Lieb algebra and the computation of Kazhdahn-Lusztig polynomials~\cite{fan.1996, fanGreen.1999}.  There is also an extensive literature on enumeration of permutations avoiding a given pattern; for an introduction, see~\citetalias{bona.permutations}.

While many have studied pattern avoidance for particular patterns, there has been relatively little attention given to the question of pattern avoidance as a general phenomenon.  Similarly, there has been a great deal of combinatorial insight into questions of pattern avoidance, it has been rare to approach pattern avoidance from an algebraic perspective.
In this chapter, we first introduce a method for reinterpreting pattern containment as equivalent to a factorization problem for certain permutation patterns.  We then use these results directly in analysing the fibers certain quotients of the $0$-Hecke monoid.

We begin by introducing the notion of a width system, which, in some cases, allows the factorization of a permutation $x$ containing a pattern $\sigma$ as $x=y\sigma'z$, where $\sigma'$ is a `shift' of $\sigma$, $y$ and $z$ satisfy certain compatibility requirements, and the $\len(x)=\len(y)+\len(\sigma)+\len(z)$.  This factorization generalizes an important result of Billey, Jockusch, and Stanley~\cite{BilleyJockuschStanley.1993}, which states that any permutation $x$ containing a $[321]$-pattern contains a braid; that is, some reduced word for $x$ in the simple transpositions contains a contiguous subword $s_i s_{i+1} s_i$.  (This subword, in our context, plays the role of the $\sigma'$.)  Equivalently, a permutation that is $[321]$-avoiding is fully commutative, meaning that every reduced word may be obtained by commutation relations.  These permutations have been extensively studied, with major contributions by Fan and Green~\cite{fan.1996, fanGreen.1999} and Stembridge~\cite{stembridge.1996}, who associated a certain poset to each fully commutative element, where linear extensions of the poset are in bijection with reduced words for the permutation. 

Width systems allow us to extend this notion of subword containment considerably, and give an algebraic condition for pattern containment for certain patterns.  The width system is simply a measure of various widths of a pattern occurrence within a permutation (called an `instance').  For certain width systems, an instance of minimal width implies a factorization of the form discussed above.  These width systems tend to exist for relatively long permutations.  The main results are contained in Propositions~\ref{prop:s2widthSystems}, \ref{prop:s3widthSystems}, \ref{prop:widthSystemExtend}, \ref{prop:widthSystemExtend2}, and Corollary~\ref{cor:bountifulPerms}.  

We then apply these ideas directly, and study pattern avoidance of certain patterns (most interestingly  $[321]$-avoidance) in the context of quotients of the $0$-Hecke monoid.
Non-decreasing parking functions $\NDPF_N$ may be realized as a quotient of the $0$-Hecke monoid for the symmetric group $S_N$, and coincide with the set of order-preserving regressive functions on a poset when the poset is a chain.  These functions are enumerated by the Catalan numbers; for example, if one represents $f \in \NDPF_N$ as a step function, its graph will be a (rotated) Dyck path.  These functions form a $\JJ$-trivial monoid under composition, and may be realized as a quotient of the $0$-Hecke monoid.  We show that the fibers of this quotient each contain a unique $[321]$-avoiding permutation of minimal length, and a $[231]$-avoiding permutation of maximal length (Theorem ~\ref{thm:ndpfFibers231}).  We then show that a slightly modified quotient has fibers containing a unique $[321]$-avoiding permutation of minimal length, and a $[312]$-avoiding permutation of maximal length (Theorem~\ref{thm:ndpfFibers312}).

This provides a bijection between $[312]$ and $[321]$-avoiding permutations that is very similar in spirit to the bijection of Simion and Schmidt between $[132]$-avoiding permutations and $[123]$-avoiding permutations~\cite{simion.schmidt.1985}.  (The patterns $[312]$ and $[123]$ are the respective ``complements'' of the patterns $[312]$ and $[321]$.)

We then combine these results to obtain a bijection between $[4321]$-avoiding permutations and elements of a submonoid of $\NDPF_{2N}$ (Theorem~\ref{thm:ndpfFibers4321}), which we consider as a parabolic submonoid of a type $B$ generalization of non-decreasing parking functions.

We then expand our discussion to the affine symmetric group and affine $0$-Hecke monoid.  The affine symmetric group was introduced originally by Lusztig~\cite{lusztig.1983}, and questions concerning pattern avoidance in the affine symmetric group have recently been studied by Lam~\cite{Lam06affinestanley}, Green~\cite{Green.2002}, Billey and Crites~\cite{billeyCrites.2011}.  Lam and Green separately showed that an affine permutation contains a $[321]$-pattern if and only if it contains a braid, in the same sense as in the finite case.  

We introduce a definition for affine non-decreasing parking functions $\ANDPF_N$, and demonstrate that this monoid of functions may be obtained as a quotient of the affine symmetric group.  We obtain a combinatorial map from affine permutations to $\ANDPF_N$ and demonstrate that this map coincides with the definition of $\ANDPF_N$ by generators and relations as a quotient of $\ASn_N$.  Finally, we prove that each fiber of this quotient contains a unique $[321]$-avoiding element of minimal length (Theorem~\ref{thm:affNdpfFibers321}).

\section{Some Results on Crystal Bases}

Crystal bases were originally introduced by Kashiwara~\cite{kashiwara.90, kashiwara.91, kashiwara.nakashima.94} to describe the internal structure of representations of a semi-simple Lie algebra $\ggc$, and over time the theory was expanded to include crystals for representations of affine Lie algebras~\cite{kang.kashiwara.misra.miwa.nakashima.nakayashiki.91, kang.kashiwara.misra.miwa.nakashima.nakayashiki.92, kang.kashiwara.misra.94}.  The construction of the crystal first involves a $q$-deformation of the enveloping algebra of the Lie algebra, yielding a quantum group $U_q(\ggc)$.  Then one chooses a representation of the quantum group and a certain lattice within this representation.  Finally, by taking $q$ to zero, this lattice yields a crystal basis, which has the structure of an edge-colored digraph called a crystal graph.  The vertices of the digraph index a basis of the representation, and the edges are colored according to the action of the $\tilde{e}_i$ and $\tilde{f}_i$ operators on this basis.  A crystal basis can often be pulled back to a global basis for $U_q(\ggc)$ at any $q$.  In particular, at $q=1$, one recovers a basis for the representation of $U(\ggc)$.

\subsection{Promotion Operators }

Highest weight representations for $U_q(\ggh)$, where $\ggh$ is an affine Lie algebra, are infinite-dimensional.  As a result, the crystal bases for these representations are infinite.  However, a modification of the weight lattice used in the definition of the quantum group yields an object $U_q'(\ggh)$ which admits finite dimensional representations.  These representations are no longer highest-weight representations, but often (when the crystal is ``perfect'') can be used as building blocks to construct crystals for the infinite-dimensional highest weight representations.   This construction is known as the Kyoto path model~\cite{kang.kashiwara.misra.miwa.nakashima.nakayashiki.91, kang.kashiwara.misra.miwa.nakashima.nakayashiki.92, kang.kashiwara.misra.94}.  The most important of these finite crystals are the Kirillov-Reshetikhin crystals, which have been extensively studied~\cite{fourier.okado.shilling.09}.

By removing the affine node from the Dynkin diagram for an affine Lie algebra, one may restrict back to finite type.  On the level of crystal bases, this process restricts a finite-dimensional $U_q'(\ggh)$ crystal to a classical crystal by removing all edges labeled $0$.  On the other hand, starting with a classical crystal, one may ask whether there exists a way to insert $0$-arrows to obtain a crystal graph for $U_q'(\ggh)$.

For $\ggc=\sln_n$, each crystal basis for a highest weight representation is indexed by a partition $\lambda$, and the crystal graph has vertices labeled by semi-standard Young tableaux.  Applying the crystal operators $\tilde{f_i}$ is a combinatorial operation on a tableau.

We define a \textbf{promotion operator} on a crystal graph $B$ to be a map $\pr: B\rightarrow B$ satisfying the properties:
\begin{itemize}
\item If $\wt(b) = (w_1, \ldots, w_{n+1})$, then $\wt(\pr(b)) = (w_{n+1}, w_1, \ldots, w_n)$,
\item $\pr^{n+1} = \operatorname{id}$, and 
\item For all $i \in \{1, \ldots, n \}$, we have: 
\[
    \pr \circ f_i = f_{i+1}\circ \pr , \text{ and } \pr \circ e_i =e_{i+1} \circ \pr.
\]
\end{itemize}
Given a promotion operator on a classical crystal $B$, one may define an affine structure on $B$ by placing the $0$-arrows according to:
\[
    \pr \circ f_n = f_{0}\circ \pr , \text{ and } \pr \circ e_n =e_{0} \circ \pr.
\]
We call a promotion operator \textbf{connected} if the resulting affine crystal is connected.

When $\lambda$ is a rectangular partition, there exists a combinatorial operation on the crystal called the (canonical) promotion operator $\cpr$ which implements the affine Dynkin diagram automorphism.  Shimozono showed that one may use this promotion operator to insert $0$-arrows into the crystal graph and obtain a finite-dimensional Kirillov-Reshetikhin crystal.  Shimozono further showed that on a classical crystal $B$ of shape $\lambda$, $B$ admits a promotion operator only if $\lambda$ is rectangular, in which case $\cpr$ is the unique promotion operator~\cite{Shimozono.2002}.

Later, Bandlow, Schilling and Thi\'ery showed that on a two-tensor of crystals of rectangular tableaux of type $A_n$ with $n\geq 2$, there exists a unique connected promotion operator, given by the canonical promotion operator acting diagonally on the tensor product~\cite{BST.2010}.

In the case when $n=1$, there exist non-canonical connected promotion operators that give crystals non-isomorphic to the Kirillov-Reshetikhin crystals.  The first goal of this chapter is to show that these non-canonical promotion operators yield crystals that in fact do not arise from representations of $U_q'(\hat{\sln}_2)$.  The main tool here is the classification of representations of $U_q'(\hat{\sln}_2)$ by Chari and Pressley using evaluation representations~\cite{chariPressley.95}.

\subsection{Computer Implementation of Stembridge Local Axioms}

In Section~\ref{sec:stembridge}, we provide a computer implementation of Stembridge's local axioms for crystals arising from highest weight representations in the Sage computer algebra system.  We first review Stembridge's results, then discuss the design of the Sage system, and finally provide code which checks a simply-laced crystal in Sage for compliance with the local axioms.  This base of code could also be extended to check local axioms for non-simply-laced typNes.  The code provided is about to be integrated into the main distribution of Sage.

%% file: chapter-idempotents.tex
Our goal in this chapter is to discuss the salient features of the representation theory of the $0$-Hecke algebra $\CC H_0(W)$, which arises as the algebra over a certain monoid $H_0(W)$ obtained by deformation of a Coxeter group $W$.  We review the construction and features of the algebra, and then give a construction of orthogonal idempotents in the algebra $\CC H_0(S_N)$.  

The results in this chapter originally appeared in the Electronic Journal of Combinatorics~\cite{Denton.2011}, and appeared in an extended abstract for the Formal Power Series and Algebraic Combinatorics conference~\cite{Denton.2010.FPSAC}.

Section~\ref{sec:idempotentsbg} establishes notation and describes the relevant background necessary for the rest of the paper.  For further background information on the properties of the symmetric group, one can refer to the books of \cite{Humphreys.1990} and \cite{stanley97}.  Section~\ref{sec:repTheoryZeroHecke} reviews the essential facts of the representation theory of $H_0(S_N)$.  
Section~\ref{sec:cand} gives the construction of the diagram demipotents.  
Section~\ref{sec:bra} describes the branching rule the diagram demipotents obey, and also establishes the Sibling Rivalry Lemma, which is useful in proving the main results, in Theorem~\ref{thm:main}.  Section~\ref{sec:nilp} establishes bounds on the power to which the diagram demipotents must be raised to obtain an idempotent.  
Finally, remaining questions are discussed in Section~\ref{sec:quest}.

\section{Definitions and Background}
\label{sec:idempotentsbg}

Let $W$ be a finite Coxeter group, which is to say a group generated by reflections.  For $W$ finite, these groups are classified by \emph{Dynkin Diagrams} which encode a system of generators and relations.  Namely, the Dynkin Diagram is a graph on vertex set $I=\{1, \ldots, n\}$, with multiple edges allowed.  Dynkin diagrams also encode data about many other types of objects, such as Hecke algebras and certain Lie algebras.  The number $n=|I|$ is called the \emph{rank} of the Coxeter group (or other object) associated to $D$.  For each pair of indices $(i,j)$ we associate a positive integer $m(i,j)$, equal to two plus the number of edges connecting $i$ and $j$ in the Dynkin diagram.  (Technically, this is the Coxeter diagram, but in our case the two coincide.)  Then $W$ has a generating set $S=\{s_1, \ldots, s_n\}$, satisfying relations:
\begin{equation}
  \begin{alignedat}{2}
    s_i^2 &=1 &\quad& \text{ for all $i\in I$}\,,\\
    \underbrace{s_is_js_is_js_i \cdots}_{m(i,j)} &=
    \underbrace{s_js_is_js_is_j\cdots}_{m(i,j)} && \text{ for all $i,j\in I$}\, ,
  \end{alignedat}
\end{equation}

For example, the Dynkin diagram for the symmetric group $S_{n+1}$ is simply the chain on $n$ vertices.   Then $S_{n+1}$ is generated by a collection of $\{s_1, \ldots, s_n\}$ which can be identified with the simple transpositions (in disjoint cycle notation, given by $(i, i+1)$).  The relations for $S_{n+1}$ are:
\begin{itemize}
\item Reflection: $s_i^2=1$,
\item Commutation: $s_i s_j=s_j s_i$ for $|i-j|>1$,
\item Braid relation: $s_i s_{i+1} s_i=s_{i+1}s_i s_{i+1}$.
\end{itemize}

The Hecke Algebra $\CC H_q(W)$ is a $q$-deformation of the group algebra of $W$, generated by elements $\{T_i \mid i\in I\}$ satisfying relations dependent on a complex parameter $q$:
\begin{equation}
  \begin{alignedat}{2}
    T_i^2 &= ((q-1)T_i + q &\quad& \text{ for all $i\in I$}\,,\\
    \underbrace{T_iT_jT_iT_jT_i \cdots}_{m(i,j)} &=
    \underbrace{T_jT_iT_jT_iT_j\cdots}_{m(i,j)} && \text{ for all $i,j\in I$}\, ,
  \end{alignedat}
\end{equation}
At $q=0$, we obtain the \emph{$0$-Hecke Algebra}, denoted $\CC H_0(W)$.  By making the substitution
$\pi_i:=-T_i$ and considering the monoid generated by the $\pi_i$,
we obtain the \emph{$0$-Hecke monoid}, which we will denote simply by $H$ or $H_0(W)$ if there is any chance of confusion over the originating Coxeter group.  It is clear that the monoid-algebra of $H_0(W)$ is $\CC H_0(W)$.

\begin{remark}[Words for $W$ and $H_0(W)$ Elements.]
Given a list $w=\{w_1,w_2, \ldots, w_k\}$ with $w_i\in I$, and a given collection of generators $\{g_i\}$ indexed by $I$, we can form the word $g_w=g_{w_1}\cdots g_{w_k}$.  For compactness of notation, we will often write words as sequences subscripting the symbol for the generating set.  Thus, $\pi_1\pi_2\pi_3=\pi_{123}$.  (This notation is unambiguous, as we will not explicitly compute any examples of rank greater than nine.)

Elements of the $0$-Hecke monoid are indexed by elements of $W$:  Any reduced word $s=s_{i_1}\cdots s_{i_k}$ for $\sigma\in W$ is also a reduced word in the $0$-Hecke monoid, $\pi_{i_1}\cdots \pi_{i_k}$.  A well-known property of Coxeter groups is that given two reduced words $w$ and $v$ for an element $\sigma$, $w$ is related to $v$ by a sequence of braid and commutation relations~\cite{Bjorner_Brenti.2005}.  These relations still hold in the $0$-Hecke monoid, so $\pi_w=\pi_v$.  From this, we can see that the $0$-Hecke monoid has $|W|$ elements, and that the $0$-Hecke algebra has dimension $|W|$ as a vector space.  

We can obtain a \emph{parabolic subgroup} (resp. submonoid, subalgebra) by considering the object whose generators are indexed by a subset $J\subset I$, retaining the original relations.  Such subgroups will be denoted $W_J$.  The Dynkin diagram of the corresponding object is obtained by deleting the relevant nodes (and incident edges) from the original Dynkin diagram.  It is well known that every Coxeter group (and thus $0$-Hecke monoid) contains a unique longest element, being an element whose length is maximal amongst all elements of the subgroup.  Since the parabolic subgroups (submonoids) are still Coxeter groups, there is a unique longest element in each parabolic subgroup and corresponding submonoid, which we will denote by $s_J \in W$ or $\pi_J^+ \in H_0(W)$.  We will use $\hat{J}$ to denote the complement of $J$ in $I$.  For example, in $H_0(S_8)$ with $J=\{1,2,6\}$, then $w_J^+=\pi_{1216}$, and $w_{\hat{J}}^+=\pi_{3453437}$.

The $0$-Hecke monoid is aperiodic, meaning that for any $x\in H_0(W)$, there exists a finite positive integer $k$ such that $x^k=x^{k+1}$.  In particular, for any element $x\in H_0(W)$ we may define: 
\[
J(x):=\{ i\in I \mid \text{ s.t. $i$ appears in some reduced word for $x$} \}.
\]
This set is well defined because if $i$ appears in some reduced word for $x$ then it must appear in every reduced word for $x$.  Then $x^\omega=w_{J(x)}^+$.  Additonally, this element is, by construction, idempotent.
\end{remark}

\begin{remark}[The Algebra Automorphism $\Psi$ of $\mathbb{C}H_0(S_N)$.]

The algebra $\CC H_0(W)$ is alternatively generated by elements $\bpi{i}:=(1-\pi_i)$, which satisfy the same relations as the $\pi_i$ generators.  There is a unique automorphism $\Psi$ of $\CC H_0(W)$ defined by sending $\pi_i \rightarrow (1-\pi_i)$.

For any longest element $w_J^+$, the image $\Psi(w_J^+)$ is a longest element in the $(1-\pi_i)$ generators; this element is denoted $w_J^-$.
\end{remark}

\begin{remark}[Dynkin Diagram Automorphisms of $\mathbb{C}H_0(W)$.]

Any automorphism of the underlying graph of a Dynkin diagram induces an automorphism of the Hecke algebra.  For the Dynkin diagram of $S_N$, there is exactly one non-trivial automorphism, sending the node $i$ to $N-i+1$.

This diagram automorphism induces an automorphism of the symmetric group, sending the generator $s_i$ to $s_{N-i}$.  Similarly, there is an automorphism of the $0$-Hecke monoid sending the generator $\pi_i$ to $\pi_{N-i}$.
\end{remark}

\section{Representation Theory of $H_0(W)$}
\label{sec:repTheoryZeroHecke}

The representation theory of $\mathbb{C}H_0(S_N)$ was described in \cite{Norton.1979} and expanded to any finite Coxeter groups in \cite{Carter.1986}.  A more general approach to the representation theory can be taken by approaching the $0$-Hecke algebra as a monoid algebra, as per  \cite{Ganyushkin_Mazorchuk_Steinberg.2009}.  The main results are reproduced here for ease of reference.

For any subset $J\subset I$, let $\lambda_J$ denote the one-dimensional representation of $\CC H_0(W)$ defined by the action of the generators:
\begin{equation*}
    \lambda_J(\pi_i) = \begin{cases}
      0 & \text{if $i\in J$},\\
      1 & \text{if $i\notin J.$}
    \end{cases}
\end{equation*}
For $W$ of rank $n$, the $\lambda_J$ are $2^{n}$ non-isomorphic representations, all one-dimensional and thus simple.  In fact, these are all of the simple representations of $\mathbb{C}H_0(W)$, which can be verified by forming a composition series for $H_0(W)$.

\begin{definition}
For each $i\in I$, define the \emph{evaluation maps} $\Phi_i^+$ and $\Phi_i^+$ on generators by:
\begin{eqnarray*}
\Phi_N^+ &:& \mathbb{C}H_0(W) \rightarrow \CC H_0(W_{I\setminus \{i\}}) \\
\Phi_N^+(\pi_i) &=&     \begin{cases}
      1          & \text{if $i=N$,}\\
      \pi_i & \text{if $i \neq N$.}
    \end{cases}\\
\Phi_N^- &:& \mathbb{C}H_0(W) \rightarrow \mathbb{C}H_0(W_{I\setminus \{i\}}) \\
\Phi_N^-(\pi_i) &=&     \begin{cases}
      0          & \text{if $i=N$,}\\
      \pi_i & \text{if $i \neq N$.}
    \end{cases}
\end{eqnarray*}
\end{definition}
One can easily check that these maps extend to algebra morphisms from $H_0(W)\rightarrow H_0(W_{I\setminus i})$.  For any $J$, define $\Phi_J^+$ as the composition of the maps $\Phi_i^+$ for $i\in J$, and define $\Phi_J^-$ analogously.  Then the simple representations of $H_0(W)$ are given by the maps $\lambda_J = \Phi_J^+ \circ \Phi_{\hat{J}}^-$, where $\hat{J}=I\setminus J$.

The map $\Phi_J^+$ is also known as the parabolic map~\cite{Billey_Fan_Lsonczy.1999}, which sends an element $x$ to an element $y$ such that $y$ is the longest element less than $x$ in Bruhat order in the parabolic submonoid with generators indexed by $J$.

The nilpotent radical $\mathcal{N}$ in $\mathbb{C}H_0(S_N)$ is spanned by elements of the form $x-w_{J(x)}^+$, where $x \in H_0(W)$.  This element $w_{J(x)}^+$ is always idempotent.  If $y$ is already idempotent, then $y=w_{J(y)}^+$, and so $y-w_{J(y)}^+=0$ contributes nothing to $\mathcal{N}$.  However, all other elements $x-w_{J(x)}^+$ for $x$ not idempotent are linearly independent, and thus give a basis of $\mathcal{N}$.

Norton further showed that 
\begin{equation*}
\mathbb{C}H_0(S_N)=\bigoplus_{J\subset I} H_0(S_N)w_J^-w_{\hat{J}}^+
\end{equation*}
is a direct sum decomposition of $\mathbb{C}H_0(S_N)$ into indecomposable left ideals, a result which Carter expanded to general Coxeter groups.

\begin{theorem}[Norton, 1979]
\label{thm:norton}
Let $\{p_J | J\subset I\}$ be a family of mutually orthogonal primitive idempotents with $p_J \in \mathbb{C}H_0(S_N)w_J^-w_{\hat{J}}^+$ for all $J\subset I$ such that $\sum_{J\subset I} p_J = 1$.  

Then $\mathbb{C}H_0(S_N)w_J^-w_{\hat{J}}^+=\mathbb{C}H_0(S_N)p_J$, and if $\mathcal{N}$ is the nilpotent radical of $\mathbb{C}H_0(S_N)$, $\mathcal{N}w_J^-w_{\hat{J}}^+=\mathcal{N}p_J$ is the unique maximal left ideal of $\mathbb{C}H_0(S_N)p_J$, and $\mathbb{C}H_0(S_N)p_J/\mathcal{N}p_J$ affords the representation $\lambda_J$.

Finally, the semisimple quotient is commutative and may be described thusly:
\[
\mathbb{C}H_0(S_N)/\mathcal{N} \stackrel{~}{=} \bigoplus_{J\subset I}\mathbb{C}H_0(S_N)p_J/\mathcal{N}p_J 
\stackrel{~}{=} \mathbb{C}^{2^{N-1}}.
\]
\end{theorem}

The elements $w_J^-w_{\hat{J}}^+$ are neither orthogonal nor idempotent; the proof of Norton's theorem is non-constructive, and does not give a formula for the idempotents.

\section{Diagram Demipotents}
\label{sec:cand}

The elements $\pi_i$ and $(1-\pi_i)$ are idempotent.  There are actually $2^{N-1}$ idempotents in $H_0(S_N)$, namely the elements $w_J^+$ for any $J\subset I$.  These idempotents are clearly not orthogonal, though.  The goal of this chapter is to give a formula for a collection of \emph{orthogonal} idempotents in $\mathbb{C}H_0(S_N)$.  

For our purposes, it will be convenient to index subsets of the index set $I$ (and thus also simple and projective representations) by \emph{signed diagrams.}

\begin{definition}
A \emph{signed diagram} is a Dynkin diagram in which each vertex is labeled with a $+$ or $-$.  
\end{definition}

\begin{figure}
\centering
\begin{tikzpicture}[scale=.75]
\node (one) at ( 0,0) [circle,draw] {$1^+$};
\node (two) at ( 2,0) [circle,draw] {$2^+$};
\node (three) at ( 4,0) [circle,draw] {$3^-$};
\node (four) at ( 6,0) [circle,draw] {$4^-$};
\node (five) at ( 8,0) [circle,draw] {$5^-$};
\node (six) at ( 10,0) [circle,draw] {$6^+$};
\node (seven) at ( 12,0) [circle,draw] {$7^-$};

\draw [-] (one.east) -- (two.west);
\draw [-] (two.east) -- (three.west);
\draw [-] (three.east) -- (four.west);
\draw [-] (four.east) -- (five.west);
\draw [-] (five.east) -- (six.west);
\draw [-] (six.east) -- (seven.west);
\end{tikzpicture}
\caption{A signed Dynkin diagram for $S_8$.}
\label{fig:signedDiagram}
\end{figure}
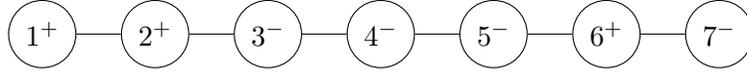

Figure~\ref{fig:signedDiagram} depicts the signed diagram for type $A_7$, corresponding to $H_0(S_8)$ with $J=\{1, 2, 6\}$.  For brevity, a diagram can be written as just a string of signs.  For example, the signed diagram in the Figure is written $++---+-$.  For $k \in I$, and a fixed choice of signed diagram $D$, the $\sgn(k)$ is the sign labeling $k$ in $D$.

We now construct a \emph{diagram demipotent} corresponding to each signed diagram.  Let $P$ be a set partition of the index set $I$ obtained from a signed diagram $D$ by grouping together sets of adjacent pluses and minuses.  For the diagram in Figure~\ref{fig:signedDiagram}, we would have $P = \{ \{1,2\} , \{3,4,5\}, \{6\},\{7\} \}$.  Let $P_k$ denote the $k$th subset in $P$.  For each $P_k$, let $w_{P_k}^{sgn(k)}$ be the longest element of the parabolic sub-monoid associated to the index set $P_k$, constructed with the generators $\pi_i$ if $sgn(k)=+$ and constructed with the $(1-\pi_i)$ generators if $sgn(k)=-$.

\begin{definition}
\label{def:ddemipotents}
Let $D$ be a signed diagram with associated composition $P= P_1 \cup \cdots \cup P_m $.  Set: 
\begin{eqnarray*}
L_D&=&w_{P_1}^{sgn(1)}w_{P_2}^{sgn(2)}\cdots w_{P_m}^{sgn(m)},\text{ and }\\
R_D&=&w_{P_m}^{sgn(m)}w_{P_{m-1}}^{sgn(m-1)}\cdots w_{P_1}^{sgn(1)}.
\end{eqnarray*}
The \emph{diagram demipotent} $C_D$ associated to the signed diagram $D$ is then $L_DR_D$.  The \emph{opposite diagram demipotent} $C_D'$ is $R_DL_D$.
\end{definition}

Thus, the diagram demipotent for the diagram in Figure~\ref{fig:signedDiagram} is 
\[
\pi_{121}^+\pi_{345343}^-\pi_{6}^+\pi_{7}^- \pi_{6}^+\pi_{345343}^-\pi_{121}^+.
\]

It is not immediately obvious that these elements are demipotent; this is a direct result of Lemma~\ref{thm:sibRiv}, below.

For $N=1$, there is only the empty diagram, and the diagram demipotent is just the identity.

For $N=2$, there are two diagrams, $+$ and $-$, and the two diagram demipotents are $\pi_1$ and $1-\pi_1$ respectively.  Notice that these form a decomposition of the identity, as $\pi_i+(1-\pi_i)=1$.

For $N=3$, we have the following list of diagram demipotents.  The first column gives the diagram, the second gives the element written as a product, and the third expands the element as a sum.  For brevity, words in the $\pi_i$ or $\bpi{i}$ generators are written as strings in the subscripts.  Thus, $\pi_1\pi_2$ is abbreviated to $\pi_{12}$.

\begin{equation*}
\begin{array}[b]{|c|c|c|}
D & C_D & \text{$C_D$ Expanded} \\ \hline
++ & \pi_{121} 			& \pi_{121} \\
+- & \pi_1\bpi{2}\pi_1 		& \pi_1 - \pi_{121} \\
-+ & \bpi{1}\pi_2\bpi{1} 	& \pi_2 - \pi_{12} - \pi_{21} + \pi_{121} \\
- -& \bpi{121} 			& 1 - \pi_1 - \pi_2 + \pi_{12} + \pi_{21} - \pi_{121}\\
\end{array}
\end{equation*}

Observations:

\begin{itemize}
\item The idempotent $\bpi{121}$ is an alternating sum over the monoid.  This is a general phenomenon: By \cite{Norton.1979}, $w_J^-$ is the length-alternating signed sum over the elements of the parabolic sub-monoid with generators indexed by $J$.

\item The shortest element in each expanded sum is an idempotent in the monoid with $\pi_i$ generators; this is also a general phenomenon.  The shortest term is just the product of longest elements in nonadjacent parabolic sub-monoids, and is thus idempotent.  Then the shortest term of $C_D$ is $\pi_J^+$, where $J$ is the set of nodes in $D$ marked with a $+$.  Each diagram yields a different leading term, so we can immediately see that the $2^{N-1}$ idempotents in the monoid appear as a leading term for exactly one of the diagram demipotents, and that they are linearly independent.

\item For many purposes, one only needs to explicitly compute half of the list of diagram demipotents; the other half can be obtained via the automorphism $\Psi$.  A given diagram demipotent $x$ is orthogonal to $\Psi(x)$, since one has left and right $\pi_1$ descents, and the other has left and right $\bpi{1}$ descents, and $\pi_1 \bpi{1}=0$.

\item The diagram demipotents are fixed under the automorphism determined by $\pi_{\sigma}\rightarrow \pi_{\sigma^{-1}}$.  In particular, $L_D$ is the reverse of $R_D$, and $C_D$ can be expressed as a palindrome in the alphabet $\{\pi_i, \bpi{i}\}$.

\item The diagram demipotents $C^D$ and $C^E$ for $D\neq E$ do not necessarily commute.  Non-commuting demipotents first arise with $N=6$.  However, the idempotents obtained from the demipotents are orthogonal and do commute.

\item It should also be noted that these demipotents (and the resulting idempotents) are not in the projective modules constructed by Norton, but generate projective modules isomorphic to Norton's.

\item The diagram demipotents $C_D$ listed here are not fixed under the automorphism induced by the Dynkin diagram automorphism.  In particular, the ``opposite'' diagram demipotents $C_D'=R_DL_D$ really are different elements of the algebra, and yield an equally valid but different set of orthogonal idempotents.  For purposes of comparison, the diagram demipotents for the reversed Dynkin diagram are listed below for $N=3$.

\begin{equation*}
\begin{array}[b]{|c|c|c|}
D & C_D' & \text{$C_D'$ Expanded} \\ \hline
++ & \pi_{212} 			& \pi_{212} \\
+- & \pi_2\bpi{1}\pi_2 		& \pi_2 - \pi_{212} \\
-+ & \bpi{2}\pi_1\bpi{2} & \pi_1 - \pi_{12} - \pi_{21} + \pi_{212} \\
- -& \bpi{212} 	& 1 - \pi_1 - \pi_2 + \pi_{12} + \pi_{21} - \pi_{212}\\
\end{array}
\end{equation*}
\end{itemize}

For $N\leq 4$, the diagram demipotents are actually idempotent and orthogonal.  For larger $N$, raising the diagram demipotent to a sufficiently large power yields an idempotent (see below, Theorem~\ref{thm:main}); in other words, the diagram demipotents are indeed demipotent.  The power that an diagram demipotent must be raised to in order to obtain an actual idempotent is called its \emph{nilpotence degree}; we demonstrate below that the nilpotence degree is always $\leq N-3$.

For $N=5$, two of the diagram demipotents need to be squared to obtain an idempotent.  For $N=6$, eight elements must be squared.  For $N=7$, there are four elements that must be cubed, and many others must be squared.  Some pretty good upper bounds on the nilpotence degree of the diagram demipotents are given in Section~\ref{sec:nilp}.  As a preview, for $N>4$ the nilpotence degree is always $\leq N-3$, and conditions on the diagram can often greatly reduce this bound.

As an alternative to raising the demipotent to some power, we can express the idempotent as a product of diagram demipotents for smaller diagrams.  Let $D_k$ be the signed diagram obtained by taking only the first $k$ nodes of $D$.  Then, as we will see, the idempotent can also be expressed as the product $C_{D_1}C_{D_2}C_{D_3}\cdots C_{D_{N-1}}=C_{D}^N$.

\begin{remark}[Right Weak Order.]
Let $m$ be a standard basis element of the $0$-Hecke algebra in the $\pi_i$ basis.  Then for any $i\in D_L(m)$, $\pi_i m=m$, and for any $i\not \in D_L(m)$, $\pi_im\geq_R m$ in left weak order.  This is an adaptation of a standard fact in the theory of Coxeter groups to the $0$-Hecke setting.
\end{remark}

\begin{corollary}[Diagram Demipotent Triangularity]
\label{cor:triangularity}
Let $C_D$ be a diagram demipotent and $m$ an element of the $0$-Hecke monoid in the $\pi_i$ generators.  Then $C_Dm = \lambda m + x$, where $x$ is an element of $H_0(S_N)$ spanned by monoid elements lower in right weak order than $m$, and $\lambda \in \{0,1\}$.  Furthermore, $\lambda=1$ if and only if $D_L(m)$ is exactly the set of nodes in $D$ marked with pluses.
\end{corollary}

\begin{proof}
The diagram demipotent $C_D$ is a product of $\pi_i$'s and $(1-\pi_i)$'s.
\end{proof}

\begin{proposition}
Each diagram demipotent is the sum of a non-zero idempotent part and a nilpotent part.  That is, all eigenvalues of a diagram demipotent are either $1$ or $0$.
\end{proposition}

\begin{proof}
Assign a total ordering to $H_0(S_N)$ as generated by the $\pi_i$ respecting Bruhat order.  Then by Corollary~\ref{cor:triangularity}, the matrix $M_D$ of any diagram demipotent $C_D$ is lower triangular, and each diagonal entry of $M_D$ is either one or zero.  A lower triangular matrix with diagonal entries in $\{0,1\}$ has eigenvalues in $\{0,1\}$; thus $C_D$ is the sum of an idempotent and a nilpotent part.

To show that the idempotent part is non-zero, consider any element $m$ of the monoid such that $D_L(m)$ is exactly the set of nodes in $D$ marked with pluses.  Then $C_Dm = m + x$ shows that $C_D$ has a $1$ on the diagonal, and thus has $1$ as an eigenvalue.  Hence the idempotent part of $C_D$ is non-zero.  (This argument still works if $D$ has no plusses, since the associated diagram demipotent fixes the identity.)
\end{proof}

\section{Branching}
\label{sec:bra}

There is a convenient and useful branching of the diagram demipotents for $H_0(S_N)$ into diagram demipotents for $H_0(S_{N+1})$.

\begin{lemma}
\label{lem:longWords}
Let $J=\{i,i+1,\ldots,N-1\}$
Then $w_J^{+}\pi_N w_J^{+}$ is the longest element in the generators $i$ through $N$.  Likewise, $w_J^{+}\pi_{i-1}w_J^{+}$ is the longest element in the generators $i-1$ through $N-1$.  Similar statements hold for $w_J^{-}\bpi{N} w_J^{-}$ and $w_J^{-}\bpi{i-1} w_J^{-}$.
\end{lemma}

\begin{proof}
Let $J=\{i,i+1,\ldots,N-1\}$.

The lexicographically minimal reduced word for the longest element in consecutive generators $1$ through $k$ is obtained by concatenating the ascending sequences $\pi_{1 \ldots {k-i}}$ for all $0<i<k$.
For example, the longest element in generators $1$ through $4$ is $\pi_{1234123121}$.

Now form the product $m=w_J^+ \pi_N w_J^+$ (for example $\pi_{1234123121}\pi_5\pi_{1234123121}$).  This contains a reduced word for $w_J^+$ as a subword, and is thus $m\geq w_J^+$ in the (strong) Bruhat Order.  But since $w_J^+$ is the longest element in the given generators, $m$ and $w_J^+$ must be equal.

For the second statement, apply the same methods using the lexicographically maximal word for the longest elements.

The analogous statement follows directly by applying the automorphism $\Psi$.
\end{proof}


Recall that each diagram demipotent $C_D$ is the product of two elements $L_D$ and $R_D$.  For a signed diagram $D$, let $D+$ denote the diagram with an extra $+$ adjoined at the end.  Define $D-$ analogously.  

\begin{corollary}
Let $C_D=L_DR_D$ be the diagram demipotent associated to the signed diagram $D$ for $S_N$.  Then $C_{D+}=L_D \pi_N R_D$ and $C_{D-}=L_D \bpi{N} R_D$.  In particular, $C_{D+} + C_{D-} = C_D$.  Finally, the sum of all diagram demipotents for $H_0(S_N)$ is the identity.
\end{corollary}

\begin{proof}
The identities 
\[
C_{D+}=L_D \pi_N R_D \text{ and } C_{D-}=L_D \bpi{N} R_D
\]
are consequences of Lemma~\ref{lem:longWords}, and the identity $C_{D+} + C_{D-} = C_D$ follows directly.

To show that the sum of all diagram demipotents for fixed $N$ is the identity, recall that the diagram demipotent for the empty diagram is the identity, then apply the identity $C_{D+} + C_{D-} = C_D$ repeatedly.
\end{proof}

Next we have a key lemma for proving many of the remaining results in this paper:

\begin{lemma}[Sibling Rivalry]
\label{thm:sibRiv}
Sibling diagram demipotents commute and are orthogonal: $C_{D-}C_{D+}=C_{D+}C_{D-}=0$.  
Equivalently, 
\[
C_{D}C_{D+}=C_{D+}C_{D}=C_{D+}^2 
\text{ and } C_{D}C_{D-}=C_{D-}C_{D}=C_{D-}^2.
\]
\end{lemma}

\begin{proof}

\begin{figure}
\centering
\begin{tikzpicture}[scale=.5]
\node (epsilon) at ( 2,6) [circle,draw] {$r$};

\node (phi) at ( 0,4) [circle,draw] {$q$};

\node (p) at ( -3, 2 )  [circle,draw] {$p$};
\node (bp) at ( 3, 2 )  [circle,draw] {$\barp$};

\node (x) at ( -5, 0 )  [circle,draw] {$x$};
\node (y) at ( -1, 0 )  [circle,draw] {$y$};

\node (bx) at ( 1, 0 )  [circle,draw] {$\barx$};
\node (by) at ( 5, 0 )  [circle,draw] {$\bary$};

\draw [-] (epsilon.west) -- (phi.north) node [above,midway]  {$+$};

\draw [-] (phi.west) -- (p.north) node [above,midway]  {$+$};
\draw [-] (phi.east) -- (bp.north) node  [above,midway] {$-$};

\draw [-] (p.west) -- (x.north) node  [above,midway] {$+$};
\draw [-] (p.east) -- (y.north) node [above,midway] {$-$};

\draw [-] (bp.west) -- (bx.north) node [above,midway] {$+$};
\draw [-] (bp.east) -- (by.north) node [above,midway] {$-$};
\end{tikzpicture}
\caption{Relationship of Elements in the Proof of the Sibling Rivalry Lemma.}
\label{fig:sibRivalry}
\end{figure}
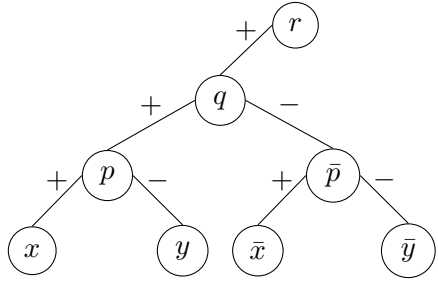

We proceed by induction, using two levels of branching.  Thus, we want to show the orthogonality of two diagram demipotents $x$ and $y$ which are branched from a parent $p$ and grandparent $q$.  Without loss of generality, let $q$ be the positive child of an element $r$.  Call $q$'s other child $\bar{p}$, which in turn has children $\bar{x}$ and $\bar{y}$.  The relations between the elements is summarized in Figure~\ref{fig:sibRivalry}.

The goal, then, is to prove that $yx=0$ and $\bary \barx=0$.  Since $p=x+y$, we have that $yx=(p-x)x=px-x^2$.  Thus, we can equivalently go about proving that $px=x^2$ or $py=y^2$.  It will be easier to show $px=x^2$.  We will also show that $\barp\barx = \barx^2$.  Once this is done, we will have proven the result for diagrams ending in $+++$, $++-$, $+-+$, and $+--$.  By applying the automorphism $\Psi$, we obtain the result for the other four cases.  

One can obtain the reverse equalities $xy=0, \barx\barp=0,$ and so on, either by performing equivalent computations, or else by another use of the $\Psi$ automorphism.  For the latter, suppose that we know $C_{D+}C_{D-}=0$ for arbitrary $D$.  Then applying $\Psi$ to this equation gives $C_{\hat{D}-}C_{\hat{D}+}=0$, where $\hat{D}$ is the signed diagram $D$ with all signs reversed.  Since $D$ was arbitrary, $\hat{D}$ is also arbitrary, so $C_{D-}C_{D+}=0$ for arbitrary $D$.

The remainder of this proof will provide the induction argument.  For the base case, we have $C_{\emptyset}=1$, and $C_{+}=\pi_1$, so clearly $C_{\emptyset}C_+=C_{\emptyset}C_+=C_+=C_+^2$, with analagous statement for $C_{-}$.  For the rank two cases, one can confirm the statement manually using the diagram demipotents listed in Section~\ref{sec:cand}.

Let $r=LR$, dropping the $D$ subscript for convenience, generated with $i$ in the index set $I$.  Let the three new generators be $\pi_a, \pi_b$ and $\pi_c$.  Notice that $\pi_b$, $\bpi{b}$, $\pi_c$, and $\bpi{c}$ all commute with $L$ and $R$.  

The inductive hypothesis tells us that $pq=q p=p^2$ and $\barp q=q \barp=\barp^2$.  We also have the following identities:
\begin{itemize}
\item $q = L \pi_a R$,
\item $p = L \pi_a \pi_b \pi_a R = \pi_b q \pi_b$,
\item $x = L \pi_{aba} \pi_c \pi_{aba} R = \pi_{cbc} q \pi_{cbc}$,
\item $pq = q \pi_b q \pi_b = p^2 = \pi_b q \pi_b q \pi_b$.
\end{itemize}

Then we compute directly:
\begin{eqnarray*}
px & = & \pi_b q \pi_b \pi_{cbc} q \pi_{cbc} \\
   & = & \pi_b q \pi_{cbc} q \pi_{cbc} \\
   & = & \pi_{bc} (q \pi_{b} q \pi_b) \pi_{cbc} \\
   & = & \pi_{bc} ( \pi_b q \pi_{b} q \pi_b) \pi_{cbc} \\
   & = & \pi_{bcb} ( q \pi_{b} q ) \pi_{cbc} \\
   & = & \pi_{cbc} ( q \pi_{cbc} q ) \pi_{cbc} \\
   & = & x^2.
\end{eqnarray*}

To complete the proof, we need to show that $\barp\barx = \barx^2$.  To do so, we use the following identities:
\begin{itemize}
\item $q = L \pi_a R$,
\item $\barp = L \pi_a (1-\pi_b) \pi_a R$,
\item $\barx = L \pi_a (1-\pi_b) \pi_c (1-\pi_b) \pi_a R$.
\end{itemize}

Then we expand the following equation:
\[
\barp \barx = L \pi_a (1-\pi_b) \pi_a R L \pi_a (1-\pi_b) \pi_c (1-\pi_b) \pi_a R.
\]
We expand this as follows:
\[
\barp \barx = q^2 \pi_c 
- q \barp \pi_c 
- q \pi_c \barp 
+ q \pi_c \barp \pi_c 
- \barp q \pi_c  
+ \barp^2 \pi_c  
+ \barp \pi_c \barp
- \barp \pi_c \barp \pi_c.
\]

Meanwhile, 
\begin{eqnarray*}
\barx & = & L ( \pi_{ac} - \pi_{abca} - \pi_{acba} + \pi_{abcba} ) R \\
  & = & \pi_cq - \barp\pi_c - \pi_c \barp + \pi_c \barp \pi_c 
\end{eqnarray*}

Expanding $\barx^2$ in terms of $\barp$ and $q$ is a lengthy but straightforward calculation, which yields:
\begin{eqnarray*}
\barx^2 & = & q^2 \pi_c 
- q \barp \pi_c 
- q \pi_c \barp 
+ q \pi_c \barp \pi_c 
- \barp q \pi_c  
+ \barp^2 \pi_c  
+ \barp \pi_c \barp
- \barp \pi_c \barp \pi_c \\
& = & \barp \barx
\end{eqnarray*}

This completes the proof of the lemma.
\end{proof}


\begin{corollary}
The diagram demipotents $C_D$ are demipotent.
\end{corollary}

This follows immediately by induction: if $C_D^k=C_D^{k+1}$, then $C_{D+}C_D^k=C_{D+}C_D^{k+1}$, and by sibling rivalry, $C_{D+}^{k+1}=C_{D+}^{k+2}$.

Now we can say a bit more about the structure of the diagram demipotents.  

\begin{proposition}
\label{thm:algebra}
Let $p=C_D, x=C_{D+}, y=C_{D-}$, so $p=x+y$ and $xy=0$.  Let $v$ be an element of $H$.  Furthermore, let $p, x$, and $y$ have abstract Jordan decomposition $p=p_i+p_n$, $x=x_i+x_n$, $y=y_i+y_n$, with $p_ip_n=p_np_i$ and $p_i^2=p_i$, $p_n^k=0$ for some $k$, and similar relations for the Jordan decompositions of $x$ and $y$.

Then we have the following relations:
\begin{enumerate}
\item If there exists $k$ such that $p^kv=0$, then $x^{k+1}v=y^{k+1}v=0$.
\item If $pv=v$, then $x(x-1)v=0$
\item If $(x-1)^kv=0$, then $(x-1)v=0$
\item If $pv=v$ and $x^kv=0$ for some $k$, then $yv=v$.
\item If $xv=v$, then $yv=0$ and $pv=v$.
\item Let $u^x_i$ be a basis of the $1$-space of $x$, so that $xu^x_i=u^x_i$, $yu^x_i=0$ and $pu^x_i=v$, and $u^y_j$ a basis of the $1$-space of $y$.  Then the collection $\{u^x_i, u^y_j\}$ is a basis for the $1$-space of $p$.
\item $p_i=x_i+y_i$, $p_n=x_n+y_n$, $x_iy_i=0$.
\end{enumerate}
\end{proposition}

\begin{proof}
\begin{enumerate}
\item Multiply the relation $pv=(x+y)v=0$ by $x$, and recall that $xy=0$.

\item Multiply the relation $pv=(x+y)v=v$ by $x$, and recall that $xy=0$.

\item Multiply $(x-1)^kv=0$ by $y$ to get $yv=0$.  Then $pv=xv$.  Then $(x-1)^kv=(p-1)^kv=0$.  By the induction hypothesis, $(p-1)^kv=(p-1)v$ implies that $pv=v$, but then $xv=pv=v$, so the result holds.

\item By $(2)$, we have $x^2v=xv$, so in fact, $x^kv=xv=0$.  Then $v=pv=xv+yv=yv$.

\item If $xv=v$, then multiplying by $y$ immediately gives $0=yxv=yv$.  Since $yv=0$, then $pv=(x+y)v=xv=v$.

\item From the previous item, it is clear that the bases $v_x^i$ and $v_y^j$ exist with the desired properties.  All that remains to show is that they form a basis for the $1$-space of $p$.  

Suppose $v$ is in the $1$-space of $p$, so $pv=v$.  Then let $xv=a$ and $yv=b$ so that $pv=(x+y)v=a+b=v$.  
Then $a=xv=x(a+b)=x^2v+xyv=x^2v=xa$.  Then $a$ is in the $1$-space of $x$, and, simlarly, $b$ is in the $1$-space of $y$.  Then the $1$-space of $p$ is spanned by the $1$-spaces of $x$ and $y$, as desired.

\item Let $M_p$, $M_x$ and $M_y$ be matrices for the action of $p$, $x$ and $y$ on $H$.  Then the above results imply that the $0$-eigenspace of $p$ is inherited by $x$ and $y$, and that the $1$-eigenspace of $p$ splits between $x$ and $y$.

We can thus find a basis $\{u^x_k, u^y_l, u^0_m \}$ of $H$ such that: $pu^0_k=xu^0_k=yu^0_k=0$, $xu^x_k=u^x_k$, $pu^x_k=u^x_k$, $yu^x_k=0$, $yu^y_k=u^y_k$, $pu^y_k=u^y_k$, and $xu^y_k=0$.  In this basis, $p$ acts as the identity on $\{u^x_k, u^y_l\}$, and $x$ and $y$ act as orthogonal idempotents.  This proves that $p_i=x_i+y_i$ and $x_iy_i=0$.  Since $p=p_i+p_n=x_i+x_n+y_i+y_n$, then it follows that $p_n=x_n+y_n$.
\end{enumerate}
\end{proof}

\begin{corollary}
There exists a linear basis $v_D^j$ of $\mathbb{C}H_0(S_N)$, indexed by a signed diagram $D$ and some numbers $j$, such that the idempotent $I_D$ obtained from the abstract Jordan decomposition of $C_D$ fixes every $v_D^j$.  For every signed diagram $E\neq D$, the idempotent $I_E$ kills $v_D^j$.
\end{corollary}

The proof of this corollary further shows that this basis respects the branching from $H_0(S_{N-1})$ to $H_0(S_{N})$.  In particular, finding this linear basis for $H_0(S_{N})$ allows the easy recovery of the bases for the indecomposable modules for any $M<N$.

\begin{proof}
Any two sibling idempotents have a linear basis for their $1$-spaces as desired, such that the union of these two bases form a basis for their parent's $1$-space.  Then the union of all such bases gives a basis for the $1$-space of the identity element, which is all of $H$.  

All that remains to show is that for every signed diagram $E\neq D$ with a fixed number of nodes, the idempotent $I_E$ kills $v_D^j$.  Let $F$ be last the common ancestor of $D$ and $E$ under the branching of signed diagrams, so that $F+$ is an ancestor of (or equal to) $D$ and $F-$ is an ancestor of (or equal to) $E$.  Then $I_{F+}$ fixes every $v_D^j$, since the collection $v_D^j$ extends to a basis of the $1$-space of $I_{F+}$.  Likewise, $I_{F-}$ kills every $v_D^j$, by the previous theorem.  
\end{proof}

We now state the main result.  For $D$ a signed diagram, let $D_i$ be the signed sub-diagram consisting of the first $i$ entries of $D$.  

\begin{theorem}
\label{thm:main}
Each diagram demipotent $C_D$ (see Definition \ref{def:ddemipotents}) for $H_0(S_N)$ is demipotent, and yields an idempotent $I_D=C_{D_1}C_{D_2}\cdots C_{D}=C_D^N$.  The collection of these idempotents $\{I_D\}$ form an orthogonal set of primitive idempotents that sum to $1$.
\end{theorem}

\begin{proof}
We can completely determine an element of $\mathbb{C}H_0(S_N)$ by examining its natural action on all of $\mathbb{C}H_0(S_N)$, since if $xv=yv$ for all $v\in \mathbb{C}H_0(S_N)$, then $(x-y)v=0$ for every $v$, and $0$ is the only element of $\mathbb{C}H_0(S_N)$ that kills every element of $\mathbb{C}H_0(S_N)$.

The previous results show that the characteristic polynomial of each diagram demipotent is $X^a(X-1)^b$ for some non-negative integers $a$ and $b$, with all nilpotence associated with the $0$-eigenvalue.  This establishes that the diagram demipotents $C_D$ are actually demipotent, in the sense that there exists some $k$ such that $(C_D)^k$ is idempotent.  Theorem~\ref{thm:algebra} shows that this $k$ grows by at most one with each branching, and thus $k\leq N$.  A prior corollary shows that the idempotents sum to the identity.

The previous corollary establishes a basis for $\mathbb{C}H_0(S_N)$ such that each idempotent $I_D$ either kills or fixes each element of the basis, and that for each $E\neq D$, $I_E$ kills the $1$-space of $I_D$.  Since $I_D$ is in the $1$-space of $I_D$, then $I_E$ must also kill $I_D$.  This shows that the idempotents are orthogonal, and completes the theorem.
\end{proof}

\section{Nilpotence Degree of Diagram Demipotents}
\label{sec:nilp}

Take any $m$ in the $0$-Hecke monoid whose descent set is exactly the set of positive nodes in the signed diagram $D$.  Then $C_Dm=m + (\text{lower order terms})$, by a previous lemma, and $I_Dm= (C_D)^k(m)=m + (\text{lower order terms})$.  The set $\{ I_Dm | D_L(m) = \{ \text{positive nodes in D} \}\}$ is thus linearly independent in $H_0(S_N)$, and gives a basis for the projective module corresponding to the idempotent $I_D$.

\begin{figure}
\centering
\begin{tikzpicture}[scale=.4]

\node (empty) at ( 8,18) [circle,draw] {$1$};

\node (p) at ( -1,15) [circle,draw] {$1$};
\node (m) at ( 17,15) [circle,draw] {$\ldots$};

\draw [-] (empty.south) -- (p.north) node [above,midway]  {$+$};
\draw [-] (empty.south) -- (m.north) node [above,midway]  {$-$};

\node (pp) at ( -9,12) [circle,draw] {$1$};
\node (pm) at (  7,12) [circle,draw] {$1$};

\draw [-] (p.south) -- (pp.north) node [above,midway]  {$+$};
\draw [-] (p.south) -- (pm.north) node [above,midway]  {$-$};

\node (ppp) at (-13,9) [circle,draw] {$1$};
\node (ppm) at ( -5,9) [circle,draw] {$1$};
\node (pmp) at (  3,9) [circle,draw] {$1$};
\node (pmm) at ( 11,9) [circle,draw] {$1$};

\draw [-] (pp.south) -- (ppp.north) node [above,midway]  {$+$};
\draw [-] (pp.south) -- (ppm.north) node [above,midway]  {$-$};
\draw [-] (pm.south) -- (pmp.north) node [above,midway]  {$+$};
\draw [-] (pm.south) -- (pmm.north) node [above,midway]  {$-$};

\node (pppp) at (-15,6) [circle,draw] {$1$};
\node (pppm) at (-11,6) [circle,draw] {$1$};
\node (ppmp) at ( -7,6) [circle,draw] {$1$};
\node (ppmm) at ( -3,6) [circle,draw] {$1$};
\node (pmpp) at (  1,6) [circle,draw] {$2$};
\node (pmpm) at (  5,6) [circle,draw] {$2$};
\node (pmmp) at (  9,6) [circle,draw] {$1$};
\node (pmmm) at ( 13,6) [circle,draw] {$1$};

\draw [-] (ppp.south) -- (pppp.north) node [left,midway]  {$+$};
\draw [-] (ppp.south) -- (pppm.north) node [right,midway]  {$-$};
\draw [-] (ppm.south) -- (ppmp.north) node [left,midway]  {$+$};
\draw [-] (ppm.south) -- (ppmm.north) node [right,midway]  {$-$};
\draw [-] (pmp.south) -- (pmpp.north) node [left,midway]  {$+$};
\draw [-] (pmp.south) -- (pmpm.north) node [right,midway]  {$-$};
\draw [-] (pmm.south) -- (pmmp.north) node [left,midway]  {$+$};
\draw [-] (pmm.south) -- (pmmm.north) node [right,midway]  {$-$};

\node (ppppp) at (-16,3) [circle,draw] {$1$};
\node (ppppm) at (-14,3) [circle,draw] {$1$};
\node (pppmp) at (-12,3) [circle,draw] {$1$};
\node (pppmm) at (-10,3) [circle,draw] {$1$};
\node (ppmpp) at ( -8,3) [circle,draw] {$2$};
\node (ppmpm) at ( -6,3) [circle,draw] {$2$};
\node (ppmmp) at ( -4,3) [circle,draw] {$1$};
\node (ppmmm) at ( -2,3) [circle,draw] {$1$};
\node (pmppp) at (  0,3) [circle,draw] {$2$};
\node (pmppm) at (  2,3) [circle,draw] {$2$};
\node (pmpmp) at (  4,3) [circle,draw] {$2$};
\node (pmpmm) at (  6,3) [circle,draw] {$2$};
\node (pmmpp) at (  8,3) [circle,draw] {$2$};
\node (pmmpm) at ( 10,3) [circle,draw] {$2$};
\node (pmmmp) at ( 12,3) [circle,draw] {$1$};
\node (pmmmm) at ( 14,3) [circle,draw] {$1$};

\draw [-] (pppp.south) -- (ppppp.north) node [left,midway]  {$+$};
\draw [-] (pppm.south) -- (pppmm.north) node [right,midway]  {$-$};
\draw [-] (ppmp.south) -- (ppmpp.north) node [left,midway]  {$+$};
\draw [-] (ppmm.south) -- (ppmmm.north) node [right,midway]  {$-$};
\draw [-] (pmpp.south) -- (pmppp.north) node [left,midway]  {$+$};
\draw [-] (pmpm.south) -- (pmpmm.north) node [right,midway]  {$-$};
\draw [-] (pmmp.south) -- (pmmpp.north) node [left,midway]  {$+$};
\draw [-] (pmmm.south) -- (pmmmm.north) node [right,midway]  {$-$};

\draw [-] (pppm.south) -- (pppmp.north) node [left,midway]  {$+$};
\draw [-] (pppp.south) -- (ppppm.north) node [right,midway]  {$-$};
\draw [-] (ppmm.south) -- (ppmmp.north) node [left,midway]  {$+$};
\draw [-] (ppmp.south) -- (ppmpm.north) node [right,midway]  {$-$};
\draw [-] (pmpm.south) -- (pmpmp.north) node [left,midway]  {$+$};
\draw [-] (pmpp.south) -- (pmppm.north) node [right,midway]  {$-$};
\draw [-] (pmmm.south) -- (pmmmp.north) node [left,midway]  {$+$};
\draw [-] (pmmp.south) -- (pmmpm.north) node [right,midway]  {$-$};
\node (pppppx) at (-16,0) [circle,draw] {$1$};
\node (ppppmx) at (-14,0) [circle,draw] {$1$};
\node (pppmpx) at (-12,0) [circle,draw] {$2$};
\node (pppmmx) at (-10,0) [circle,draw] {$1$};
\node (ppmppx) at ( -8,0) [circle,draw] {$3$};
\node (ppmpmx) at ( -6,0) [circle,draw] {$2$};
\node (ppmmpx) at ( -4,0) [circle,draw] {$2$};
\node (ppmmmx) at ( -2,0) [circle,draw] {$1$};
\node (pmpppx) at (  0,0) [circle,draw] {$2$};
\node (pmppmx) at (  2,0) [circle,draw] {$2$};
\node (pmpmpx) at (  4,0) [circle,draw] {$3$};
\node (pmpmmx) at (  6,0) [circle,draw] {$2$};
\node (pmmppx) at (  8,0) [circle,draw] {$2$};
\node (pmmpmx) at ( 10,0) [circle,draw] {$2$};
\node (pmmmpx) at ( 12,0) [circle,draw] {$2$};
\node (pmmmmx) at ( 14,0) [circle,draw] {$1$};

\draw [-] (ppppp.south) -- (pppppx.north) node [left,midway]  {$\pm$};
\draw [-] (ppppm.south) -- (ppppmx.north) node [left,midway]  {$\pm$};
\draw [-] (pppmp.south) -- (pppmpx.north) node [left,midway]  {$\pm$};
\draw [-] (pppmm.south) -- (pppmmx.north) node [left,midway]  {$\pm$};
\draw [-] (ppmpp.south) -- (ppmppx.north) node [left,midway]  {$\pm$};
\draw [-] (ppmpm.south) -- (ppmpmx.north) node [left,midway]  {$\pm$};
\draw [-] (ppmmp.south) -- (ppmmpx.north) node [left,midway]  {$\pm$};
\draw [-] (ppmmm.south) -- (ppmmmx.north) node [left,midway]  {$\pm$};
\draw [-] (pmppp.south) -- (pmpppx.north) node [left,midway]  {$\pm$};
\draw [-] (pmppm.south) -- (pmppmx.north) node [left,midway]  {$\pm$};
\draw [-] (pmpmp.south) -- (pmpmpx.north) node [left,midway]  {$\pm$};
\draw [-] (pmpmm.south) -- (pmpmmx.north) node [left,midway]  {$\pm$};
\draw [-] (pmmpp.south) -- (pmmppx.north) node [left,midway]  {$\pm$};
\draw [-] (pmmpm.south) -- (pmmpmx.north) node [left,midway]  {$\pm$};
\draw [-] (pmmmp.south) -- (pmmmpx.north) node [left,midway]  {$\pm$};
\draw [-] (pmmmm.south) -- (pmmmmx.north) node [left,midway]  {$\pm$};
\end{tikzpicture}
\caption{Nilpotence degree of diagram demipotents.  The root node denotes the diagram demipotent with empty diagram (the identity).  In all computed example, sibling diagram demipotents have the same nilpotence degree; the lowest row has been abbreviated accordingly for readability.}
\label{fig:nilpotence}
\end{figure}
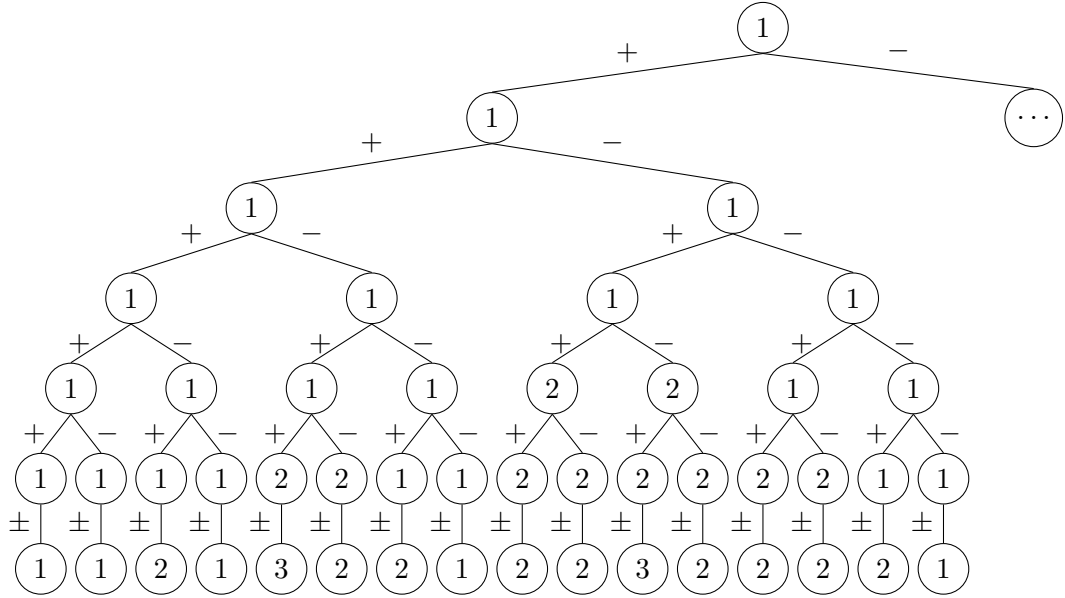

We have shown that for any diagram demipotent $C_D$, there exists a minimal integer $k$ such that $(C_D)^k$ is idempotent.  Call $k$ the \emph{nilpotence degree} of $C_D$.  The nilpotence degree of all diagram demipotents for $N\leq 7$ is summarized in Figure~\ref{fig:nilpotence}.

The diagram demipotent $C_{+\cdots +}$ with all nodes positive is given by the longest word in the $0$-Hecke monoid, and is thus already idempotent.  The same is true of the diagram demipotent $C_{-\cdots -}$ with all nodes negative.  As such, both of these elements have nilpotence degree $1$.

\begin{lemma}
The nilpotence degree of sibling diagram demipotents $C_{D+}$ and $C_{D-}$ are either equal to or one greater than the nilpotence degree $k$ of the parent $C_{D}$.  Furthermore, the nilpotence degree of sibling diagram demipotents are equal.
\end{lemma}

\begin{proof}
Let $x$ and $y$ be the sibling diagram demipotents, with parent diagram demipotent $p$, so $p=C_{D}=L_DR_D, x=C_{D+}=L_D\pi_NR_D, y=C_{D-}=L_D(1-\pi_N)R_D$.  Let $p$ have nilpotence degree $k$, so that $p^k=p^{k+1}$.  We have already seen that the nilpotence degree of $x$ and $y$ is at most $k+1$.  We first show that the nilpotence degree of $x$ or $y$ cannot be less than the nilpotence degree of $p$.  

Recall the following quotients of $\mathbb{C}H_0(S_N)$: 
\begin{eqnarray*}
\Phi_N^+ &:& \mathbb{C}H_0(S_N) \rightarrow \mathbb{C}H_0(S_{N-1}) \\
\Phi_N^+(\pi_i) &=&     \begin{cases}
      1          & \text{if $i=N$,}\\
      \pi_i & \text{if $i \neq N$.}
    \end{cases}\\
\Phi_N^- &:& \mathbb{C}H_0(S_N) \rightarrow \mathbb{C}H_0(S_{N-1}) \\
\Phi_N^-(\pi_i) &=&     \begin{cases}
      0          & \text{if $i=N$,}\\
      \pi_i & \text{if $i \neq N$.}
    \end{cases}
\end{eqnarray*}
given by introducing the relation $\pi_N=1$.  One can easily check that these are both morphisms of algebras.  
Notice that $\Phi_N^+(x)=p$, and $\Phi_N^-(y)=p$.  Then if the nilpotence degree of $x$ is $l<k$, we have 
$p^l=\Phi_N^+(x^l)=\Phi_N^+(x^{l+1})=p^{l+1}$, implying that the nilpotence degree of $p$ was actually $l$, a contradiction.  The same argument can be applied to $y$ using the quotient $\Phi_n^-$.

Suppose one of $x$ and $y$ has nilpotence degree $k$. Assume it is $x$ without loss of generality.  Then:
\begin{eqnarray*}
p^k & = & p^{k+1} \\
\Leftrightarrow x^k+y^k & = & x^{k+1}+y^{k+1} \\
\Leftrightarrow x^{k+1}+y^k & = & x^{k+1}+y^{k+1} \\
\Leftrightarrow y^k & = & y^{k+1} \\
\end{eqnarray*}  
Then the nilpotence degree of $y$ is also $k$.

Finally, if neither $x$ nor $y$ have nilpotence degree $k$, then they both must have nilpotence degree $k+1$.
\end{proof}

Computer exploration suggests that siblings always have equal nilpotence degree, and that nilpotence degree either stays the same or increases by one after each branching.

\begin{lemma}
Let $D$ be a signed diagram with a single sign change, or the sibling of such a diagram.  Then $C_D$ is idempotent (and thus has nilpotence degree $1$).
\end{lemma}

\begin{proof}
We prove the statement for a diagram with single sign change, since siblings automatically have the same nilpotence degree.
Without loss of generality let the diagram of $D$ be 
$--\cdots--++\cdots++$.  Let $L$ the subset of the index set with negative marks in $D$.  Let $i$ be the minimal element of the index set with a positive mark, and let $H=I \setminus (L\cup \{i\})$.  Then:
\[
C_D= w_L^- w_H^+ \pi_i w_H^+ w_L^-.
\]
Notice that $w_H^+$ and $w_L^-$ commute.

Set $y=w_L^- w_H^+ (1-\pi_i) w_H^+ w_L^-$, 
and $p = C_D+y = w_L^- w_H^+ w_H^+ w_L^- = w_H^+ w_L^-$.

Now $y$ is not a diagram demipotent, though $p$ could be considered a diagram demipotent for disconnected Dynkin Diagram with the $i$th node removed.

It is immediate that:
\[
p^2=p, \qquad C_Dp=C_D=pC_D \qquad yp=y=py
\]

Now we can establish orthogonality of $C_D$ and $y$:
\begin{eqnarray*}
C_Dy & = & (w_L^- w_H^+ \pi_i w_H^+ w_L^-) (w_L^- w_H^+ (1-\pi_i) w_H^+ w_L^-) \\
& = & w_L^-( w_H^+ \pi_i w_H^+)(w_L^- (1-\pi_i)w_L^-) w_H^+\\
& = & w_L^- \pi_{H\cup i}^+ \pi_{L\cup i}^- w_H^+ \\
& = & 0
\end{eqnarray*}
The product of $\pi_{H\cup i}^+$ and $\pi_{L\cup i}^-$ is zero, since $\pi_{H\cup i}^+$ has a $\pi_i$ descent, and 
$\pi_{L\cup i}^-$ has a $\barp_i$ descent.

Then $C_D = pC_D = (C_D + y)C_D = (C_D)^2$, so we see that $C_D$ is idempotent.
\end{proof}

In particular, this lemma is enough to see why there is no nilpotence before $N=5$; every signed Dynkin diagrams with three or fewer nodes has no sign change, one sign change, or is the sibling of a diagram with one sign change.

\begin{proposition}
Let $D$ be any signed diagram with $n$ nodes, and let $E$ be the largest prefix diagram such that $E$ has a single sign change, or is the sibling of a diagram with a single sign change.  Then if $E$ has $k$ nodes, the nilpotence degree of $D$ is at most $n-k$.
\end{proposition}

\begin{proof}
This result follows directly from the previous lemma and the fact that the nilpotence degree can increase by at most one with each branching.
\end{proof}

This bound is not quite sharp for $H_0(S_N)$ with $N\leq 7$:  The diagrams $+-++$, $+-+++$, and $+-++++$ all have nilpotence degree $2$.  However, at $N=7$, the highest expected nilpotence degree is $3$ (since every diagram demipotent with three or fewer nodes is idempotent), and this degree is attained by $4$ of the demipotents.  These diagram demipotents are $++-+++$, $+-+-++$, and their siblings.

An open problem is to find a formula for the nilpotence degree directly in terms of the diagram of a demipotent.

\section{Further Directions}
\label{sec:quest}

\subsection{Conjectural Demipotents with Simpler Expression}

Computer exploration has suggested a collection of demipotents that are simpler to describe than those we have presented here.  

For a word $w=(w_1w_2\cdots w_k)$ with $w_i$ in the index set and a signed diagram $D$, we obtain the \emph{masked word} $w^D$ by applying the sign of $i$ in $D$ to each instance of $i$ in $w$.  For example, for the word $w=(1,2,1,3,1,2)$ and $D=+-+$, the masked word is $w^D=(1,-2,1,3,1,-2)$.  A masked word yields an element of $H_0(S_N)$ in the obvious way: we write 
\[
\pi_w^D:=\prod \pi_{w_i}^{sgn(i)},
\]
where $sgn(i)$ is the sign of $i$ in $D$.  

Some masked words are demipotent and others are not.  We call a word \emph{universal} if:
\begin{itemize}
\item $w$ contains every letter in $I$ at least once, and
\item $w^D$ is demipotent for every signed diagram $D$.
\end{itemize} 

\begin{conjecture}
\label{conj:universality}
The word $u_N=(1,2,\ldots, N-2, N-1, N-2, \ldots, 2,1)$ is universal.
\end{conjecture}

Computer exploration has shown that $u_N$ are universal up to $\mathbb{C}H_0(S_9)$, and that the idempotents thus obtained are the same as the idempotents obtained from the diagram demipotents $C_D$.  However, these demipotents $u^D_N$, though they branch in the same way as the diagram demipotents, fail to have the sibling rivalry property.  Thus, another method should be found to show that these elements are demipotent.

An important quotient of the $0$-Hecke monoid is the monoid of \emph{Non-Decreasing Parking Functions}, $\NDPF_N$.  These are the functions $f: [N] \rightarrow [N]$ satisfying
\begin{itemize}
\item $f(i)\leq i$, and 
\item For any $i\leq j$, then $f(i)\leq f(j)$.
\end{itemize}
This monoid can be obtained from $H_0(S_N)$ by introducing the additional relation:
\[
\pi_i\pi_{i+1}\pi_i = \pi_i\pi_{i+1}.
\]
The lattice of idempotents of the monoid $\NDPF_N$ is identical to the lattice of idempotents in $H_0(S_N)$.  We have shown that every masked word $u_N^D$ is idempotent in the algebra of $\NDPF_N$, supporting Conjecture~\ref{conj:universality}.  For the full exploration of $\NDPF_N$, including the proof of the claim that $u_N^D$ is idempotent in $\mathbb{C}\NDPF_N$, see \ref{ch:jtrivial}.

\subsection{Direct Description of the Idempotents}

A number of questions remain concerning the idempotents we have constructed.

First, uniqueness of the idempotents described in this paper is unknown.  In fact, there are many families of orthogonal idempotents in $H_0(S_N)$.  The idempotents we have constructed are invariant as a set under the automorphism $\Psi$, and compatible with the branching from $S_{N-1}$ to $S_N$ according to the choice of orientation of the Dynkin diagram.  

Second, computer exploration has shown that, over the complex numbers, the idempotents obtained from the diagram demipotents have $\pm 1$ coefficients.  This phenomenon has been observed up to $N=9$.  This seems to be peculiar to the construction we have presented, as we have found other idempotents that do not have this property.  It would be interesting to have an even more direct construction of the idempotents, such as a rule for directly determining the coefficients of each idempotent.

It should be noted that a general `lifting' construction has long been known, which constructs orthogonal idempotents in the algebra.  (See~\cite[Chapter 77]{curtis_reiner.1962})  A particular implementation of this lifting construction for algebras of $\mathcal{J}$-trivial monoids is given in~\cite{dhst.2011}.  This lifting construction starts with the idempotents in the monoid, which in the semisimple quotient have the multiplicative structure of a lattice.  In the case of a $0$-Hecke algebra with index set $I$, these idempotents are just the long elements $w_J^+$, for any $J\subset I$.  Then the multiplication rule in the semisimple quotient for two such idempotents $w_J^+$, $w_K^+$ is just $w_K^+w_J^+=w_{J\cup K}^+$.  Each idempotent in the semisimple quotient is in turn lifted to an idempotent in the algebra, and forced to be orthogonal to all idempotents previously lifted.  Many sets of orthogonal idempotents can be thus obtained, but the process affords little understanding of the combinatorics of the underlying monoid.  

The $\pm 1$ coefficients that have been observed in the idempotents thus far constructed suggest that there are still interesting combinatorics to be learned from this problem.

\subsection{Generalization to Other Types}

A combinatorial construction for idempotents in the $0$-Hecke algebra for general Coxeter groups would be desirable.  It is simple to construct idempotents for any rank 2 Dynkin diagram.  The author has also constructed idempotents for type $B_3$ and $D_4$, but has not been able to find a satisfactory formula for general type $B_N$ or $D_N$.  

A major obstruction to the direct application of our construction to other types arises from our expressions for the longest elements in type $A_N$.  For the index set $J \cup \{k\} \subset I$, where $k$ is larger (or smaller) than any index in $J$ we have expressed the longest element for $J\cup \{\pi_k\}$ as $w_J^+\pi_k w_J^+$.  This expression contains only a single $\pi_k$.  In every other type, expressions for the longest element generally require at least two of any generator corresponding to a leaf of the Dynkin diagram.  This creates an obstruction to branching demipotents in the way we have described for type $A_N$.

For example, in type $D_4$, a reduced expression for the longest element is $\pi_{423124123121}$.  The generators corresponding to leaves in the Dynkin diagram are $\pi_1$, $\pi_3$, and $\pi_4$, all of which appear at least twice in this expression.  (In fact, this is true for any of the $2316$ reduced words for the longest element in $D_4$.)  Ideally, to branch easily from type $A_3$, we would be able to write the long element in the form $w_J^+\pi_4w_J^+$, where $4 \not \in J$, but this is clearly not possible.

%% file: chapter-jtrivial.tex
We describe the general representation theory of $\JJ$-trivial monoids, which includes the zero-Hecke monoids.  We analyze specific examples of $\JJ$-trivial monoids, including semi-lattices, the monoid of order-preserving functions on a poset, and non-decreasing parking functions.  The non-decreasing parking functions may be obtained as a quotient of the zero-Hecke monoid; using this fact, we obtain a formula for orthogonal idempotents and applications to pattern avoidance.  We also conjecture an algorithm for obtaining a family of orthogonal idempotents in the algebra of order-preserving functions on a poset.

The results in this chapter originally appeared in the S\'eminaire Lotharingien de Combinatoire~\cite{dhst.2011}.  Section~\ref{sec:jrep} is an abridgement of the version that appears in~\cite{dhst.2011}, but the other sections are identical.

The chapter is arranged as follows. In Section~\ref{sec:bgnot} we recall
the definition of a number of classes of monoids, including the
$\JJ$-trivial monoids, define some running examples of $\JJ$-trivial
monoids, and establish notation.

In Section~\ref{sec:jrep} we recount a few results on the
representation theory of $\JJ$-trivial monoids (a full version of this section, 
with proofs and more results, may be found in~\cite{dhst.2011}), and illustrate them
in the context of the $0$-Hecke monoid.
All the constructions and proofs involve only combinatorics in the
monoid. Due to this, the
results do not depend on the ground field $\K$. In fact, we have
checked that all the arguments pass to $\K=\ZZ$ and therefore to any
ring. It sounds likely that the theory would apply
mutatis-mutandis to semi-rings, in the spirit
of~\cite{Izhakian.2010.SemigroupsRepresentationsOverSemirings}.

Finally, in Section~\ref{sec:NPDF}, we examine the monoid of order
preserving regressive functions on a poset $P$, which generalizes the
monoid of nondecreasing parking functions on the set $\{1, \ldots,
N\}$. We give combinatorial constructions for idempotents in the
monoid and also prove that the Cartan matrix is upper triangular. In
the case where $P$ is a meet semi-lattice (or, in particular, a
lattice), we establish an idempotent generating set for the monoid,
and present a conjectural recursive formula for orthogonal idempotents
in the algebra.

\section{Background and Notation}
\label{sec:bgnot}

A \emph{monoid} is a set $M$ together with a binary operation $\cdot :
M\times M \to M$ such that we have \emph{closure} ($x \cdot y\in M$
for all $x,y \in M$), \emph{associativity} ( $(x\cdot y) \cdot z = x
\cdot ( y \cdot z)$ for all $x,y,z \in M$), and the existence of an
\emph{identity} element $1\in M$ (which satistfies $1\cdot x = x \cdot
1 = x$ for all $x\in M$). 
In this paper, unless explicitly mentioned, all monoids are \emph{finite}.
We use the convention that $A\subseteq B$ denotes $A$ a subset of $B$, and
$A\subset B$ denotes $A$ a proper subset of $B$.

\begin{figure}
  \includegraphics[scale=0.75]{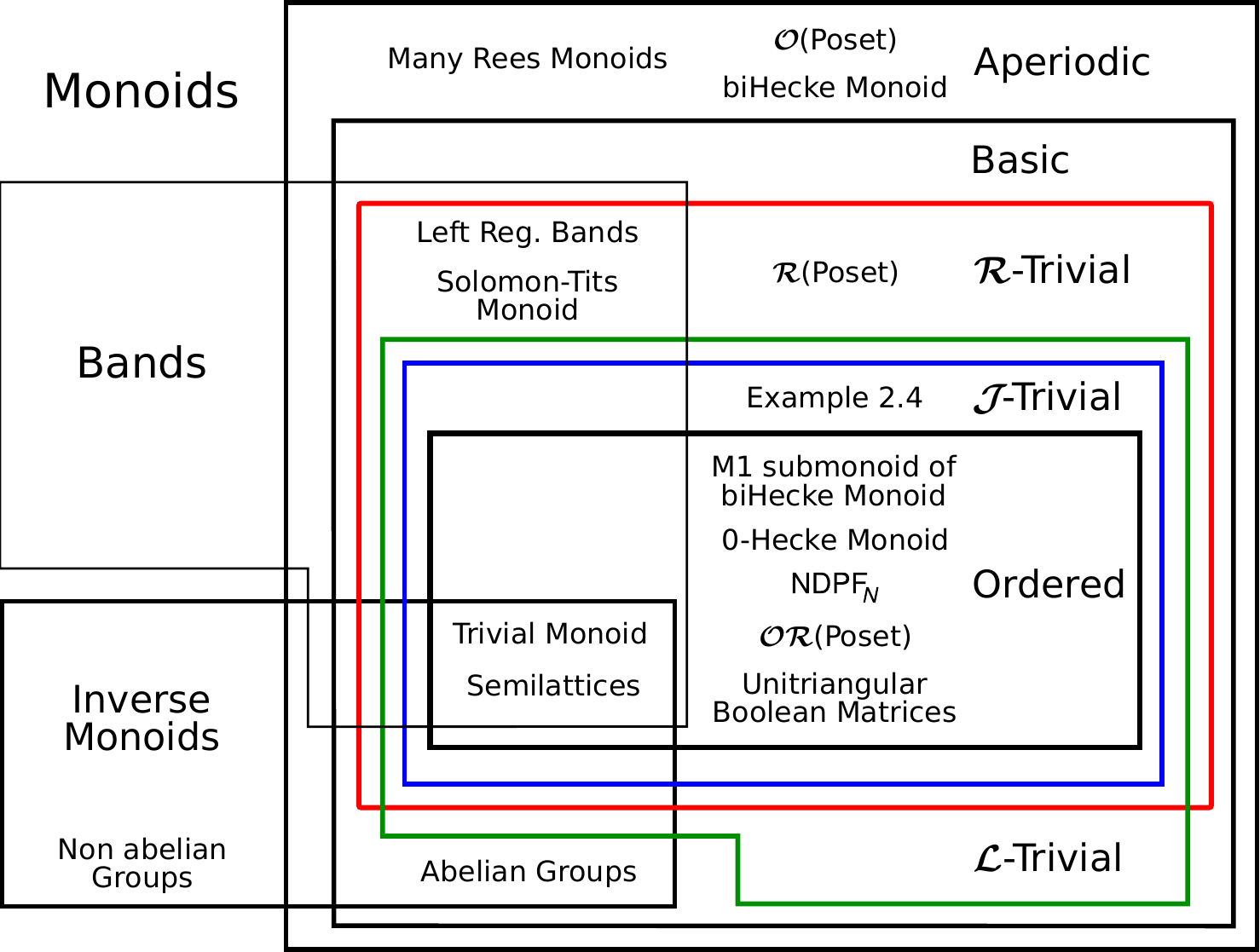}
  \caption{Classes of finite monoids, with examples}
  \label{fig.monoids}
\end{figure}

Monoids come with a far richer diversity of features than groups, but
collections of monoids can often be described as \emph{varieties}
satisfying a collection of algebraic identities and closed under
subquotients and finite products (see e.g.~\cite{Pin.1986,Pin.2009}
or~\cite[Chapter VII]{Pin.2009}).  Groups are an example of a variety
of monoids, as are all of the classes of monoids described in this
paper.  In this section, we recall the basic tools for monoids, and
describe in more detail some of the varieties of monoids that are
relevant to this paper. A summary of those is given in
Figure~\ref{fig.monoids}.

In 1951 Green introduced several preorders on monoids which are essential for the study of
their structures (see for example~\cite[Chapter V]{Pin.2009}). Let $M$ be a monoid and 
define $\le_\RR, \le_\LL, \le_\JJ, \le_\HH$ for $x,y\in M$ as follows:
\begin{equation*}
\begin{split}
	&x \le_\RR y \quad \text{if and only if $x=yu$ for some $u\in M$}\\
	&x \le_\LL y \quad \text{if and only if $x=uy$ for some $u\in M$}\\
	&x \le_\JJ y \quad \text{if and only if $x=uyv$ for some $u,v\in M$}\\
	&x \le_\HH y \quad \text{if and only if $x\le_\RR y$ and $x\le_\LL y$.}
\end{split}
\end{equation*}
These preorders give rise to equivalence relations:
\begin{equation*}
\begin{split}
	&x \; \RR \; y \quad \text{if and only if $xM = yM$}\\
	&x \; \LL \; y \quad \text{if and only if $Mx = My$}\\
	&x \; \JJ \; y \quad \text{if and only if $MxM = MyM$}\\
	&x \; \HH \; y \quad \text{if and only if $x \; \RR \; y$ and $x \; \LL \; y$.}
\end{split}
\end{equation*}

We further add the relation $\le_\BB$ (and its associated equivalence
relation $\BB$) defined as the finest preorder such that $x\leq_\BB
1$, and
\begin{equation}
\label{equation.bruhat}
  \text{$x\le_\BB y$ implies that $uxv\le_\BB uyv$ for all $x,y,u,v\in M$.}
\end{equation}
(One can view $\le_\BB$ as the intersection of all preorders with the above
property; there exists at least one such preorder, namely $x\le y$ for all $x,y\in M$).

Beware that $1$ is the largest element of these (pre)-orders. This is
the usual convention in the semi-group community, but is the converse
convention from the closely related notions of left/right/Bruhat order
in Coxeter groups.

\begin{definition}
A monoid $M$ is called $\KK$-\emph{trivial} if all $\KK$-classes are
of cardinality one, where $\KK\in \{\RR,\LL,\JJ,\HH,\BB\}$.
\end{definition}

An equivalent formulation of $\KK$-triviality is given in terms of
\emph{ordered} monoids. A monoid $M$ is called:
\begin{equation*}
\begin{aligned}
  & \text{\emph{right ordered}} && \text{if $xy\le x$ for all $x,y\in M$}\\
  & \text{\emph{left ordered}} && \text{if $xy\le y$ for all $x,y\in M$}\\
  & \text{\emph{left-right ordered}} && \text{if $xy\leq x$ and $xy\leq y$ for all $x,y\in M$}\\
  & \text{\emph{two-sided ordered}} && \text{if $xy=yz \leq y$ for all $x,y,z\in M$ with $xy=yz$}\\
  & \text{\emph{ordered with $1$ on top}} && \text{if $x\leq 1$ for all $x\in M$, and $x\le y$}\\ 
  & && \text{implies $uxv\le uyv$ for all $x,y,u,v\in M$}
\end{aligned}
\end{equation*}
for some partial order $\le$ on $M$.

\begin{proposition}
  \label{proposition.ordered}
  $M$ is right ordered (resp. left ordered, left-right ordered, two-sided ordered,
  ordered with $1$ on top) if and only if $M$ is $\RR$-trivial
  (resp. $\LL$-trivial, $\JJ$-trivial, $\HH$-trivial, $\BB$-trivial).

  When $M$ is $\KK$-trivial for $\KK\in \{\RR,\LL,\JJ,\HH,\BB\}$, then
  $\leq_\KK$ is a partial order, \emph{called
    $\KK$-order}. Furthermore, the partial order $\le$ is finer than
  $\le_\KK$: for any $x, y\in M$, $x \le_\KK y$ implies $x \leq y$.
\end{proposition}
\begin{proof}
  We give the proof for right-order as the other cases can be proved in a similar fashion.

  Suppose $M$ is right ordered and that $x,y\in M$ are in the same
  $\RR$-class. Then $x=ya$ and $y=xb$ for some $a,b\in M$. This
  implies that $x\leq y$ and $y\leq x$ so that $x=y$.

  Conversely, suppose that all $\RR$-classes are singletons. Then
  $x\le_\RR y$ and $y\le_\RR x$ imply that $x=y$, so that the
  $\RR$-preorder turns into a partial order. Hence $M$ is
  right ordered using $xy \le_\RR x$.
\end{proof}

\subsection{Aperiodic and $\RR$-trivial monoids}

The class of $\HH$-trivial monoids coincides with that of
\emph{aperiodic} monoids (see for example~\cite[Proposition
4.9]{Pin.2009}): a monoid is called \emph{aperiodic} if for any $x\in
M$, there exists some positive integer $N$ such that
$x^{N}=x^{N+1}$. The element $x^\omega := x^{N}=x^{N+1}=x^{N+2}=\cdots$
is then an idempotent (the idempotent $x^\omega$ can in fact be
defined for any element of any monoid~\cite[Chapter VI.2.3]{Pin.2009},
even infinite monoids; however, the period $k$ such that $x^N = x^{N+k}$ need no
longer be $1$). We write $\idempMon := \{x^\omega \mid x\in M\}$
for the set of idempotents of $M$.

Our favorite example of a monoid which is aperiodic, but not
$\RR$-trivial, is the biHecke monoid studied
in~\cite{Hivert_Schilling_Thiery.BiHeckeMonoid.2010,Hivert_Schilling_Thiery.BiHeckeMonoidRepresentation.2010}. 
This is the submonoid of functions from a finite Coxeter group $W$ to itself
generated simultaneously by the elementary bubble sorting and
antisorting operators $\overline{\pi}_i$ and $\pi_i$
\begin{equation}
\label{equation.bihecke}
   \biheckemonoid(W) :=
   \langle \pi_1, \pi_2, \ldots, \pi_n, \opi_1, \opi_2, \ldots, \opi_n \rangle\,.
\end{equation}
See~\cite[Definition 1.1]{Hivert_Schilling_Thiery.BiHeckeMonoid.2010}
and~\cite[Proposition 3.8]{Hivert_Schilling_Thiery.BiHeckeMonoid.2010}.

The smaller class of $\RR$-trivial monoids coincides with the class of
so-called \emph{weakly ordered monoids} as defined by
Schocker~\cite{Schocker.2008}. Also, via the right regular
representation, any $\RR$-trivial monoid can be represented as a monoid
of regressive functions on some finite poset $P$ (a function $f: P\to
P$ is called \emph{regressive} if $f(x) \le x$ for every $x\in P$);
reciprocally any such monoid is $\RR$-trivial.  We now present an
example of a monoid which is $\RR$-trivial, but not $\JJ$-trivial.
\begin{example}
  Take the free left regular band $\mathcal{B}$ generated by two idempotents
  $a,b$.  Multiplication is given by concatenation taking into account the
  idempotent relations, and then selecting only the two left factors (see for
  example~\cite{Saliola.2007}).  So $\mathcal{B} = \{1,a,b,ab,ba\}$ and
  $1\mathcal{B}=\mathcal{B}$, $a\mathcal{B} = \{a,ab\}$, $b\mathcal{B} =\{b,
  ba\}$, $ab \mathcal{B} = \{ab\}$, and $ba\mathcal{B} = \{ba\}$. This shows
  that all $\RR$-classes consist of only one element and hence $\mathcal{B}$
  is $\RR$-trivial.

  On the other hand, $\mathcal{B}$ is not $\LL$-trivial since
  $\{ab,ba\}$ forms an $\LL$-class since $b\cdot ab = ba$ and $a\cdot
  ba = ab$. Hence $\mathcal{B}$ is also not $\JJ$-trivial.
\end{example}

\subsection{$\JJ$-trivial monoids}

The most important for our paper is the class of $\JJ$-trivial
monoids.  In fact, our main motivation stems from the fact that the
submonoid $M_1=\{ f\in M \mid f(1) =1\} $ of the biHecke monoid $M$
in~\eqref{equation.bihecke} of functions that fix the identity, is
$\JJ$-trivial
(see~\cite[Corollary~4.2]{Hivert_Schilling_Thiery.BiHeckeMonoid.2010}
and~\cite{Hivert_Schilling_Thiery.BiHeckeMonoidRepresentation.2010}).

\begin{example} \label{example.J_trivial}
The following example of a $\JJ$-trivial monoid is given in~\cite{Straubing.Therien.1985}. 
Take $M=\{1,x,y,z,0\}$ with relations 
$x^2=x$, $y^2=y$, $xz=zy=z$, and all other products are equal to $0$. Then $M1M = M$, 
$MxM=\{x,z,0\}$, $MyM=\{y,z,0\}$, $MzM=\{z,0\}$, and $M0M=\{0\}$, which shows that $M$ is 
indeed $\JJ$-trivial. Note also that $M$ is left-right ordered with the order $1>x>y>z>0$,
which by Proposition~\ref{proposition.ordered} is equivalent to $\JJ$-triviality.
\end{example}

\subsection{Ordered monoids (with $1$ on top)}
Ordered monoids $M$ with $1$ on top form a subclass of $\JJ$-trivial
monoids.  To see this suppose that $x,y\in M$ are in the same
$\RR$-class, that is $x=ya$ and $y=xb$ for some $a,b\in M$. Since
$a\le 1$, this implies $x=ya \le y$ and $y=xb\le x$ so that
$x=y$. Hence $M$ is $\RR$-trivial. By analogous arguments, $M$ is also
$\LL$-trivial.  Since $M$ is finite, this implies that $M$ is
$\JJ$-trivial (see~\cite[Chapter V, Theorem 1.9]{Pin.2009}).

The next example shows that ordered monoids with 1 on top form a proper
subclass of $\JJ$-trivial monoids.
\begin{example}
The monoid $M$ of Example~\ref{example.J_trivial} is not ordered.
To see this suppose that $\le$ is an order on $M$ with maximal element $1$. 
The relation $y\le 1$ implies $0=z^2\le z = xzy \le xy =0$ which contradicts $z\neq 0$.
\end{example}

It was shown by Straubing and Th\'erien~\cite{Straubing.Therien.1985} and Henckell and 
Pin~\cite{Henckell_Pin.2000} that every $\JJ$-trivial monoid is a quotient of an ordered
monoid with $1$ on top.

In the next two subsections we present two important examples of
ordered monoids with $1$ on top: the $0$-Hecke monoid and the monoid
of regressive order preserving functions, which generalizes nondecreasing parking functions.

\subsection{$0$-Hecke monoids}
\label{ssec:zeroHeckeDefinition}
Let $W$ be a finite Coxeter group. It has a presentation
\begin{equation}
  W = \langle\, s_i \; \text{for} \; i\in I\ \suchthat\  (s_is_j)^{m(s_i,s_j)},\ \forall i,j\in I\,\rangle\,,
\end{equation}
where $I$ is a finite set, $m(s_i,s_j) \in \{1,2,\dots,\infty\}$, and $m(s_i,s_i)=1$.
The elements $s_i$ with $i\in I$ are called \emph{simple reflections}, and the
relations can be rewritten as:
\begin{equation}
  \begin{alignedat}{2}
    s_i^2 &=1 &\quad& \text{ for all $i\in I$}\,,\\
    \underbrace{s_is_js_is_js_i \cdots}_{m(s_i,s_j)} &=
    \underbrace{s_js_is_js_is_j\cdots}_{m(s_i,s_j)} && \text{ for all $i,j\in I$}\, ,
  \end{alignedat}
\end{equation}
where $1$ denotes the identity in $W$. An expression $w=s_{i_1}\cdots s_{i_\ell}$ for $w\in W$ 
is called \emph{reduced} if it is of minimal length $\ell$.
See~\cite{Bjorner_Brenti.2005, Humphreys.1990} for further details on Coxeter groups.

The Coxeter group of type $A_{n-1}$ is the symmetric group $\sg[n]$ with generators
$\{s_1,\dots,s_{n-1}\}$ and relations:
\begin{equation}
  \begin{alignedat}{2}
    s_i^2           & = 1                &     & \text{ for } 1\leq i\leq n-1\,,\\
    s_i s_j         & = s_j s_i            &     & \text{ for } |i-j|\geq2\,, \\
    s_i s_{i+1} s_i & = s_{i+1} s_i s_{i+1} &\quad& \text{ for } 1\leq i\leq n-2\,;
  \end{alignedat}
\end{equation}
the last two relations are called the \emph{braid relations}.

\begin{definition}[\textbf{$0$-Hecke monoid}]
The $0$-Hecke monoid $H_0(W) = \langle \pi_i \mid i \in I \rangle$ of a Coxeter group $W$ 
is generated by the \emph{simple projections} $\pi_i$ with relations
\begin{equation}
  \begin{alignedat}{2}
    \pi_i^2 &=\pi_i &\quad& \text{ for all $i\in I$,}\\
    \underbrace{\pi_i\pi_j\pi_i\pi_j\cdots}_{m(s_i,s_{j})} &=
    \underbrace{\pi_j\pi_i\pi_j\pi_i\cdots}_{m(s_i,s_{j})} && \text{ for all $i,j\in I$}\ .
  \end{alignedat}
\end{equation}
Thanks to these relations, the elements of $H_0(W)$ are canonically
indexed by the elements of $W$ by setting $\pi_w :=
\pi_{i_1}\cdots\pi_{i_k}$ for any reduced word $i_1 \dots i_k$ of $w$.
\end{definition}

\emph{Bruhat order} is a partial order defined on any Coxeter group $W$
and hence also the corresponding $0$-Hecke monoid $H_0(W)$.  Let
$w=s_{i_1} s_{i_2} \cdots s_{i_\ell}$ be a reduced expression for
$w\in W$. Then, in Bruhat order $\le_B$,
\begin{equation*}
	u\le_B w \quad \begin{array}[t]{l}
	\text{if there exists a reduced expression $u = s_{j_1} \cdots s_{j_k}$}\\
	\text{where $j_1 \ldots j_k$ is a subword of $i_1 \ldots i_\ell$.} \end{array}
\end{equation*}
In Bruhat order, $1$ is the minimal element. Hence, it is not hard to check that, with
reverse Bruhat order, the $0$-Hecke monoid is indeed an ordered monoid with $1$
on top.

In fact, the orders $\le_\LL$, $\le_\RR$, $\le_\JJ$, $\le_\BB$ on
$H_0(W)$ correspond exactly to the usual (reversed) left, right,
left-right, and Bruhat order on the Coxeter group $W$.

\subsection{Monoid of regressive order preserving functions}
\label{ssec:ropf}

For any partially ordered set $P$, there is a particular $\JJ$-trivial monoid
which has some very nice properties and that we investigate further in
Section~\ref{sec:NPDF}. Notice that we use the right action in this paper, so
that for $x\in P$ and a function $f:P\to P$ we write $x.f$ for the value of
$x$ under $f$.

\begin{definition}[\textbf{Monoid of regressive order preserving functions}]
  Let $(P, \leq_P)$ be a poset. The set $\OR(P)$ of functions $f: P \to P$ which are
  \begin{itemize}
  \item \emph{order preserving}, that is, for all $x,y\in P,\ x\leq_P y$ implies
    $x.f\leq_P y.f$
  \item \emph{regressive}, that is, for all $x\in P$ one has $x.f \leq_P x$
  \end{itemize}
  is a monoid under composition.
\end{definition}
\begin{proof}
  It is trivial that the identity function is order preserving and regressive and that
  the composition of two order preserving and regressive functions is as well.
\end{proof}

According to~\cite[14.5.3]{Ganyushkin_Mazorchuk.2009}, not much is
known about these monoids.

When $P$ is a chain on $N$ elements, we obtain the monoid $\NDPF_N$ of
nondecreasing parking functions on the set $\{1, \ldots, N\}$ (see
e.g.~\cite{Solomon.1996}; it also is described under the notation
$\mathcal C_n$ in e.g.~\cite[Chapter~XI.4]{Pin.2009} and, together
with many variants, in~\cite[Chapter~14]{Ganyushkin_Mazorchuk.2009}).
The unique minimal set of generators for $\NDPF_N$ is given by the
family of idempotents $(\pi_i)_{i\in\{1,\dots,n-1\}}$, where each
$\pi_i$ is defined by $(i+1).\pi_i:=i$ and $j.\pi_i:=j$ otherwise. The
relations between those generators are given by:
\begin{gather*}
  \pi_i\pi_j = \pi_j\pi_i \quad \text{ for all $|i-j|>1$}\,,\\
  \pi_i\pi_{i-1}=\pi_i\pi_{i-1}\pi_i=\pi_{i-1}\pi_i\pi_{i-1}\,.
\end{gather*}
It follows that $\NDPF_n$ is the natural quotient of $H_0(\sg[n])$ by
the relation $\pi_i\pi_{i+1}\pi_i = \pi_{i+1}\pi_i$, via the quotient
map $\pi_i\mapsto
\pi_i$~\cite{Hivert.Thiery.HeckeSg.2006,Hivert.Thiery.HeckeGroup.2007,
  Ganyushkin_Mazorchuk.2010}. Similarly, it is a natural quotient of
Kiselman's
monoid~\cite{Ganyushkin_Mazorchuk.2010,Kudryavtseva_Mazorchuk.2009}.

To see that $\OR(P)$ is indeed a subclass of ordered monoids with $1$
on top, note that we can define a partial order by saying $f\le g$ for
$f,g\in \OR(P)$ if $x.f \le_P x.g$ for all $x\in P$. By
regressiveness, this implies that $f\le \id$ for all $f\in \OR(P)$ so
that indeed $\id$ is the maximal element. Now take $f,g,h \in \OR(P)$
with $f\le g$.  By definition $x.f \le_P x.g$ for all $x\in P$ and
hence by the order preserving property $(x.f).h \le_P (x.g).h$, so
that $fh\le gh$. Similarly since $f\le g$, $(x.h).f \le_P (x.h).g$ so
that $hf\le hg$. This shows that $\OR(P)$ is ordered.

The submonoid $M_1$ of the biHecke monoid~\eqref{equation.bihecke},
and $H_0(W)\subset M_1$, are submonoids of the monoid of regressive
order preserving functions acting on the Bruhat poset.

\subsection{Monoid of unitriangular Boolean matrices}

Finally, we define the $\JJ$-trivial monoid $\unitribool_n$ of
\emph{unitriangular Boolean matrices}, that is of $n\times n$ matrices
$m$ over the Boolean semi-ring which are unitriangular: $m[i,i]=1$ and
$m[i,j]=0$ for $i>j$. Equivalently (through the adjacency matrix),
this is the monoid of the binary reflexive relations contained in the
usual order on $\{1,\dots,n\}$ (and thus antisymmetric), equipped with
the usual composition of relations. Ignoring loops, it is convenient
to depict such relations by acyclic digraphs admitting $1,\dots,n$ as
linear extension. The product of $g$ and $h$ contains the edges of
$g$, of $h$, as well as the transitivity edges $i\edge k$ obtained
from one edge $i\edge j$ in $g$ and one edge $j\edge k$ in $h$. Hence,
$g^2=g$ if and only if $g$ is transitively closed.

The family of monoids $(\unitribool_n)_n$ (resp. $(\NDPF_n)_n$) plays a
special role, because any $\JJ$-trivial monoid is a subquotient of
$\unitribool_n$ (resp. $\NDPF_n$) for $n$ large
enough~\cite[Chapter~XI.4]{Pin.2009}. In particular, $\NDPF_n$ itself
is a natural submonoid of $\unitribool_n$.
\begin{remark}
  We now demonstrate how $\NDPF_n$ can be realized as a submonoid of relations.
  For simplicity of notation, we consider the monoid $\OR(P)$
  where $P$ is the reversed chain $\{1>\dots>n\}$. Otherwise said,
  $\OR(P)$ is the monoid of functions on the chain $\{1<\dots<n\}$
  which are order preserving and extensive ($x.f\geq x$). Obviously,
  $\OR(P)$ is isomorphic to $\NDPF_n$.

  The monoid $\OR(P)$ is isomorphic to the submonoid of the
  relations $A$ in $\unitribool_n$ such that $i\edge j \in A$ implies
  $k\edge l\in A$ whenever $i\geq k\geq l\geq j$ (in
  the adjacency matrix: $(k,l)$ is to the south-west of $(i,j)$ and
  both are above the diagonal). The isomorphism is given by the map $A
  \mapsto f_A\in \OR(P)$, where
  \begin{equation*}
    u\cdot f_A := \max\{v\suchthat u\,\edge\,v\in A\}\,.
  \end{equation*}
  The inverse bijection $f\in \OR(P) \mapsto A_f\in \unitribool_n$ is given by
  \begin{equation*}
    u\,\edge\,v \in A_f \text{ if and only if } u\cdot f \leq v\,.
  \end{equation*}
  For example, here are the elements of $\OR(\{1>2>3\})$ and the
  adjacency matrices of the corresponding relations in
  $\unitribool_3$:
  \begin{equation*}
    \begin{array}{ccccc}
      \begin{tikzpicture}[->,baseline=(current bounding box.east)]
        \matrix (m) [matrix of math nodes, row sep=.5em, column sep=1.5em]{
          1      & 1 \\
          2      & 2 \\
          3      & 3 \\
        };
        \draw (m-1-1) -> (m-1-2);
        \draw (m-2-1) -> (m-2-2);
        \draw (m-3-1) -> (m-3-2);
      \end{tikzpicture}&
      \begin{tikzpicture}[->,baseline=(current bounding box.east)]
        \matrix (m) [matrix of math nodes, row sep=.5em, column sep=1.5em]{
          1      & 1 \\
          2      & 2 \\
          3      & 3 \\
        };
        \draw (m-1-1) -> (m-2-2);
        \draw (m-2-1) -> (m-2-2);
        \draw (m-3-1) -> (m-3-2);
      \end{tikzpicture}&
      \begin{tikzpicture}[->,baseline=(current bounding box.east)]
        \matrix (m) [matrix of math nodes, row sep=.5em, column sep=1.5em]{
          1      & 1 \\
          2      & 2 \\
          3      & 3 \\
        };
        \draw (m-1-1) -> (m-1-2);
        \draw (m-2-1) -> (m-3-2);
        \draw (m-3-1) -> (m-3-2);
      \end{tikzpicture}&
      \begin{tikzpicture}[->,baseline=(current bounding box.east)]
        \matrix (m) [matrix of math nodes, row sep=.5em, column sep=1.5em]{
          1      & 1 \\
          2      & 2 \\
          3      & 3 \\
        };
        \draw (m-1-1) -> (m-2-2);
        \draw (m-2-1) -> (m-3-2);
        \draw (m-3-1) -> (m-3-2);
      \end{tikzpicture}&
      \begin{tikzpicture}[->,baseline=(current bounding box.east)]
        \matrix (m) [matrix of math nodes, row sep=.5em, column sep=1.5em]{
          1      & 1 \\
          2      & 2 \\
          3      & 3 \\
        };
        \draw (m-1-1) -> (m-3-2);
        \draw (m-2-1) -> (m-3-2);
        \draw (m-3-1) -> (m-3-2);
      \end{tikzpicture}\\\\
      \begin{pmatrix}
        1 & 0 & 0\\
        0 & 1 & 0\\
        0 & 0 & 1\\
      \end{pmatrix} &
      \begin{pmatrix}
        1 & 1 & 0\\
        0 & 1 & 0\\
        0 & 0 & 1\\
      \end{pmatrix} &
      \begin{pmatrix}
        1 & 0 & 0\\
        0 & 1 & 1\\
        0 & 0 & 1\\
      \end{pmatrix} &
      \begin{pmatrix}
        1 & 1 & 0\\
        0 & 1 & 1\\
        0 & 0 & 1\\
      \end{pmatrix} &
      \begin{pmatrix}
        1 & 1 & 1\\
        0 & 1 & 1\\
        0 & 0 & 1\\
      \end{pmatrix} \; .
    \end{array}
  \end{equation*}
\end{remark}

\section{Essential Features of the Representation Theory}
\label{sec:jrep}

In this section we study the representation theory of $\JJ$-trivial
monoids $\tMonoid$, using the $0$-Hecke monoid $H_0(W)$ of a finite
Coxeter group as running example. In Section~\ref{ss.simple.radical}
we construct the simple modules of $\tMonoid$ and derive a description
of the radical $\rad\K\tMonoid$ of the monoid algebra of
$\tMonoid$. We then introduce a star product on the set $\idempMon$ of
idempotents in Theorem~\ref{theorem.star} which makes it into a
semi-lattice, and prove in
Corollary~\ref{corollary.triangular-radical} that the semi-simple
quotient of the monoid algebra $\K\tMonoid/\rad \K\tMonoid$ is the
monoid algebra of $(\idempMon,\star)$. 

\subsection{Simple modules, radical, star product, and semi-simple quotient}
\label{ss.simple.radical}

The goal of this subsection is to construct the simple modules of the
algebra of a $\JJ$-trivial monoid $\tMonoid$, and to derive a
description of its radical and its semi-simple quotient. The proof
techniques are similar to those of Norton~\cite{Norton.1979} for the
$0$-Hecke algebra. However, putting them in the context of
$\JJ$-trivial monoids makes the proofs more transparent. In fact, most
of the results in this section are already known and admit natural
generalizations in larger classes of monoids ($\RR$-trivial, ...). For
example, the description of the radical is a special case of
Almeida-Margolis-Steinberg-Volkov~\cite{Almeida_Margolis_Steinberg_Volkov.2009},
and that of the simple modules
of~\cite[Corollary~9]{Ganyushkin_Mazorchuk_Steinberg.2009}.

Also, the description of the semi-simple quotient is often derived
alternatively from the description of the radical, by noting that it
is the algebra of a monoid which is $\JJ$-trivial and idempotent
(which is equivalent to being a semi-lattice; see
e.g.~\cite[Chapter VII, Proposition 4.12]{Pin.2009}).
\begin{proposition}\label{proposition.simple}
  Let $\tMonoid$ be a $\JJ$-trivial monoid and $x\in \tMonoid$. Let $S_x$ be the
  $1$-dimensional vector space spanned by an element $\epsilon_x$, and define the
  right action of any $y\in \tMonoid$ by
  \begin{equation}
    \epsilon_x y =
    \begin{cases}
      \epsilon_x & \text{if $xy=x$,}\\
      0          & \text{otherwise.}
    \end{cases}
  \end{equation}
  Then $S_x$ is a right $\tMonoid$-module. Moreover, any simple module is isomorphic
  to $S_x$ for some $x \in \tMonoid$ and is in particular one-dimensional.
\end{proposition}

Note that some $S_x$ may be isomorphic to each other, and that the
$S_x$ can be similarly endowed with a left $\tMonoid$-module structure.

\begin{proof}
  Recall that, if $\tMonoid$ is $\JJ$-trivial, then $\leq_\JJ$ is a
  partial order called $\JJ$-order (see
  Proposition~\ref{proposition.ordered}). Let $(x_1, x_2, \ldots,
  x_n)$ be a linear extension of $\JJ$-order, that is an enumeration
  of the elements of $\tMonoid$ such that $x_i \leq_\JJ x_j$ implies
  $i\leq j$. For $0<i\leq
  n$, define $F_i = \K\{x_j \suchthat j\leq i\}$ and set
  $F_0=\{0_\K\}$.  Clearly the $F_i$'s are ideals of $\K\tMonoid$
  such that the sequence
  \begin{equation*}
    F_0 \subset F_1 \subset F_2 \subset \cdots \subset F_{n-1} \subset F_n
  \end{equation*}
  is a composition series for the regular representation
  $F_n=\K\tMonoid$ of $\tMonoid$. Moreover, for any $i>0$, the
  quotient $F_i/F_{i-1}$ is a one-dimensional $\tMonoid$-module
  isomorphic to $S_{x_i}$. Since any simple $\tMonoid$-module must
  appear in any composition series for the regular representation, it
  has to be isomorphic to $F_i/F_{i-1}\cong S_{x_i}$ for some $i$.
\end{proof}

\begin{corollary}
  Let $\tMonoid$ be a $\JJ$-trivial monoid. Then, the quotient of its
  monoid algebra $\K\tMonoid$ by its radical is commutative.
\end{corollary}
Note that the radical $\rad \K\tMonoid$ is not necessarily generated
as an ideal by $\{gh - hg \suchthat g,h\in \tMonoid\}$. For example,
in the commutative monoid $\{1, x, 0\}$ with $x^2=0$, the radical is
$\K (x - 0)$. However, thanks to the following this is true if
$\tMonoid$ is generated by idempotents (see
Corollary~\ref{corollary.rad.idemp}).

The following proposition gives an alternative description of the radical of
$\K\tMonoid$.
\begin{proposition}\label{proposition.basis.radical}
  Let $\tMonoid$ be a $\JJ$-trivial monoid. Then
  \begin{equation}
    \{ x-x^\omega \suchthat x\in \tMonoid\backslash\idempMon \}
  \end{equation}
  is a basis for $\rad\K\tMonoid$.

  Moreover $(S_e)_{e \in \idempMon}$ is a complete set of pairwise non-isomorphic 
  representatives of isomorphism classes of simple
  $\tMonoid$-modules.
\end{proposition}

\begin{proof}
  For any $x,y \in \tMonoid$, either $yx=y$ and then $yx^\omega =y$,
  or $yx <_\JJ y$ and then $yx^\omega <_\JJ y$. Therefore $x-x^\omega$
  is in $\rad\K\tMonoid$ because for any $y$ the product
  $\epsilon_y(x-x^\omega)$ vanishes. Since $x^\omega \leq x$, by
  triangularity with respect to $\JJ$-order, the family
  \begin{equation*}
    \{ x-x^\omega \suchthat x\in \tMonoid\backslash \idempMon\} \cup \idempMon
  \end{equation*}
  is a basis of $\K\tMonoid$. There remains to show that the radical
  is of dimension at most the number of non-idempotents in $\tMonoid$,
  which we do by showing that the simple modules $(S_e)_{e\in
    \idempMon}$ are not pairwise isomorphic. Assume that $S_e$ and
  $S_f$ are isomorphic. Then, since $\epsilon_e e = \epsilon_e$, it
  must be that $\epsilon_e f = \epsilon_e$ so that $ef=e$. Similarly
  $fe=f$, so that $e$ and $f$ are in the same $\JJ$-class and
  therefore equal.
\end{proof}

The following theorem elucidates the structure of the semi-simple quotient of
the monoid algebra $\K\tMonoid$.
\begin{theorem} \label{theorem.star} 
  Let $\tMonoid$ be a $\JJ$-trivial monoid. Define a product $\star$ on
  $\idempMon$ by: 
  \begin{equation}
    e \star f := (e f)^\omega\,.
  \end{equation}
  Then, the restriction of $\leq_\JJ$ on $\idempMon$ is a lattice such
  that
  \begin{equation}
    e \wedge_\JJ f = e \star f\,,
  \end{equation}
  where $e \wedge_\JJ f$ is the meet or infimum of $e$ and $f$ in the
  lattice. In particular $(\idempMon, \star)$ is an idempotent
  commutative $\JJ$-trivial monoid.
\end{theorem}

We start with two preliminary easy lemmas (which are consequences of
e.g.~\cite[Chapter VII, Proposition~4.10]{Pin.2009}).
\begin{lemma}
\label{lemma.idem_factor}
  If $e \in \idempMon$ is such $e = ab$ for some $a,b\in \tMonoid$, then 
  \[
  e=ea=be=ae=eb\,.
  \]
\end{lemma}
\begin{proof}
  For $e\in \idempMon$, one has $e=e^3$ so that $e=eabe$. As a
  consequence, $e\leq_\JJ ea\leq_\JJ e$ and $e\leq_\JJ be\leq_\JJ e$,
  so that $e=ea=be$.  In addition $e=e^2=eab=eb$ and $e=e^2=abe=ae$.
\end{proof}

\begin{lemma}\label{lemma.j.idemp}
  For $e\in \idempMon$ and $y \in \tMonoid$, the following three statements are
  equivalent:
  \begin{equation}
    e \leq_\JJ y, \qquad\qquad e = ey, \qquad\qquad e = ye \;.
  \end{equation}
\end{lemma}
\begin{proof}
  Suppose that $e,y$ are such that $e \leq_\JJ y$. Then $e=ayb$ for some $a,b\in
  \tMonoid$. Applying Lemma~\ref{lemma.idem_factor} we obtain $e=ea=be$ so that 
  $eye = eaybe = eee = e$ since $e\in \idempMon$. A second application of 
  Lemma~\ref{lemma.idem_factor} shows that
  $ey = eye =e$ and $ye = eye = e$.
  The converse implications hold by the definition of $\leq_\JJ$.
\end{proof}

\begin{proof}[Proof of Theorem~\ref{theorem.star}]
  We first show that, for any $e,f \in \idempMon$ the product $e\star
  f$ is the greatest lower bound $e\wedge_\JJ f$ of $e$ and $f$ so that
  the latter exists. It is clear that $(ef)^\omega \leq_\JJ e$ and
  $(ef)^\omega \leq_\JJ f$. Take now $z\in \idempMon$ satisfying
  $z\leq_\JJ e$ and $z\leq_\JJ f$. Applying Lemma~\ref{lemma.j.idemp},
  $z = ze = zf$, and therefore $z = z(ef)^\omega$. Applying
  Lemma~\ref{lemma.j.idemp} backward, $z\leq_\JJ (ef)^\omega$, as
  desired.

  Hence $(\idempMon, \leq_\JJ)$ is a meet semi-lattice with a greatest
  element which is the unit of $\tMonoid$. It is therefore a lattice
  (see e.g.~\cite{stanley97}). Since lower bound is a
  commutative associative operation, $(\idempMon, \star)$ is a
  commutative idempotent monoid.
\end{proof}

We can now state the main result of this section.
\begin{corollary} 
  \label{corollary.triangular-radical}
  Let $\tMonoid$ be a $\JJ$-trivial monoid. Then, $(\K\idempMon,
  \star)$ is isomorphic to $\K\tMonoid/\rad\K\tMonoid$ and $\phi: x
  \mapsto x^\omega$ is the canonical algebra morphism associated to
  this quotient.
\end{corollary}
\begin{proof}
  Denote by $\psi\ :\ \K\tMonoid \to\K\tMonoid/\rad\K\tMonoid$
  the canonical algebra morphism. It follows from
  Proposition~\ref{proposition.basis.radical} that, for any $x$
  (idempotent or not), $\psi(x) = \psi(x^\omega)$ and that
  $\{\psi(e)\suchthat e\in \idempMon\}$ is a basis for the
  quotient. Finally, $\star$ coincides with the product in the
  quotient: for any $e,f\in \idempMon$,
  \begin{displaymath}
    \psi(e)\psi(f) = \psi(ef) = \psi((ef)^\omega) = \psi(e\star f)\,.\qedhere
  \end{displaymath}
\end{proof}

\begin{corollary}
  \label{corollary.rad.idemp}
  Let $\tMonoid$ be a $\JJ$-trivial monoid generated by
  idempotents. Then the radical $\rad\K\tMonoid$ of its monoid algebra is
  generated as an ideal by 
  \begin{equation}
    \{gh - hg \suchthat g,h\in \tMonoid\}\,.
  \end{equation}
\end{corollary}
\begin{proof}
  \newcommand{\Com}{\mathcal{C}}
  Denote by $\Com$ the ideal generated by $\{gh - hg \suchthat g,h\in
  \tMonoid\}$. Since $\rad\K\tMonoid$ is the linear span of
  $(x-x^\omega)_{x\in\tMonoid}$, it is sufficient to show that for any
  $x\in\tMonoid$ one has $x\equiv x^2 \pmod \Com$. Now write $x=e_1\cdots e_n$
  where $e_i$ are all idempotent. Then,
  \begin{equation*}
    x \equiv e_1^2\cdots e_n^2 \equiv e_1\cdots e_ne_1\cdots e_n 
    \equiv x^2 \pmod \Com\,.\qedhere
  \end{equation*}
\end{proof}

\begin{example}[Representation theory of $H_0(W)$]
\label{example.zero.hecke}
Consider the $0$-Hecke monoid $H_0(W)$ of a finite Coxeter group $W$, with index set $I=\{1, 2, \ldots, n\}$. 
For any $J \subseteq I$, we can consider the parabolic submonoid $H_0(W_J)$ generated by
$\{\pi_i \mid i\in J\}$.  Each parabolic submonoid contains a unique longest element $\longest_J$.  
The collection $\{\longest_J \mid J\subseteq I\}$ is exactly the set of idempotents in $H_0(W)$.

For each $i\in I$, we can construct the \emph{evaluation maps} $\Phi_i^+$ and $\Phi_i^-$ defined 
on generators by:
\begin{eqnarray*}
\Phi_i^+ &:& \mathbb{C}H_0(W) \rightarrow \mathbb{C}H_0(W_{I\setminus \{i\}}) \\
\Phi_i^+(\pi_j) &=&     \begin{cases}
      1          & \text{if $i=j$,}\\
      \pi_j & \text{if $i \neq j$,}
    \end{cases}
\end{eqnarray*}
and
\begin{eqnarray*}
\Phi_i^- &:& \mathbb{C}H_0(W) \rightarrow \mathbb{C}H_0(W_{I\setminus \{i\}}) \\
\Phi_i^-(\pi_j) &=&     \begin{cases}
      0          & \text{if $i=j$,}\\
      \pi_j & \text{if $i \neq j$.}
    \end{cases}
\end{eqnarray*}
One can easily check that these maps extend to algebra morphisms from 
$H_0(W)\rightarrow H_0(W_{I\setminus \{i\}})$.  For any $J$, define $\Phi_J^+$ as the composition 
of the maps $\Phi_i^+$ for $i\in J$, and define $\Phi_J^-$ analogously (the map $\Phi_J^+$ is 
the \emph{parabolic map} studied by Billey, Fan, and Losonczy~\cite{Billey_Fan_Lsonczy.1999}).  
Then, the simple representations of $H_0(W)$ are given by the maps 
$\lambda_J = \Phi_J^+ \circ \Phi_{\hat{J}}^-$, where $\hat{J}=I\setminus J$.  
This is clearly a one-dimensional representation.
\end{example}

\subsubsection{Projective modules}

An important result is that the projective modules for a $\JJ$-trivial monoid are \emph{combinatorial}.  This result uses the derivation of the Cartan matrix, which can be found in the full paper~\cite{dhst.2011}.  Thus, we state the theorem here without proof, and then examine the result in the context of the zero-Hecke monoid.

\begin{theorem}
  \label{theorem.projective_modules}
  For any idempotent $e$ denote by $R(e) = eM$,
  \begin{equation*}
    R_=(e) = \{x\in eM\suchthat \lfix{x} = e\}
    \quad\text{and}\quad 
    R_<(e) = \{x\in eM\suchthat \lfix{x} <_\RR e\} \, .
  \end{equation*}
  Then, the projective module $P_e$ associated to $S_e$ is isomorphic
  to $\K R(e)/\K R_<(e)$. In particular, the projective module $P_e$
  is combinatorial: taking as basis the image of $R_=(e)$ in the
  quotient, the action of $m\in\tMonoid$ on $x\in R_=(e)$ is given by:
  \begin{equation}
    x \cdot m =
    \begin{cases}
      xm &\text{if $\lfix{xm} = e$,}\\
      0  & \text{otherwise}.
    \end{cases}
  \end{equation}
\end{theorem}

\begin{corollary}
  The family $\{b_x \suchthat \lfix{x} = e\}$ is a basis for the right
  projective module associated to $S_e$.
\end{corollary}

\begin{figure}
  \includegraphics[width=\textwidth]{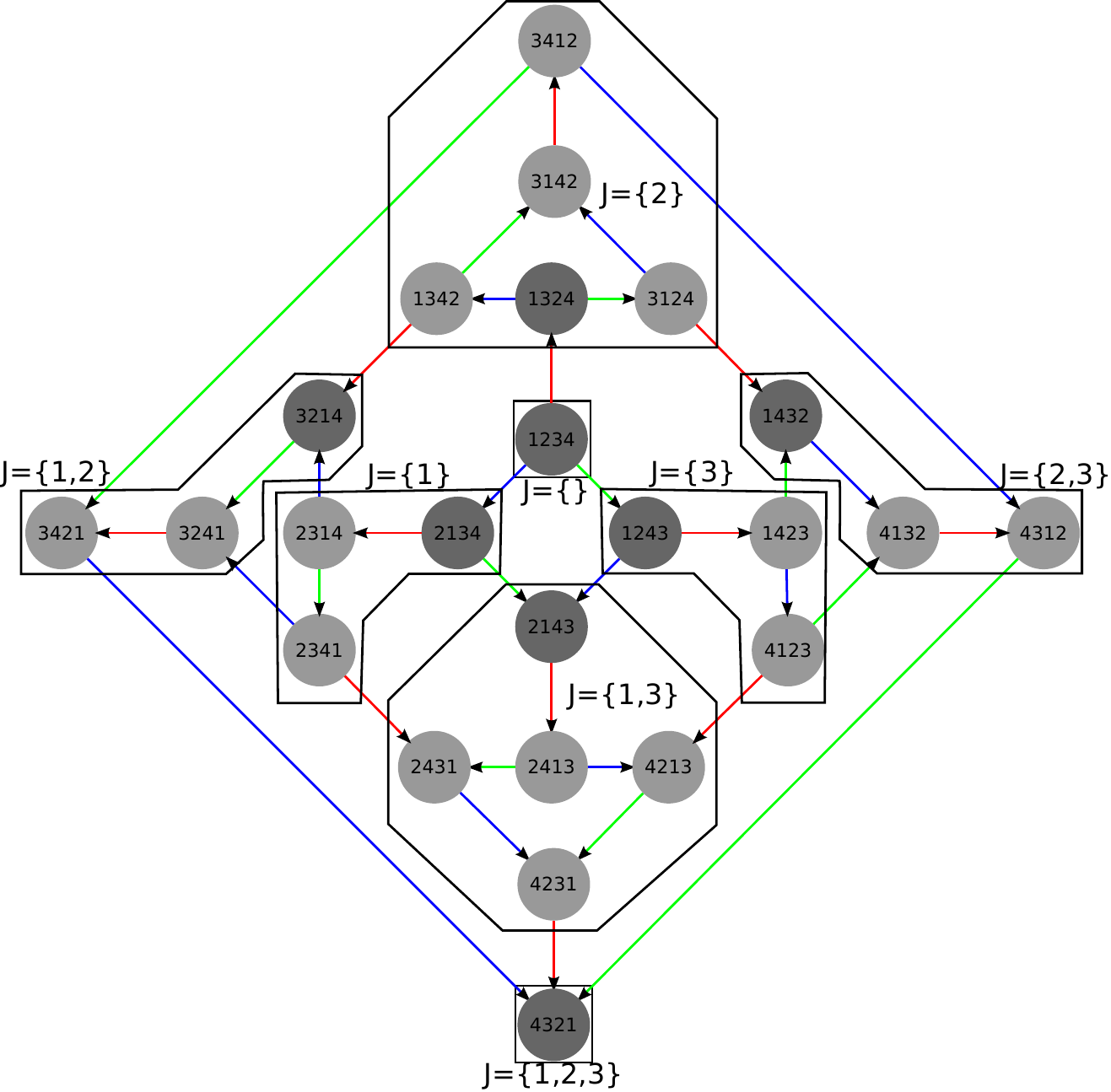}
  \caption{The decomposition of $H_0(\sg[4])$ into indecomposable
    right projective modules. This decomposition follows the partition
    of $\sg[4]$ into left descent classes, each labelled by its
    descent set $J$.  The blue, red, and green lines indicate the
    action of $\pi_1, \pi_2,$ and $\pi_3$ respectively.  The darker circles
    indicate idempotent elements of the monoid. }
  \label{h0s4.projectives}
\end{figure}

\begin{example}[Representation theory of $H_0(W)$, continued]
  \label{example.zero.hecke.projectives}
  The right projective modules of $H_0(W)$ are combinatorial, and
  described by the decomposition of the right order along left descent
  classes, as illustrated in Figure~\ref{h0s4.projectives}. Namely,
  let $P_J$ be the right projective module of $H_0(W)$ corresponding
  to the idempotent $\pi_J$. Its basis $b_w$ is indexed by the
  elements of $w$ having $J$ as left descent set. The action of
  $\pi_i$ coincides with the usual right action, except that
  $b_w.\pi_i=0$ if $w.\pi_i$ has a strictly larger left descent set than
  $w$.

  Here we reproduce Norton's construction of $P_J$~\cite{Norton.1979},
  as it is close to an explicit description of the isomorphism in the
  proof of Theorem~\ref{theorem.projective_modules}. First, notice
  that the elements $\{\pi_i^-=(1-\pi_i) \mid i\in I\}$ are idempotent
  and satisfy the same Coxeter relations as the $\pi_i$.  Thus, the
  set $\{\pi_i^-\}$ generates a monoid isomorphic to $H_0(W)$. For each
  $J\subseteq I$, let $\longest_J^-$ be the longest element in the
  parabolic submonoid associated to $J$ generated by the $\pi_i^-$
  generators, and $\longest_J^+=\longest_J$. For each subset
  $J\subseteq I$, let $\hat{J}=I\setminus J$. Define
  $f_J=\longest_{\hat{J}}^-\longest_J^+$. Then, $f_J \pi_w=0$ if
  $J\subset D_L(w)$. It follows that the right module $f_J H_0(W)$ is
  isomorphic to $P_J$ and its basis $\{ f_J\pi_w\suchthat D_L(w)=J\}$
  realizes the combinatorial module of $P_J$.

  One should notice that the elements
  $\longest_{\hat{J}}^-\longest_J^+$ are, in general, neither
  idempotent nor orthogonal. Furthermore,
  $\longest_{\hat{J}}^-\longest_J^+H_0(W)$ is not a submodule of
  $\pi_J H_0(W)$ as in the proof of
  Theorem~\ref{theorem.projective_modules}.

  The description of left projective modules is symmetric.
\end{example}

\section{Monoid of order preserving regressive functions on a poset $P$}
\label{sec:NPDF}

In this section, we discuss the monoid $\OR(P)$ of order preserving
regressive functions on a poset $P$.  Recall that this is the monoid
of functions $f$ on $P$ such that for any $x\leq y \in P$, $x.f\leq x$ and 
$x.f\leq y.f$.

In Section~\ref{sub:combIdem}, we discuss constructions for idempotents in 
$\OR(P)$ in terms of the image sets of the idempotents, as well as methods for 
obtaining $\lfix f$ and $\rfix f$ for any given function $f$.  In Section~\ref{sub:cartNDPF}, 
we show that the Cartan matrix for $\OR(P)$ is upper uni-triangular with respect to the lexicographic order associated to any linear 
extension of $P$.  In Section~\ref{sub:ndpfSemis}, we specialize to $\OR(L)$ where $L$ is a 
meet semi-lattice, describing a minimal generating set of idempotents.  Finally, in 
Section~\ref{sub:ndpfOrthIdem}, we describe a simple construction for a set of orthogonal 
idempotents in $\NDPF_N$, and present a conjectural construction for orthogonal idempotents 
for $\OR(L)$.

\subsection{Combinatorics of idempotents}
\label{sub:combIdem}

The goal of this section is to describe the idempotents in $\OR(P)$ using
order considerations. We begin by giving the definition of joins, even in the setting when
the poset $P$ is not a lattice.

\begin{definition}
  Let $P$ be a finite poset and $S\subseteq P$. Then $z\in P$ is called \emph{a
    join} of $S$ if $x \leq z$ holds for any $x\in S$, and $z$ is minimal with
  that property.

  We denote $\joins(S)$ the set of joins of $S$, and $\joins(x,y)$ for short
  if $S=\{x,y\}$. If $\joins(S)$ (resp. $\joins(x,y)$) is a singleton (for
  example because $P$ is a lattice) then we denote $\bigvee S$ (resp. $x\vee y$)
  the unique join.  Finally, we define $\joins(\emptyset)$ to be the set of minimal elements in $P$.
\end{definition}

\begin{lemma}
  \label{lemma.fix}
  Let $P$ be some poset, and $f\in \OR(P)$. If $x$ and $y$ are fixed
  points of $f$, and $z$ is a join of $x$ and $y$, then $z$ is a fixed
  point of $f$.
\end{lemma}
\begin{proof}
  Since $x\leq z$ and $y\leq z$, one has $x=x.f \leq z.f$ and $y=y.f \leq z.f$.
  Since furthermore $z.f\le z$, by minimality of $z$ the equality $z.f = z$ must hold.
\end{proof}

\begin{lemma}
  \label{lemma.sup}
  Let $I$ be a subset of $P$ which contains all the minimal elements
  of $P$ and is stable under joins. Then, for any $x\in P$, the set
  $\{y\in I\suchthat y\leq x\}$ admits a unique maximal element which
  we denote by $\sup_I(x)\in I$. Furthermore, the map $\sup_I:
  x\mapsto \sup_I(x)$ is an idempotent in $\OR(P)$.
\end{lemma}
\begin{proof}
   For the first statement, suppose for some $x \not \in I$ there are two 
   maximal elements $y_1$ and $y_2$ in $\{y\in I\suchthat y\leq x\}$.  
   Then the join $y_1 \wedge y_2 < x$, since otherwise $x$ would be a join of $y_1$ 
   and $y_2$, and thus $x\in I$ since $I$ is join-closed.  But this contradicts the maximality 
   of $y_1$ and $y_2$, so the first statement holds.

  Using that $\sup_I(x)\leq x$ and $\sup_I(x)\in I$, $e:=\sup_I$ is
  a regressive idempotent by construction. Furthermore, it is is order
  preserving: for $x\leq z$, $x.e$ and $z.e$ must be comparable or
  else there would be two maximal elements in $I$ under $z$.  Since
  $z.e$ is maximal under $z$, we have $z.e\geq x.e$.
\end{proof}

Reciprocally, all idempotents are of this form:
\begin{lemma}
  \label{lemma.idem_image}
  Let $P$ be some poset, and $f\in \OR(P)$ be an idempotent. Then the
  image $\im(f)$ of $f$ satisfies the following:
  \begin{enumerate}
  \item \label{item.min} All minimal elements of $P$ are contained in $\im(f)$.
  \item \label{item.fix} Each $x \in \im(f)$ is a fixed point of $f$.
  \item \label{item.join} The set $\im(f)$ is stable under joins: if
    $S\subseteq \im(f)$ then $\joins(S)\subseteq \im(f)$ .
  \item \label{item.image} For any $x\in P$, the image $x.f$ is the upper bound $\sup_{\im(f)}(x)$.
  \end{enumerate}
\end{lemma}
\begin{proof}  
  Statement~\eqref{item.min} follows from the fact that $x.f\le x$ so that minimal elements
  must be fixed points and hence in $\im(f)$.

  For any $x=a.f$, if $x$ is not a fixed point then $x.f=(a.f).f\neq a.f$, contradicting the 
  idempotence of $f$.  Thus, the second statement holds.
  
  Statement~\eqref{item.join} follows directly from the second statement and Lemma~\ref{lemma.fix}.
  
  If $y\in \im(f)$ and $y\leq x$ then $y = y.f \leq x.f$.  Since this holds for every element of 
  $\{y \in \im(f) \mid y\leq x \}$ and $x.f$ is itself in this set, statement~\eqref{item.image} holds.
\end{proof}

Thus, putting together Lemmas~\ref{lemma.sup}
and~\ref{lemma.idem_image} one obtains a complete description of the
idempotents of $\OR(P)$.
\begin{proposition}
  \label{proposition.idempotents.OO}
  The idempotents of $\OR(P)$ are given by the maps $\sup_I$, where $I$
  ranges through the subsets of $P$ which contain the minimal elements
  and are stable under joins.
\end{proposition}

For $f\in \OR(P)$ and $y\in P$, let $f^{-1}(y)$ be the fiber of $y$ under $f$, that is, the set of
all $x\in P$ such that $x.f=y$.
\begin{definition}
Given $S$ a subset of a finite poset $P$, set $C_0(S)=S$ and $C_{i+1}(S)=C_i(S) \cup 
\{x\in P \mid x \text{ is a join of some elements in } C_i(S) \}$.  Since $P$ is finite, there exists some $N$ 
such that $C_N(S)=C_{N+1}(S)$.  The \emph{join closure} is defined as this stable set, 
and denoted $C(S)$.  A set is \emph{join-closed} if $C(S)=S$. Define
\[
F(f) := \bigcup_{y\in P} \{ x\in f^{-1}(y) \mid x \text{ minimal in } f^{-1}(y) \}
\]
to be the collection of minimal points in the fibers of $f$.
\end{definition}

\begin{corollary}
Let $X$ be the join-closure of the set of minimal points of $P$.  Then $X$ is fixed by every 
$f\in \OR(P)$.
\end{corollary}

\begin{lemma}[Description of left and right symbols]
For any $f\in \OR(P)$, there exists a minimal idempotent $f_r$ whose image set is $C(\im(f))$, 
and $f_r=\rfix{f}$.  There also exists a minimal idempotent $f_l$ whose image set is $C(F(f))$, 
and $f_l = \lfix{f}$.
\end{lemma}

\begin{proof}
The $\rfix{f}$ must fix every element of $\im(f)$, and the image of $\rfix{f}$ must be join-closed
by Lemma~\ref{lemma.idem_image}.
$f_r$ is the smallest idempotent satisfying these requirements, and is thus the $\rfix{f}$.

Likewise, $\lfix{f}$ must fix the minimal elements of each fiber of $f$, and so must fix all of $C(F(f))$.  
For any $y \not \in F(f)$, find $x\leq y$ such that $x.f=y.f$ and $x\in F(f)$.  Then 
$x= x.f_l \leq y.f_l \leq y$.
For any $z$ with $x\leq z \leq y$, we have $x.f\leq z.f \leq y.f=x.f$, so $z$ is in the 
same fiber as $y$.  Then we have $(y.f_l).f =y.f$, so $f_l$ fixes $f$ on the left.  
Minimality then ensures that $f_l=\lfix{f}$.
\end{proof}

Let $P$ be a poset, and $P'$ be the poset obtained by removing a
maximal element $x$ of $P$. Then, the following rule holds:
\begin{proposition}[Branching of idempotents]
  Let $e=\sup_I$ be an idempotent in $\OR(P')$. If $I\subseteq P$ is
  still stable under joins in $P$, then there exist two idempotents in
  $\OR(P)$ with respective image sets $I$ and $I\cup \{x\}$.
  Otherwise, there exists an idempotent in $\OR(P)$ with image set
  $I\cup \{x\}$. Every idempotent in $\OR(P)$ is uniquely obtained by
  this branching.
\end{proposition}
\begin{proof}
  This follows from straightforward reasoning on the subsets $I$ which contain the
  minimal elements and are stable under joins, in $P$ and in $P'$.
\end{proof}

\subsection{The Cartan matrix for $\OR(P)$ is upper uni-triangular}
\label{sub:cartNDPF}

We have seen that the left and right fix of an element of $\OR(P)$
can be identified with the subsets of $P$ closed under joins. We put a
total order $\leq_\lex$ on such subsets by writing them as bit vectors
along a linear extension $p_1,\dots,p_n$ of $P$, and comparing those
bit vectors lexicographically.

\begin{proposition}
  Let $f\in \OR(P)$. Then, $\im(\lfix f) \leq_\lex \im(\rfix f)$, with
  equality if and only if $f$ is an idempotent.
\end{proposition}
\begin{proof}
  Let $n=|P|$ and $p_1,\ldots,p_n$ a linear extension of $P$.
  For $k\in \{0,\dots,n\}$ set respectively $L_k=\im(\lfix f) \cap
  \{p_1,\dots,p_k\}$ and $R_k = \im(\rfix f) \cap \{p_1,\dots,p_k\}$.

  As a first step, we prove the property $(H_k)$: if $L_k=R_k$ then
  $f$ restricted to $\{p_1,\dots,p_k\}$ is an idempotent with image
  set $R_k$. Obviously, $(H_0)$ holds. Take now $k>0$ such that
  $L_k=R_k$; then $L_{k-1}=R_{k-1}$ and we may use by induction
  $(H_{k-1})$.

  Case 1: $p_k\in F(f)$, and is thus the smallest point in its fiber.  This implies that
  $p_k\in L_k$, and by assumption, $L_k=R_k$. 
  By $(H_{k-1})$, $p_k.f<_\lex p_k$ gives a contradiction:
  $p_k.f\in R_{k-1}$, and therefore $p_k.f$ is in the same fiber as
  $p_k$. Hence $p_k.f=p_k$.

  Case 2: $p_k \in C(F(f))=\im(\lfix{f})$, but $p_k \not \in F(f)$.   Then $p_k$ is a join of 
  two smaller elements $x$ and $y$ of $L_k=R_k$; in particular, $p_k\in R_k$. By induction, 
  $x$ and $y$ are fixed by $f$, and therefore $p_k.f=p_k$ by Lemma~\ref{lemma.fix}.

  Case 3: $p_k\not \in C(F(f))=\im(\lfix f)$; then $p_k$ is not a minimal element
  in its fiber; taking $p_i<_\lex p_k$ in the same fiber, we have $(p_k.f).f =
  (p_i.f).f = p_i.f = p_k.f$. Furthermore, $R_k=R_{k-1} =
  \{p_1,\dots,p_{k-1}\}.f = \{p_1,\dots,p_k\}.f$.

  In all three cases above, we deduce that $f$ restricted to
  $\{p_1,\dots,p_k\}$ is an idempotent with image set $R_k$, as desired.

  If $L_n=R_n$, we are done. Otherwise, take $k$ minimal such that
  $L_k\ne R_k$. Assume that $p_k\in L_k$ but not in $R_k$. In
  particular, $p_k$ is not a join of two elements $x$ and $y$ in
  $L_{k-1}=R_{k-1}$; hence $p_k$ is minimal in its fiber, and by the
  same argument as in Case 3 above, we get a contradiction.
\end{proof}

\begin{corollary}
  The Cartan matrix of $\OR(P)$ is upper uni-triangular with respect to the
  lexicographic order associated to any linear extension of $P$.
\end{corollary}

\begin{problem}
  Find larger classes of monoids where this property still holds.
  Note that this fails for the $0$-Hecke monoid which is a submonoid
  of an $\OR(B)$ where $B$ is Bruhat order.
\end{problem}

\subsection{Restriction to meet semi-lattices}
\label{sub:ndpfSemis}

For the remainder of this section, let $L$ be a \emph{meet semi-lattice} and we consider 
the monoid $\OR(L)$. Recall that $L$ is a meet semi-lattice if every pair of elements $x,y\in L$
has a unique meet.

For $a\geq b$, define an idempotent $e_{a,b}$ in $\OR(L)$ by:
\[
    x.e_{a,b} =
    \begin{cases}
        x\wedge b & \text{if $x\leq a$,}\\
        x         & \text{otherwise.}
    \end{cases}
\]

\begin{remark}
  \label{remark.oo.eab}
  The function $e_{a,b}$ is the (pointwise) largest element of $\OR(L)$ such that $a.f=b$.

  For $a\geq b\geq c$, $e_{a,b} e_{b,c} = e_{a,c}$. In the case where
  $L$ is a chain, that is $\OR(L)=\NDPF_{|L|}$, those idempotents
  further satisfy the following braid-like relation: $e_{b,c} e_{a,b}
  e_{b,c} = e_{a,b} e_{b,c} e_{a,b} = e_{a,c}$.
\end{remark}
\begin{proof}
  The first statement is clear. Take now $a\geq b\geq c$ in a
  meet semi-lattice. For any $x\leq a$, we have $x.e_{a,b}=x\wedge b
  \leq b,$ so $x.(e_{a,b}e_{b,c}) = x\wedge b \wedge c = x \wedge c$,
  since $b\ge c$.  On the other hand, $x.e_{a,c} = x\wedge c$, which
  proves the desired equality.

  Now consider the braid-like relation in $\NDPF_{|L|}$.  Using the previous
  result, one gets that $e_{b,c} e_{a,b} e_{b,c}=e_{b,c} e_{a,c}$ and $e_{a,b}
  e_{b,c} e_{a,b}=e_{a,c} e_{a,b}$.  For $x> a$, $x$ is fixed by
  $e_{a,c}$, $e_{a,b}$ and $e_{b,c}$, and is thus fixed by the
  composition. The other cases can be checked analogously.
\end{proof}

\begin{proposition}
  The family $(e_{a,b})_{a,b}$, where $(a,b)$ runs through the covers of $L$, minimally generates the idempotents of $\OR(L)$.
\end{proposition}

\begin{proof}
Given $f$ idempotent in $\OR(L)$, we can factorize $f$ as a product of the idempotents $e_{a,b}$.  Take a 
linear extension of $L$, and recursively assume that $f$ is the identity on all 
elements above some least element $a$ of the linear extension.  
Then define a function $g$ by: 
\[
x.g= \begin{cases}
	a & \text{if $x=a$,} \\
	x.f & \text{otherwise.}
	\end{cases}
\]
We claim that $f = ge_{a,a.f},$ and $g \in \OR(L)$.  There are a number of cases that must be checked:
\begin{itemize}
\item Suppose $x<a$.  Then $x.ge_{a,a.f} = (x.f).e_{a,a.f} = x.f \wedge a.f = x.f$, since $x<a$ implies $x.f<a.f$.
\item Suppose $x>a$.  Then $x.ge_{a,a.f} = (x.f).e_{a,a.f} = x.e_{a,a.f} = x=x.f$, since $x$ is fixed by $f$ by assumption.
\item Suppose $x$ not related to $a$, and $x.f\leq a.f$.  Then $x.ge_{a,a.f} = (x.f).e_{a,a.f} = x.f$.
\item Suppose $x$ not related to $a$, and $a.f\leq x.f\leq a$.  By the idempotence of $f$ we have $a.f=a.f.f\le x.f.f\le a.f$, 
so $x.f=a.f$, which reduces to the previous case.
\item Suppose $x$ not related to $a$, but $x.f\leq a$.  Then by idempotence of $f$ we have $x.f=x.f.f\leq a.f$, 
reducing to a previous case.
\item For $x$ not related to $a$, and $x.f$ not related to $a$ or $x.f>a$, we have $x.f$ fixed by $e_{a,a.f}$, which 
implies that $x.ge_{a,a.f}=x.f$.
\item Finally for $x=a$ we have $a.ge_{a,a.f} = a.e_{a,a.f}=a\wedge a.f=a.f$.
\end{itemize}
Thus, $f = ge_{a,a.f}$.

For all $x\le a$, we have $x.f\le a.f\le a$, so that $x.g\le a.g=a$.  For all
$x>a$, we have $x$ fixed by $g$ by assumption, and for all other $x$, the
$\OR(L)$ conditions are inherited from $f$.  Thus $g$ is in $\OR(L)$.

For all $x\neq a$, we have $x.g=x.f=x.f.f$.  Since all $x>a$ are fixed by $f$, there is no $y$ such that $y.f=a$.  Then $x.f.f=x.g.g$ for all $x\neq a$.  Finally, $a$ is fixed by $g$, so $a=a.g.g$.  Thus $g$ is idempotent.

Applying this procedure recursively gives a factorization of $f$ into a composition of functions $e_{a,a.f}$.  We can further refine this factorization using Remark~\ref{remark.oo.eab} on each $e_{a,a.f}$ by
$e_{a,a.f}=e_{a_0,a_1}e_{a_1,a_2}\cdots e_{a_{k-1},a_k}$, where $a_0=a$, $a_k=a.f$, and
$a_i$ covers $a_{i-1}$ for each $i$.  Then we can express $f$ as a product of functions $e_{a,b}$ where $a$ covers $b$.

This set of generators is minimal because $e_{a,b}$ where $a$ covers $b$ is the pointwise largest function in $\OR(L)$ mapping $a$ to $b$.
\end{proof}

As a byproduct of the proof, we obtain a canonical factorization of any idempotent $f\in \OR(L)$.

\begin{example}
The set of functions $e_{a,b}$ do not in general generate $\OR(L)$.  Let $L$ be the Boolean lattice on three elements.  Label the nodes of $L$ by triples $ijk$ with $i,j,k\in \{0,1\}$, and $abc\geq ijk$ if $a\leq i, b\leq j, c\leq k$.

Define $f$ by $f(000)=000$, $f(100)=110, f(010)=011, f(001)=101$, and $f(x)=111$ for all other $x$.  Simple inspection shows that $f\neq ge_{a,a.f}$ for any choice of $g$ and $a$.  
\end{example}

\subsection{Orthogonal idempotents}
\label{sub:ndpfOrthIdem}

For $\{1,2,\ldots, N\}$ a chain, one can explicitly write down
orthogonal idempotents for $\NDPF_N$. Recall that the minimal
generators for $\NDPF_N$ are the elements $\pi_i =e_{i+1,i}$ and that
$\NDPF_N$ is the quotient of $H_0(\sg[n])$ by the extra relation
$\pi_i\pi_{i+1}\pi_i = \pi_{i+1}\pi_i$, via the quotient map
$\pi_i\mapsto \pi_i$. By analogy with the $0$-Hecke algebra, set
$\pi_i^+=\pi_i$ and $\pi_i^-=1-\pi_i$.

We observe the following relations, which can be checked easily.
\begin{lemma}\label{lem:ndpfrels}
Let $k=i-1$.  Then the following relations hold:
    \begin{enumerate}
      \item $\pi_{i-1}^+ \pi_i^+ \pi_{i-1}^+ = \pi_i^+ \pi_{i-1}^+$,
      \item $\pi_{i-1}^- \pi_i^- \pi_{i-1}^- = \pi_{i-1}^- \pi_i^-$,
      \item $\pi_i^+ \pi_{i-1}^- \pi_i^+ = \pi_i^+ \pi_{i-1}^-$,
      \item $\pi_i^- \pi_{i-1}^+ \pi_i^- = \pi_{i-1}^+ \pi_i^-$,
      \item $\pi_{i-1}^+ \pi_i^- \pi_{i-1}^+ = \pi_i^- \pi_{i-1}^+$,
      \item $\pi_{i-1}^- \pi_i^+ \pi_{i-1}^- = \pi_{i-1}^- \pi_i^+$.
    \end{enumerate}
\end{lemma}

\begin{definition}
  Let $D$ be a \emph{signed diagram}, that is an assignment of a $+$
  or $-$ to each of the generators of $\NDPF_N$. By abuse of notation,
  we will write $i\in D$ if the generator $\pi_i$ is assigned a $+$
  sign.  Let $P=\{ P_1, P_2, \ldots, P_k\}$ be the partition of the
  generators such that adjacent generators with the same sign are in
  the same set, and generators with different signs are in different
  sets.  Set $\epsilon(P_i)\in \{+,-\}$ to be the sign of the subset $P_i$.
  Let $\longest_{P_i}^{\epsilon(P_i)}$ be the longest element in the generators in
  $P_i$, according to the sign in $D$.  Define:
\begin{itemize}
  \item $L_D := \longest_{P_1}^{\epsilon(P_1)}\longest_{P_2}^{\epsilon(P_2)}\cdots \longest_{P_k}^{\epsilon(P_k)}$, 
  \item $R_D := \longest_{P_k}^{\epsilon(P_k)}\longest_{P_k-1}^{\epsilon(P_{k-1})}\cdots \longest_{P_1}^{\epsilon(P_1)}$, 
  \item and $C_D := L_DR_D$.
\end{itemize}
\end{definition}

\begin{example}
\label{example.D}
Let $D=++++---++$.  Then $P=\{ \{1,2,3,4\}, \{5,6,7\}, \{8,9\}  \}$, and the associated 
long elements are: $\longest_{P_1}^+=\pi_4^+ \pi_3^+ \pi_2^+ \pi_1^+$, 
$\longest_{P_2}^-=\pi_5^- \pi_6^- \pi_7^-$, and $\longest_{P_3}^+=\pi_9^+ \pi_8^+$. Then
\begin{equation*}
\begin{split}
	L_D &= \longest_{P_1}^+ \longest_{P_2}^- \longest_{P_3}^+ 
	= (\pi_4^+ \pi_3^+ \pi_2^+ \pi_1^+) (\pi_5^- \pi_6^- \pi_7^-) (\pi_9^+ \pi_8^+),\\
	R_D &= \longest_{P_3}^+ \longest_{P_2}^- \longest_{P_1}^+
	= (\pi_9^+ \pi_8^+) (\pi_5^- \pi_6^- \pi_7^-) (\pi_4^+ \pi_3^+ \pi_2^+ \pi_1^+).
\end{split}
\end{equation*}
\end{example}

The elements $C_D$ are the images, under the natural quotient map from
the $0$-Hecke algebra, of the \emph{diagram demipotents} constructed
in~\cite{Denton.2010.FPSAC,Denton.2011}. An element $x$ of an algebra
is \emph{demipotent} if there exists some finite integer $n$ such that
$x^n=x^{n+1}$ is idempotent. It was shown
in~\cite{Denton.2010.FPSAC,Denton.2011} that, in the $0$-Hecke
algebra, raising the diagram demipotents to the power $N$ yields a set
of primitive orthogonal idempotents for the $0$-Hecke algebra. It
turns out that, under the quotient to $\NDPF_N$, these elements $C_D$
are right away orthogonal idempotents, which we prove now.

\begin{remark}
  \label{remark.fix_i}
  Fix $i$, and assume that $f$ is an element in the monoid generated
  by $\pi^-_{i+1},...,\pi^-_N$ and $\pi^+_{i+1},...,\pi^+_N$. Then,
  applying repeatedly Lemma~\ref{lem:ndpfrels} yields
  $$\pi^-_i f \pi^-_i = \pi^-_i f \qquad \text{and} \qquad \pi^+_i f \pi^+_i = f \pi^+_i\,.$$
\end{remark}

The following proposition states that the elements $C_D$ are also the
images of Norton's generators of the projective modules of the
$0$-Hecke algebra through the natural quotient map to $\NDPF_N$.
\begin{proposition}
  \label{proposition.norton_ndpf}
  Let $D$ be a signed diagram. Then,
  \begin{displaymath}
    C_D = \prod_{i=1,\dots,n,\  i \not\in D} \pi^-_i  \quad \prod_{i=n,\dots,1,\  i\in D} \pi^+_i\,.
  \end{displaymath}
  In other words $C_D$ reduces to one of the following two forms:
  \begin{itemize}
  \item $C_D = (\longest_{P_1}^-\longest_{P_3}^-\cdots \longest_{P_{2k\pm 1}}^-) (\longest_{P_2}^+\longest_{P_4}^+\cdots \longest_{P_{2k}}^+)$, or
  \item $C_D = (\longest_{P_2}^-\longest_{P_4}^-\cdots \longest_{P_{2k}}^-) (\longest_{P_1}^+\longest_{P_3}^+\cdots \longest_{P_{2k\pm 1}}^+)$.
  \end{itemize}
\end{proposition}
\begin{proof}
  Let $D$ be a signed diagram. If it is of the form $-E$, where $E$ is
  a signed diagram for the generators $\pi_2,\dots,\pi_{N-1}$, then
  using Remark~\ref{remark.fix_i},
  $$C_D = \pi^-_1 C_E \pi^-_1 = \pi^-_1 C_E\,.$$
  Similarly, if it is of the form $+E$, then:
  $$C_D = \pi^+_1 C_E \pi^+_1 = C_E \pi^+_1\,.$$
  Using induction on the isomorphic copy of $\NDPF_{N-1}$ generated by
  $\pi_2,\dots,\pi_{N-1}$ yields the desired formula.
\end{proof}

\begin{proposition}
  \label{proposition.idempotents.ndpf}
  The collection of all $C_D$ forms a complete set of orthogonal
  idempotents for $\NDPF_N$.
\end{proposition}
\begin{proof}
  First note that $C_D$ is never zero; for example, it is clear from
  Proposition~\ref{proposition.norton_ndpf} that the full expansion of
  $C_D$ has coefficient $1$ on $\prod_{i=n,\dots,1,\ i\in D} \pi^+_i$.

  Take now $D$ and $D'$ two signed diagrams. If they differ in the
  first position, it is clear that $C_D C_{D'} =0$. Otherwise, write
  $D = \epsilon E$, and $D' = \epsilon E'$. Then, using
  Remark~\ref{remark.fix_i} and induction,
  \begin{equation*}
    \begin{split}
      C_D C_D' & = \pi^\epsilon_1 C_E \pi^\epsilon_1 \pi^\epsilon_1 C_{E'} \pi^\epsilon_1
      = \pi^\epsilon_1 C_E \pi^\epsilon_1 C_{E'} \pi^\epsilon_1\\
      &= \pi^\epsilon_1 C_E C_{E'} \pi^\epsilon_1
      = \pi^\epsilon_1 \delta_{E,E'} C_E \pi^\epsilon_1
      = \delta_{D,D'} C_D\,.
    \end{split}
  \end{equation*}
  Therefore, the $C_D$'s form a collection of $2^{N-1}$ nonzero
  orthogonal idempotents, which has to be complete by cardinality.
\end{proof}

One can interpret the diagram demipotents for $\NDPF_N$ as branching
from the diagram demipotents for $\NDPF_{N-1}$ in the following way.
For any $C_D=L_DR_D$ in $\NDPF_{N-1}$, the leading term of $C_D$ will
be the longest element in the generators marked by plusses in $D$.
This leading idempotent has an image set which we will denote $\im(D)$
by abuse of notation. Now in $\NDPF_N$ we can associated two
`children' to $C_D$:
\[
C_{D+}=L_D\pi_{N}^+R_D \text{ and } C_{D-}=L_D\pi_{N}^-R_D.
\]
Then we have 
\[
C_{D+}+C_{D-}=C_D, \im(D+)=\im(D) \text{and} \im(D-)=\im(D)\bigcup \{N\}.
\]
\bigskip

We now generalize this branching construction to any meet semi-lattice
to derive a conjectural recursive formula for a decomposition of the
identity into orthogonal idempotents. This construction relies on the
branching rule for the idempotents of $\OR(L)$, and the existence of
the maximal idempotents $e_{a,b}$ of Remark~\ref{remark.oo.eab}.
\medskip

Let $L$ be a meet semi-lattice, and fix a linear extension of $L$. For
simplicity, we assume that the elements of $L$ are labelled
$1,\dots,N$ along this linear extension. Recall that, by
Proposition~\ref{proposition.idempotents.OO}, the idempotents are
indexed by the subsets of $L$ which contain the minimal elements of
$L$ and are stable under joins. In order to distinguish subsets of
$\{1,\dots,N\}$ and subsets of, say, $\{1,\dots,N-1\}$, even if they
have the same elements, it is convenient to identify them with $+-$
diagrams as we did for $\NDPF_N$. The \emph{valid diagrams} are those
corresponding to subsets which contain the minimal elements and are
stable under joins. A prefix of length $k$ of a valid diagram is still
a valid diagram (for $L$ restricted to $\{1,\dots,k\}$), and they are
therefore naturally organized in a binary prefix tree.

Let $D$ be a valid diagram, $e=\sup_D$ be the corresponding
idempotent. If $L$ is empty, $D=\{\}$, and we set
$L_{\{\}}=R_{\{\}}=1$.  Otherwise, let $L'$ be the meet semi-lattice
obtained by restriction of $L$ to $\{1,\dots,N-1\}$, and $D'$ the
restriction of $D$ to $\{1,\dots,N-1\}$.

\begin{itemize}
\item[Case 1] $N$ is the join of two elements of $\im(D')$ (and in particular,
  $N\in \im(D)$). Then, set $L_D=L_{D'}$ and $R_D=R_{D'}$.
\item[Case 2] $N\in \im(D)$. Then, set $L_D=L_{D'} \pi_{N, N.e}$ and
  $R_D=\pi_{N, N.e}R_{D'}$.
\item[Case 3] $N\not \in \im(D)$. Then, set $L_D=L_{D'} (1-\pi_{N, N.e})$
  and $R_D=(1-\pi_{N, N.e})R_{D'}$.
\end{itemize}
Finally, set $C_D=L_D R_D$.

\begin{remark}[Branching rule]
  Fix now $D'$ a valid diagram for $L'$. If $N$ is the join of two
  elements of $I'$, then $C_{D'}=C_{D'+}$. Otherwise $C_{D'}=C_{D'-} +
  C_{D'+}$.

  Hence, in the prefix tree of valid diagrams, the two sums of all
  $C_D$'s at depth $k$ and at depth $k+1$ respectively
  coincide. Branching recursively all the way down to the root of the
  prefix tree, it follows that the elements $C_D$ form a decomposition
  of the identity. Namely,
  \begin{equation*}
    1 = \sum_{D \text{ valid diagram}} C_D\,.
  \end{equation*}
\end{remark}

\begin{conjecture}
\label{conjecture.demi}
  Let $L$ be a meet semi-lattice. Then, the set $\{C_D\suchthat
  D \text{ valid diagram}\}$ forms a set of demipotent elements for
  $\OR(L)$ which, raised each to a sufficiently high power, yield a set
  of primitive orthogonal idempotents.
\end{conjecture}

This conjecture is supported by
Proposition~\ref{proposition.idempotents.ndpf}, as well as by computer
exploration on all $1377$ meet semi-lattices with at most $8$ elements
and on a set of meet semi-lattices of larger size which were considered
likely to be problematic by the authors. In all cases, the demipotents
were directly idempotents, which might suggest that
Conjecture~\ref{conjecture.demi} could be strengthened to state that
the collection $\{C_D\suchthat D \text{ valid diagram}\}$ forms
directly a set of primitive orthogonal idempotents for $\OR(L)$.

%% file: chapter-ndpfavoid.tex
In this chapter, we discuss some results relating the $\NDPF$ quotient of the $0$-Hecke monoid for the symmetric group to pattern avoidance results.  

In Section~\ref{sec:widthSystems} we introduce \textbf{width systems} on permutation patterns as a potential system for understanding pattern containment algebraically.  The main results of this section describe a class of permutation patterns $\sigma$ such that any permutation $x$ containing $\sigma$ factors as $x=y\sigma' z$, with $\len(x)=\len(y)+\len(\sigma)+\len(z)$.  Here $\sigma'$ is a ``shift'' of $\sigma$, and some significant restrictions on $y$ and $z$ are established.  
The main results are contained in Propositions~\ref{prop:s2widthSystems}, \ref{prop:s3widthSystems}, \ref{prop:widthSystemExtend},\ref{prop:widthSystemExtend2}, and Corollary~\ref{cor:bountifulPerms}.  

Pattern containment also has connection to the strong Bruhat order; in particular, Tenner showed that a principal order ideal of a permutation is Boolean if and only if the permutation avoids the patterns $[321]$ and $[3412]$~\cite{tenner.2007}.  

We apply these ideas directly in Section~\ref{sec:ndpfPattAvoid} while analyzing the fiber of a certain quotient of the $0$-Hecke monoid of the symmetric group.  In Theorem~\ref{thm:ndpfFibers231}, we show that each fiber of the quotient contains a unique $[321]$-avoiding permutation and a unique $[231]$-avoiding permutation.  We then apply an involution and study a slightly different quotient in which fibers contain a unique $[321]$-avoiding permutation and a unique $[312]$-avoiding permutation (Theorem~\ref{thm:ndpfFibers312}).  In Section~\ref{sec:bndpfPattAvoid}, we consider a different monoid-morphism of the $0$-Hecke monoid for which each fiber contains a unique $[4321]$-avoiding permutation (Theorem~\ref{thm:ndpfFibers4321}).

We then define the Affine Nondecreasing Parking Functions in Section~\ref{sec:affNdpfPattAvoid}, and establish these as a quotient of the $0$-Hecke monoid of the affine symmetric group.  We prove the existence of a unique $[321]$-avoiding affine permutation in each fiber of this quotient (Theorem~\ref{thm:affNdpfFibers321}).

\section{Background on Pattern Avoidance}

Pattern avoidance phenomena have been studied extensively, originally by Knuth in his 1973 classic, The Art of Computer Programming~\cite{knuth.TAOCP1}.  A thorough introduction to the subject may be found in the book ``Combinatorics of Permutations'' by Bona~\citetalias{bona.permutations}.  A \textbf{pattern} $\sigma$ is a permutation in $S_k$ for some $k$; given a permutation $x \in S_N$, we say that $x$ \textbf{contains the pattern $\sigma$} if, in the one-line notation for $x=[x_1, \ldots, x_N]$, there exists a subsequence $[x_{i_1}, \ldots, x_{i_k}]$ whose elements are in the same relative order as the elements in $p$.  If $x$ does not contain $\sigma$, then we say that 
$x$ \textbf{avoids $\sigma$}, or that $x$ is \textbf{$\sigma$-avoiding.}  (Note that if $k>N$, $x$ must avoid $\sigma$.)

For example, the pattern $[1, 2]$ appears in any $x$ such that there exists a $x_i < x_j$ for some $i<j$.  The only $[12]$-avoiding permutation in $S_N$, then, is the long element, which is strictly decreasing in one-line notation.  As a larger example, the permutation $[\mathbf{3}, \mathbf{4}, 5, \mathbf{2}, 1, 6]$ contains the pattern $[231]$ at the bold positions.  In fact, this permutation contains six distinct instances of the pattern $[231]$.

An interesting and natural question is, given a pattern $\sigma$, how many permutations in $S_N$ avoid $\sigma$?  It has been known since Knuth's original work that for any pattern in $S_3$, there are Catalan-many permutations in $S_N$ avoiding $\sigma$~\cite{knuth.TAOCP1}.

The $[321]$-avoiding permutations are of particular importance.  It was shown in~\cite{BilleyJockuschStanley.1993} that a permutation $x\in S_N$ is $[321]$-avoiding if and only if $x$ is `braid free.'  In particular, this means that there is no reduced word for $x$ containing the consecutive subsequence of $s_i s_{i+1} s_i$ (or $s_{i+1} s_i s_{i+1}$, equivalently), where the $s_i$ are the simple transpositions generating $S_N$.  Such permutations are called \textbf{fully commutative}. 

Lam~\cite{Lam06affinestanley} and Green~\cite{Green.2002} separately showed that this result extends to the  affine symmetric group.  The affine symmetric group (see Definition~\ref{def:affSn}) is a subset of the permutations of $\ZZ$, satisfying some periodicity conditions.  Pattern avoidance for the affine symmetric group works exactly as in a finite symmetric group.  
The one-line notation for $x$ is the doubly infinite sequence $x=[\ldots, x_{-1}, x_0, x_1, \ldots, x_N, x_{N+1}, \ldots]$.
Then $x$ contains a pattern $\sigma$ if any subsequence of $x$ in one-line notation has the same relative order as $\sigma$.  \textbf{Fully commutative elements} of the affine symmetric group are those which have no reduced word containing the consecutive subsequence $s_i s_{i+1} s_i$, where the indices are considered modulo $N$.
Green showed that the fully commutative elements of the affine symmetric group coincide with the $[321]$-avoiding affine permutations.

Fan and Green~\cite{fan.1996, fanGreen.1999} previously studied the quotient of the full Hecke algebra $H_q(W)$ for $W$ simply-laced, by the ideal $I$ generated by $T_{sts} + T_{st} + T_{ts} + T_s + T_t + 1$ for $s$ and $t$ generators of $W$ satisfying a braid relation $sts=tst$.  This quotient $H/I$ yields the \textbf{Temperley-Lieb Algebra}.  Fan showed that this quotient has a basis indexed by fully commutative elements of $W$, and in further work with Richard Green derived information relating this quotient to the Kazhdan-Lusztig basis for $H_q(W)$.

A further application of pattern avoidance occurs in the study of rational smoothness of Schubert varieties; an introduction to this topic may be found in~\cite{billeyLakshmibai.2000}.  The Schubert varieties $X_w$ in Type $A$ are indexed by permutations; a result of Billey~\cite{Billey98patternavoidance} shows that $X_w$ is smooth if and only if $w$ is simultaneously $[3412]$- and $[4231]$-avoiding.  More recently, Billey and Crites have extended this result to affine Schubert varieties (for affine Type A)~\cite{billeyCrites.2011}, showing that an affine Schubert variety $X_w$ is rationally smooth if and only if $w$ is simultaneously $[3412]$- and $[4231]$-avoiding or is a special kind of affine permutation, called a ``twisted spiral.''  

\section{Width Systems, Pattern Containment, and Factorizations.}
\label{sec:widthSystems}

In this section we introduce width systems on permutation patterns, which sometimes provide useful factorizations of a permutation containing a given pattern.  The results established here will be directly applied in Sections~\ref{sec:ndpfPattAvoid} and~\ref{sec:bndpfPattAvoid}. 

\begin{definition}
Let $x$ be a permutation and $\sigma \in S_k$ a pattern.  We say that $x$ \textbf{factorizes over $\sigma$} if there exist permutations $y$, $z$, and $\sigma'$ such that:
\begin{enumerate}
\item $x = y \sigma' z$,
\item $\sigma'$ has a reduced word matching a reduced word for $\sigma$ with indices shifted by some $j$,
\item The permutation $y$ satisfies $y^{-1}(j)<\cdots < y^{-1}(j+k)$, 
\item The permutation $z$ satisfies $z(j)<\cdots < z(j+k)$,
\item $\len(x) = \len(y) + \len(\sigma') + \len(z)$.
\end{enumerate}
\end{definition}

Set $W=S_N$ and $J\subset I$, with $I$ the generating set of $W$.  An element $x\in W$ has a \textbf{right descent} $i$ if $\len(x s_i)<\len(x)$, and has a \textbf{left descent} $i$ if $\len(s_i x)<\len(x)$.  Equivalently, $x$ has a right (resp., left) descent at $i$ if and only if some reduced word for $x$ ends (resp., begins) with $i$.  Let $W^J$ be the set of elements in $W$ with no right descents in $J$.  Similarly, $\leftexp{J}{W}$ consists of those elements with no left descents in $J$.  Finally, $W_J$ is the \textbf{parabolic subgroup} of $W$ generated by $\{ s_i \mid i \in J\}$.

Recall that a \textbf{reduced word} or \textbf{reduced expression} for a permutation $x$ is a minimal-length expression for $x$ as a product of the simple transpositions $s_i$.  Throughout this chapter, we will use double parentheses enclosing a sequence of indices to denote words.  For example, $((1,3,2))$ corresponds to the element  $s_1s_2s_3$ in $S_4$.  Note that same expression can also indicate an element of $H_0(S_4)$, with $((1,3,2))$ corresponding to the element $\pi_1\pi_2\pi_3$.  Context should make usage clear.

\begin{definition}
Let $\sigma$ be a permutation pattern in $S_k$, with reduced word $((i_1, \ldots, i_m))$.  Let $J=\{j, j+1, \ldots, j+l\}$ for some $l\geq k-1$ and $\sigma' \in W_J$ with reduced word $((i_1+j, \ldots, i_m+j))$.  Then we call $\sigma'$ a \textbf{$J$-shift} or \textbf{shift} of $\sigma$.
\end{definition}

\begin{figure}
  \includegraphics[scale=1]{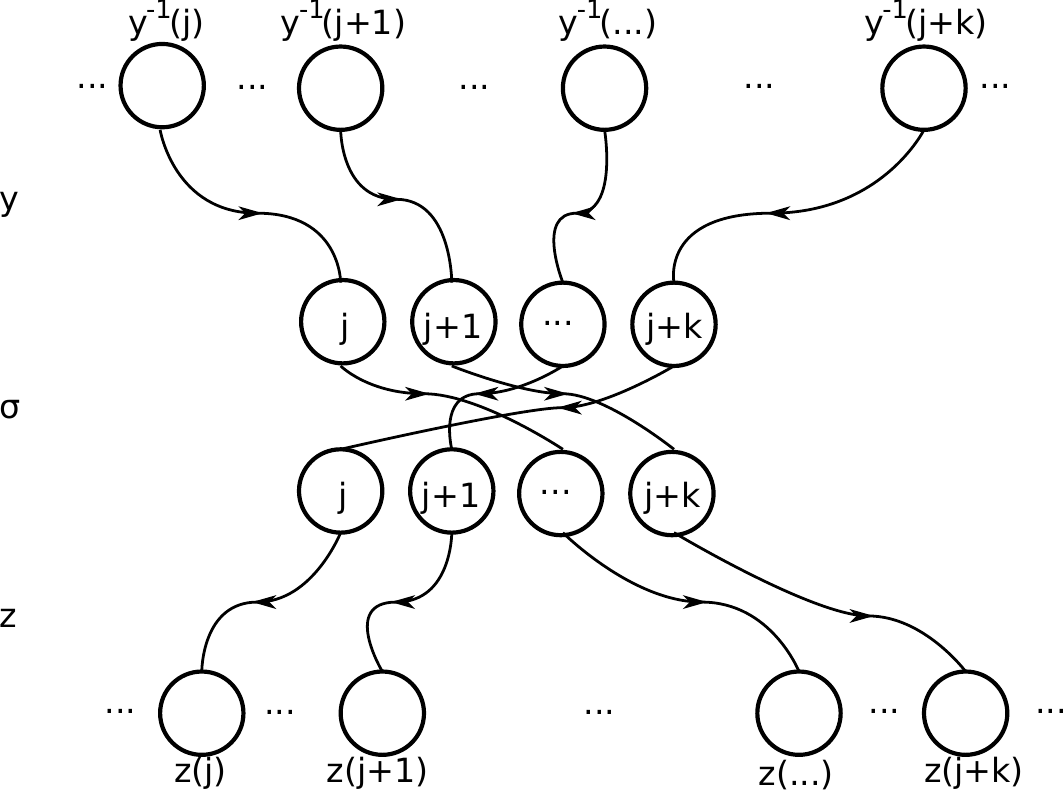}
  \caption{Diagrammatic representation of a permutation $x$ factorizing over a pattern $\sigma$ as $x=y\sigma z$ by composition of string diagrams. }
  \label{fig.patternContainment}
\end{figure}

\begin{proposition}
A permutation $x\in S_N$ factorizes over $\sigma$ if and only if $x$ admits a factorization $x=y\sigma' z$ with $y \in W^J, \sigma' \in W_J$, and $z\in \leftexp{J}{W}$, and $\len(x) = \len(y) + \len(\sigma') + \len(z)$.
\end{proposition}
\begin{proof}
This is simply a restatement of the definition of factorization over $\sigma$.  In particular, $y\in W^J$ and $z\in \leftexp{J}{W}$.
\end{proof}
This condition is illustrated diagrammatically in Figure~\ref{fig.patternContainment} using a string-diagram for the permutation $x$ factorized as $y \sigma' z$.  In the string diagram of a permutation $x$, a vertical string connects each $j$ to $x(j)$, with strings arranged so as to have as few crossings as possible.
Composition of permutations is accomplished by vertical concatenation of string diagrams.    In the diagram, $x$ is the vertical concatenation (and product of) of $y$, $\sigma'$ and $z$.

The permutation $y^{-1}$ preserves the order of $\{j, j+1, \ldots, j+k\}$, and thus the strings leading into the elements $\{j, j+1, \ldots, j+k\}$ do not cross.  Likewise, $z$ preserves the order of $\{j, j+1, \ldots, j+k\}$, and thus the strings leading out of $\{j, j+1, \ldots, j+k\}$ in $z$ do not cross.  In between, $\sigma'$ rearranges $\{j, j+1, \ldots, j+k\}$ according to the pattern $\sigma$.

By the above discussion, it is clear that if $x$ admits a factorization $y\sigma' z$ with $y \in W^J, \sigma' \in W_J$, and $z\in \leftexp{J}{W}$ then $x$ contains $\sigma$.  The question, then, is when this condition is sharp.  This question is interesting because it provides an algebraic description of pattern containment.  For example, a permutation $x$ which contains a $[321]$-pattern is guaranteed to have a reduced expression which contains a braid.  Braid containment can be re-stated as a factorization over $[321]$.  When the factorization question is sharp, (ie, $x$ contains $\sigma$ if and only if $x$ factorizes over $\sigma$) one obtains an algebraic description of $\sigma$-containment.  The class of patterns with this property is rather larger than just $[321]$, as we will see in Propositions~\ref{prop:s2widthSystems}, \ref{prop:s3widthSystems}, and~\ref{prop:widthSystemExtend}.

\begin{problem}
For which patterns $\sigma$ does $x$ contain $\sigma$ if and only if $x\in W^J \sigma' \leftexp{J}{W}$, where $\sigma'$ is a $J$-shift of $\sigma$ for some $J$?
\end{problem}

As a tool for attacking this problem, we introduce the notion of a width system for a pattern.

\begin{definition}
Suppose $x$ contains $\sigma$ at positions $(i_1, \ldots, i_k)$; the tuple $P=(P_1, \ldots, P_k)$ is called an \textbf{instance} of the pattern $\sigma$, and we denote the set of all instances of $\sigma$ in $x$ by $P_x$.
\end{definition}

\begin{definition}
A \textbf{width} on an instance $P$ of $\sigma$ is a difference $P_j-P_i$ with $j>i$.  
A \textbf{width system} $w$ for a permutation pattern $\sigma \in S_k$ is a function assigning a tuple of widths to each instance of $\sigma$ in $x$.  An instance $P$ of a pattern in $x$ is \textbf{minimal} (with respect to $\sigma$ and $w$) if $w(P)$ is lexicographically minimal amongst all instances of $\sigma$ in $x$.  Finally, an instance $P=(P_1, \ldots, P_k)$ is \textbf{locally minimal} if $P$ is the minimal instance of $\sigma$ in the partial permutation $[x_{P_1}, x_{P_1+1}, \ldots, x_{P_k-1}, x_{P_k}]$.
\end{definition}

\begin{example}
Consider the pattern $[231]$ and let $P=(p, q, r)$ be an arbitrary instance of $\sigma$ in a permutation $x$.  We choose to consider the width system $w(P)=(r-p, q-p)$.  (Other width systems include $u(P)=(r-q, q-p)$ and $v(P)=(r-q)$, for example.)

The permutation $x=[3, 4, 5, 2, 1, 6]$ contains six $[231]$ patterns.  The following table records each $[231]$-instance $P$ and the width of the instance $w(P)$:
\begin{equation*}
\begin{array}[b]{|c c c|}
 \hline           &         P     &  w(P)  \\ \hline
 \left[\mathbf{3}, \mathbf{4}, 5, \mathbf{2}, 1, 6\right]       & (1, 2, 4) & (3, 1) \\
 \left[\mathbf{3}, \mathbf{4}, 5, 2, \mathbf{1}, 6\right]       & (1, 2, 5) & (4, 1) \\
 \left[\mathbf{3}, 4, \mathbf{5}, \mathbf{2}, 1, 6\right]       & (1, 3, 4) & (3, 2) \\
 \left[\mathbf{3}, 4, \mathbf{5}, 2, \mathbf{1}, 6\right]       & (1, 3, 5) & (4, 2) \\
 \left[3, \mathbf{4}, \mathbf{5}, \mathbf{2}, 1, 6\right]       & (2, 3, 4) & (2, 1) \\
 \left[3, \mathbf{4}, \mathbf{5}, 2, \mathbf{1}, 6\right]       & (2, 3, 5) & (3, 1) \\  \hline
\end{array}
\end{equation*}
Thus, under the width system $w$ the instance $(2, 3, 4)$ is the minimal $[231]$-instance; it is also the only locally minimal $[231]$-instance.

In the permutation $y=[1, 4, 8, 5, 2, 7, 6, 3]$, we have the following instances and widths of the pattern $[231]$:
\begin{equation*}
\begin{array}[b]{|c c c|}
  \hline          &         P     &  w(P)  \\ \hline
 \left[ 1, \mathbf{4}, \mathbf{8}, 5, \mathbf{2}, 7, 6, 3 \right]  & (2, 3, 5) & (3, 1) \\
 \left[ 1, \mathbf{4}, \mathbf{8}, 5, 2, 7, 6, \mathbf{3} \right]  & (2, 3, 8) & (6, 1) \\
 \left[ 1, \mathbf{4}, 8, \mathbf{5}, \mathbf{2}, 7, 6, 3 \right]  & (2, 4, 5) & (3, 2) \\
 \left[ 1, \mathbf{4}, 8, \mathbf{5}, 2, 7, 6, \mathbf{3} \right]  & (2, 4, 8) & (6, 2) \\
 \left[ 1, \mathbf{4}, 8, 5, 2, \mathbf{7}, 6, \mathbf{3} \right]  & (2, 6, 8) & (6, 4) \\
 \left[ 1, \mathbf{4}, 8, 5, 2, 7, \mathbf{6}, \mathbf{3} \right]  & (2, 7, 8) & (6, 5) \\
 \left[ 1, 4, 8, \mathbf{5}, 2, \mathbf{7}, 6, \mathbf{3} \right]  & (4, 6, 8) & (4, 2) \\
 \left[ 1, 4, 8, \mathbf{5}, 2, 7, \mathbf{6}, \mathbf{3} \right]  & (4, 7, 8) & (4, 3) \\  \hline
\end{array}
\end{equation*}
Here, the instance $(2, 3, 5)$ is minimal under $w$.  Additionally, the instance $(4, 6, 8)$ is locally minimal, since it is the minimal instance of $[231]$ in the partial permutation 
$\left[ \mathbf{5}, 2, \mathbf{7}, 6, \mathbf{3} \right]$.
\end{example}

For certain width systems, minimality provides a natural factorization of $x$ over $\sigma$.

\begin{example}
\label{ex:231bountiful}
\begin{figure}
  \includegraphics[scale=1]{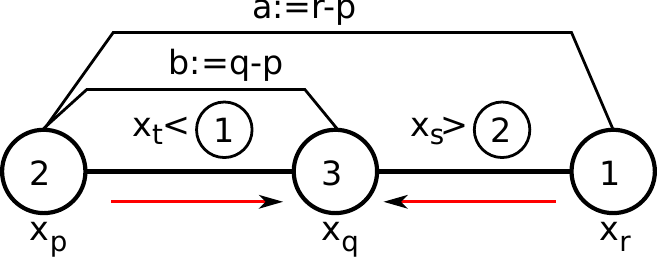}
  \caption{A diagram of a minimal $[231]$ pattern.  The circled numbers represent elements $(x_p, x_q, x_r)$ filling the roles of the pattern; the widths are denoted $a$ and $b$, and the restrictions on $x_t$ with $p<t<q$ and $x_s$ with $s<q<r$ implied by minimality of the pair $(a,b)$ are also recorded.  The red arrows record the fact that shifting the end elements towards the center using a sequence of simple transpositions reduces the length of the permutation.}
  \label{fig.minimal231}
\end{figure}
We consider the width system for the pattern $[231]$ depicted in Figure~\ref{fig.minimal231}.

Let $x=[x_1,x_2,\ldots,x_N] \in S_N$ containing a $[231]$-pattern, and let $(p<q<r)$ be the indices of a minimal-width $[231]$-pattern in $x$ under the width system $w=(r-p, q-p)$.  (So $x_r<x_p<x_q$.)

Minimality of the total width $(r-p)$ implies that for every $s$ with $q<s<r$, we have $x_s>x_p(>x_r)$, as otherwise $(x_p, x_q, x_s)$ would be a $[231]$-pattern of smaller width.  Then multiplying $x$ on the right by $u_1=s_{r-1}s_{r-2}\ldots s_{q+1}$ yields a permutation of length $\len(x)-(r-q-1)$, with 
\[
xu_1=[x_1,\ldots,x_p,\ldots,x_q,x_r,x_{q+1}\ldots,x_N].
\]

Minimality of the inner width $(q-p)$ implies that for every $t$ with $p<t<q$, then $x_t<x_r$.  (If $x_p<x_t<x_q$, then $(x_t, x_q, x_r)$ would form a $[231]$-pattern of lower width.  If $x_p>x_t$, then $q$ was not chosen minimally.)  Then multiplying $xu_1$ on the right by $u_2=s_ps_{p+1}\ldots s_{q-2}$ yields a permutation of length $\len(xu_1)-(q-p-1) = \len(x) - r + p + 2)$.  This permutation is:
\[
xu_1u_2=[x_1,\ldots,x_{q-1},x_p,x_q,x_r,x_{q+1}\ldots,x_N].
\]

Since $[x_p, x_q, x_r]$ form a $[231]$-pattern, we may further reduce the length of this permutation by multiplying on the right by $s_{q}s_{q-1}$.  The resulting permutation has no right descents in the set $J:=\{q-1, q\}$.

We then set $y=xu_1u_2s_{q}s_{q-1}$, $\sigma'=s_{q-1}s_{q}$, and $z=(u_1u_2)^{-1}$.  Notice that $z$ has no left descents in $\{q-1, q\}$ by construction, since it preserved the left-to-right order of $x_p, x_q$ and $x_r$.   Then $x=y \sigma' z$ is a factorization of $x$ over $\sigma$.
\end{example}

\begin{figure}
  \includegraphics[scale=1]{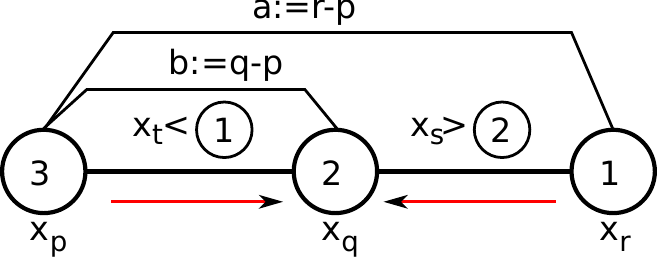}
  \caption{A diagram of a left-minimal $[321]$ pattern, labeled analogously to the labeling in Figure~\ref{fig.minimal231}.}
  \label{fig.minimal321}
\end{figure}

One may use a similar system of minimal widths to show that any permutation containing a $[321]$-pattern contains a braid, replicating a result of Billey, Jockusch, and Stanley~\cite{BilleyJockuschStanley.1993}.  The corresponding system of widths is depicted in Figure~\ref{fig.minimal321}.

\begin{definition}
Let $\sigma$ be a permutation with a width system.  The width system is \textbf{bountiful} if for any $x$ containing a locally minimal $\sigma$ at positions $(p_1, \ldots, p_k)$, any $x_t$ with $p_i<t<p_{i+1}$ has either $x_t<x_{p_k}$ for all $p_k<t$ or $x_t>x_{p_k}$ for all $p_k>t$.
\end{definition}

\begin{proposition}
If a pattern $\sigma$ admits a bountiful width system, then any $x$ containing $\sigma$ factorizes over $\sigma$.
\end{proposition}
\begin{proof}
By definition, any $x_t$ with $p_i<t<p_{i+1}$ has either $x_t<x_{p_k}$ for all $p_k<t$ or $x_t>x_{p_k}$ for all $p_k>t$.  Then using methods exactly as in Example~\ref{ex:231bountiful}, we may vacate the elements $x_t$ by multiplying on the right by simple transpositions, moving ``small'' $x_t$ out to the left and moving ``large'' $x_t$ out to the right.  This brings the minimal instance of the pattern $\sigma$ together into adjacent positions $(j, j+1, \ldots, j+k)$, while simultaneously creating a reduced word for the right factor $z$ in the factorization.  Then we set $J=\{j, j+1, \ldots, j+k-1\}$, and let $\sigma'$ be the $J$-shift of $\sigma$.  Set $y=x z^{-1} \sigma'^{-1}$.  Then by construction $x=y \sigma' z$ is a factorization of $x$ over $\sigma$.
\end{proof}

Thus, establishing bountiful width systems allows the direct factorization of $x$ containing $\sigma$ as an element of $W^J \sigma' \leftexp{J}{W}$.

\begin{problem}
Characterize the patterns which admit bountiful width systems.
\end{problem}

\begin{example}
\label{ex:123broken}
The permutation $x = [1324] = s_2$ contains a $[123]$-pattern, but does not factor over $[123]$.  To factor over $[123]$, we have $x\in W^J 1_J \leftexp{J}{W}$, with $J=\{1,2\}$ or $J=\{2,3\}$.  Both choices for $J$ contain $2$, so it is impossible to write $x$ as such a product.
\end{example}

\begin{proposition}
\label{prop:s2widthSystems}
Both patterns in $S_2$ admit bountiful width systems.
\end{proposition}
\begin{proof}
Any minimal $[12]$- or $[21]$-pattern must be adjacent, and so the conditions for a bountiful width system hold vacuously.
\end{proof}

\begin{figure}
  \includegraphics[scale=1]{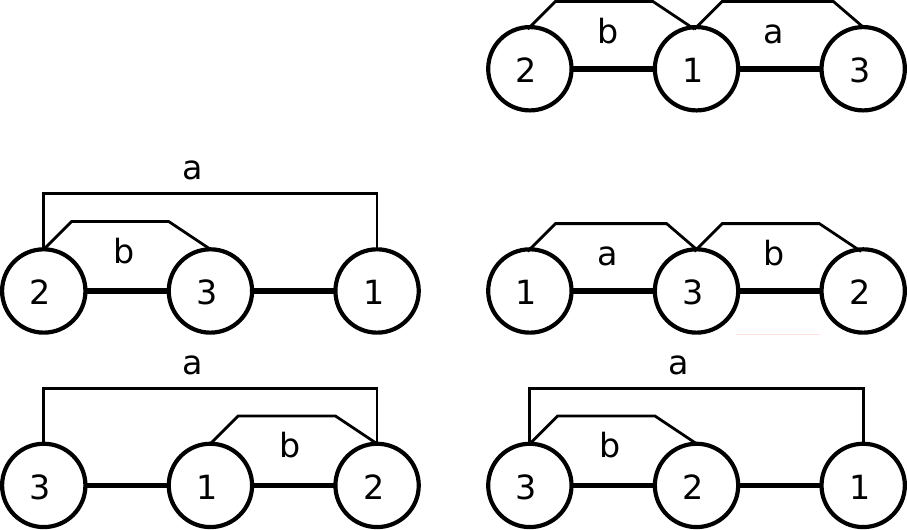}
  \caption{Diagrams of bountiful width systems for the five patterns in $S_3$ which admit bountiful width systems.}
  \label{fig.s3widthSystems}
\end{figure}

\begin{proposition}
\label{prop:s3widthSystems}
All of the patterns in $S_3$ except $[123]$ admit a bountiful width system, as depicted in Figure~\ref{fig.s3widthSystems}.
\end{proposition}

\begin{proof}
A bountiful width systems has already been provided for the pattern $[231]$.  We only provide the details of the proof that the $[213]$ pattern is bountiful, as the proofs that the width systems for the patterns $[132]$, $[312]$ and $[321]$ are bountiful are analogous.

Let $x\in S_N$ contain a $[213]$ pattern at positions $(x_p, x_q, x_r)$, and choose the width system $(a,b)=(r-q, q-p)$.  

Suppose that $(x_p, x_q, x_r)$ is lexicographically minimal in this width system, and consider $x_t$ with $p<t<q$ and $x_s$ with $q<s<r$. Then $a=1$:
\begin{itemize}
\item If $x_s<x_p$, then $(x_p, x_s, x_r)$ is a $[213]$ pattern with $a$ smaller.
\item If $x_p<x_s$, then $(x_p, x_q, x_s)$ is a $[213]$ pattern with $a$ smaller.
\end{itemize}
Thus, we must have $r-q=1$.

Since $b$ is minimal, we must also have that $x_t>x_q$ or $x_t<x_p$ for every $t$ with $p<t<q$.  This completes the proof that the width system is bountiful.
\end{proof}

\begin{proposition}
\label{prop:widthSystemExtend}
Let $\sigma$ be a pattern in $S_{K-1}$ with a bountiful width system, and let $\sigma_+ = [K, \sigma_1, \ldots, \sigma_{K-1}]$. Then $\sigma_+$ admits a bountiful width system.  

Similarly, let $\sigma_-= [\sigma_1+1, \ldots, \sigma_{K-1}+1, 1]$.  Then $\sigma_-$ admits a bountiful width system.
\end{proposition}
\begin{proof}
Let $w=(w_1, w_2, \ldots, w_{k-2})$ be a bountiful width system on $\sigma$ (so $w_i$ is the difference between indices of an instance of $\sigma$ in a given permutation).  Let $x$ contain $\sigma_+$ in positions 
$( x_p, \ldots, x_q )$.  For $\sigma_+$, we show that the width system 
$w_+=(w_1, w_2, \ldots, w_{k-2}, q-p)$ is bountiful, where $w_i$ measures widths of elements in $\sigma$ as in $w$.

Consider a $\sigma_+$-pattern in a permutation $x$ that is minimal under the width system $w_+$, appearing at indices given by the tuple $p:=(i_1, \ldots, i_{k+1})$.  Then $x$ contains a $\sigma$-pattern at positions $(i_2, \ldots, i_{k+1})$.  This pattern may not be minimal under $w$ but, by the choice of width system, is as close as possible to being $w$-minimal, in the following sense.

We examine two cases.
\begin{itemize}
\item If there are no indices $t$ with $i_2<t<i_{k+1}$ such that $x_t>x_{i_1}$, then $\sigma$ must be $w$-minimal on the range $i_2, \ldots, i_{k+1}$.  (Otherwise, a $w$-minimal $\sigma$-pattern in that space would extend to a pattern that was less than $p$ in the $w_+$ width system.)  Then bountifulness of the $\sigma$ pattern ensures that for any $t$ with $i_j<t<i_{j+1}$ with $j\geq 2$; then $x_t<x_{i_k}$ for all $i_k<t$ or $x_t>x_{i_k}$ for all $i_k>t$.  (The ``small'' elements are still smaller than the ``large'' element $x_{i_1}$.)



\item On the other hand, if there exist some $t$ with $i_2<t<i_{k+1}$ such that $x_t>x_{i_1}$, we may move these $x_t$ out of the $\sigma$ pattern to the right by a sequence of simple transpositions, each decreasing the length of the permutation by one.  Let $u$ be the product of this sequence of simple transpositions.  Then $xu$ fulfills the previous case.  Each of the $x_t$ were larger than all pattern elements to the right, so we see that $\sigma_+$ fulfills the requirements of a bountiful pattern.

\end{itemize}

The proof that $\sigma_-$ admits a bountiful width system is similar.
\end{proof}

\begin{corollary}
\label{cor:bountifulPerms}
Let $\sigma \in S_K$ be a permutation pattern, where the length of $\sigma$ is at most one less than the length of the long element in $S_K$.  Then $\sigma$ admits a bountiful width system.
\end{corollary}
\begin{proof}
This follows inductively from Proposition~\ref{prop:widthSystemExtend}, and the fact that the patterns $[12]$ and $[21]$ both admit bountiful width systems.
\end{proof}
\begin{figure}
  \includegraphics[scale=1]{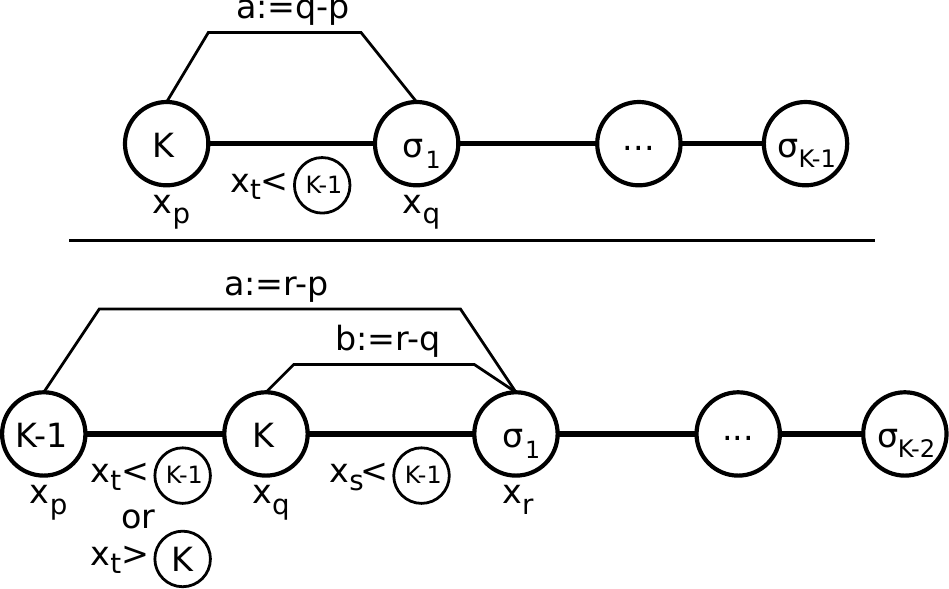}
  \caption{Diagram of extensions of a bountiful width system $w$ by the additional widths $a$ or $(a, b)$, as described in the proofs of Propositions~\ref{prop:widthSystemExtend} and~\ref{prop:widthSystemExtend2}. }
  \label{fig.widthSystemExtend}
\end{figure}

\begin{proposition}
\label{prop:widthSystemExtend2}
Let $\sigma$ be a pattern in $S_{K-2}$ with a bountiful width system, and let $\sigma_{++} = [K-1, K, \sigma_1, \ldots, \sigma_{K-1}]$. Then $\sigma_{++}$ admits a bountiful width system.  

Similarly, let $\sigma_{--}= [\sigma_1+2, \ldots, \sigma_{K-2}+2, 1, 2]$.  Then $\sigma_{--}$ admits a bountiful width system.
\end{proposition}
\begin{proof}
The proof of this proposition closely mirrors the proof of Proposition~\ref{prop:widthSystemExtend}.  Let $w=(w_1, w_2, \ldots, w_{k-2})$ a bountiful width system on $\sigma$.  Let $x$ contain $\sigma_{++}$ in positions 
$( x_p, x_r, x_s, \ldots, x_q )$.  For $\sigma_{++}$, we claim that the width system 
$w_{++}=(w_1, w_2, \ldots, w_{k-2}, q-p, s-r)$ is bountiful, where $w_i$ measures widths of elements in $\sigma$ as in $w$.  (The width system $w_{++}$ is depicted in Figure~\ref{fig.widthSystemExtend}.)

Again, local minimality of $\sigma$ ensures that all $x_t$ with $s<t<q$ with $x_t$ not in the instance of $\sigma_{++}$ are either smaller than all pattern elements to the left of $x_t$, or larger than all pattern elements to the right of $x_t$.  The choice of $w_{++}$ ensures that all $x_t$ with $p<t<r$ are either less than $x_p$ or larger than $x_r$, and that all $x_t$ with $r<t<s$ are less than $x_p$.  Then $w_{++}$ is bountiful.

The proof that $\sigma_{--}$ is bountiful is analogous.
\end{proof}

\subsection{Further Directions}
Preliminary investigation suggests that patterns admitting a bountiful width system are somewhat rare, though there are more than those described by Corollary~\ref{cor:bountifulPerms}.  Weakening the definition of a factorization over a permutation may provide an additional avenue of investigation, though.

\begin{definition}
A permutation $x\in W=S_N$ \textbf{left-factorizes} over a pattern $\sigma\in S_K$ if $x=y\sigma' z$ with:
\begin{itemize}
\item $\sigma' \in W_J$, with $J=\{j, j+1, \ldots, j+k\}$ and $\sigma'$ containing a $\sigma$-pattern,
\item $y \in W^J$,
\item $\len(x)=\len(y) + \len(\sigma) + \len(z)$.
\end{itemize}
\end{definition}

This definition drops the requirement that $z \in \leftexp{J}{W}$.  This definition may be too weak, though, since one can show that any permutation containing the pattern $[K, K-1, \ldots, 1]$ left-factors over every pattern in $S_K$.

On the other hand, consider Example~\ref{ex:123broken}.  The permutation $x=[1, 3, 2, 4]=s_2$ admits a factorization $S^{\{1,3\}} 1_{\{1,3\}} \leftexp{\{1,3\}}{S}$, and the element $1_{\{1,3\}}$ contains a $[123]$-pattern.  Allowing factorizations over arbitrary subgroups -- and obtaining a combinatorial characterization of these factorizations -- may provide a way forward.

\begin{problem}
Find a general characterization of pattern containment in terms of factorizations of a permutation.
\end{problem}

\section{Pattern Avoidance and the $\NDPF$ Quotient}
\label{sec:ndpfPattAvoid}

In this section, we consider certain quotients of the $0$-Hecke monoid of the symmetric group, and relate the fibers of the quotient to pattern-avoidance.  The $0$-Hecke monoid $H_0(S_N)$ is defined in Definition~\ref{ssec:zeroHeckeDefinition}, and the Non-decreasing Parking Function $\NDPF_N$ quotient is discussed in Section~\ref{ssec:ropf}, in its guise as the the monoid of order-preserving regressive functions on a chain.

\begin{definition}
For $x\in H_0(S_N)$, we say $x$ {\bf contains a braid} if some reduced word for $x$ contains a contiguous subword $\pi_i \pi_{i+1} \pi_i$.  

The permutation $x$ contains an {\bf unmatched ascent} if some reduced word for $x$ contains a contiguous subword $\pi_i \pi_{i+1}$ that is not part of a braid.  More precisely, if inserting a $\pi_i$ directly after the $\pi_i \pi_{i+1}$ increases the length of $x$, then $x$ contains an unmatched ascent.  Equivalently, $x$ may be factorized as $x=y \pi_i \pi_{i+1} z$, where $y$ has no right descents in $\{i, i+1\}$, and $z$ has no left descents in $\{i, i+1\}$, and $\len(x)=\len(y)+2+\len(z)$.

An {\bf unmatched descent} is analogously defined as a contiguous subword $\pi_{i+1}\pi_i$ such that insertion of a $\pi_i$ immediately before this subword increases the length of $x$.  Equivalently, $x$ may be factorized as $x=y \pi_{i+1}\pi_i z$, where $y$ has no right descents in $\{i, i+1\}$, and $z$ has no left descents in $\{i, i+1\}$, and $\len(x)=\len(y)+2+\len(z)$.
\end{definition}

\begin{lemma}
\label{lemma:231unmatched}
For $x\in S_N$, $x$ contains a $[231]$-pattern if and only if $x$ has an unmatched ascent.  Likewise, $x$ contains a $[312]$-pattern if and only if $x$ has an unmatched descent.
\end{lemma}
\begin{proof}
This is a straightforward application of the bountiful width system for the patterns $[231]$ and $[312]$.  The resulting factorization contains an unmatched ascent (resp., descent).
\end{proof}

This process of inserting an $s_i$ can be made more precise in the symmetric group setting: suppose $s_{j_1}\ldots s_{i}s_{i+1}\ldots s_{j_k}$ is a reduced expression for $x \in S_N$.  Then write $x=x_1 s_{i}s_{i+1} x_2$.  To insert $s_i$, multiply $x$ on the right by $x_2^{-1}s_ix_2$.  As such, this insertion can be realized as multiplication by some reflection.

This insertion is generally not a valid operation in $H_0(S_N)$, since inverses do not exist.  However, the operation does make sense in the $\NDPF$ setting: the $\NDPF$ relation simply allows one to exchange a braid for an unmatched ascent or vice-versa.

\begin{theorem}
\label{thm:ndpfFibers231}
Each fiber of the map $\phi: H_0(S_N)\to \NDPF_N$ contains a unique $[321]$-avoiding element of minimal length and a unique $[231]$-avoiding element of maximal length.
\end{theorem}

\begin{proof}
The first part of the theorem follows directly from a result of Billey, Jockusch, and Stanley~\cite{BilleyJockuschStanley.1993}, which states that a symmetric group element contains a braid if and only if the corresponding permutation contains a $[321]$.  Alternatively, one can use the width system for $[321]$ established in Proposition~\ref{prop:s3widthSystems} to obtain a factorization including a braid.  Then for any $x$ in the fiber of $\phi$, one can remove braids obtained from minimal-width $[321]$-patterns using the $\NDPF$ relation and obtain a $[321]$-avoiding element.  Each application of the $\NDPF$-relation reduces the length of the permutation by one, so this process must eventually terminate in a $[321]$-avoiding element.  Uniqueness follows since there are exactly $C_N$ $[321]$-avoiding elements in $S_N$, where $C_N$ is the $N$th Catalan number, and are thus in bijection with elements of $\NDPF_N$.

For the second part, we use the bountiful $[231]$ width system established in Example~\ref{ex:231bountiful}.  Let $x$ contain a $[231]$-pattern.  The width system allows us to write a factorization $x=y \pi_i \pi_{i+1} z$, where $y$ has no right descents in $\{i, i+1\}$ and $z$ has no left descents in $\{i, i+1\}$.  Then we may apply the $\NDPF$ relation to insert a $\pi_i$, turning the $[231]$-pattern into a $[321]$ pattern, and increasing the length of $x$ by one.  Since we are in a finite symmetric group, there is an upper bound on the length one may obtain by this process, and so the process must terminate with a $[231]$-avoiding element.  Recall that $[231]$-avoiding permutations are also counted by the Catalan numbers~\cite{knuth.TAOCP1}, and apply the same reasoning as above to complete the theorem.
\end{proof}

Recall that the \emph{right action} of $S_N$ acts on positions.  A permutation $y$ has a \emph{right descent at position i} if the two consecutive elements $y_i,y_{i+1}$ are out of order in one-line notation.  Then multiplying on the right by $s_i$ puts these two positions back in order and reduces the length of $y$ by one.  Likewise, if $y$ does not have a right descent at $i$, multiplying by $s_i$ increases the length by one.

\begin{example}[Fibers of the $\NDPF$ quotient]
For $S_4$, the fibers of the $\NDPF$ quotient can be found in Figure~\ref{fig.s4ndpfFibers}.  

\begin{figure}
  \includegraphics[scale=1]{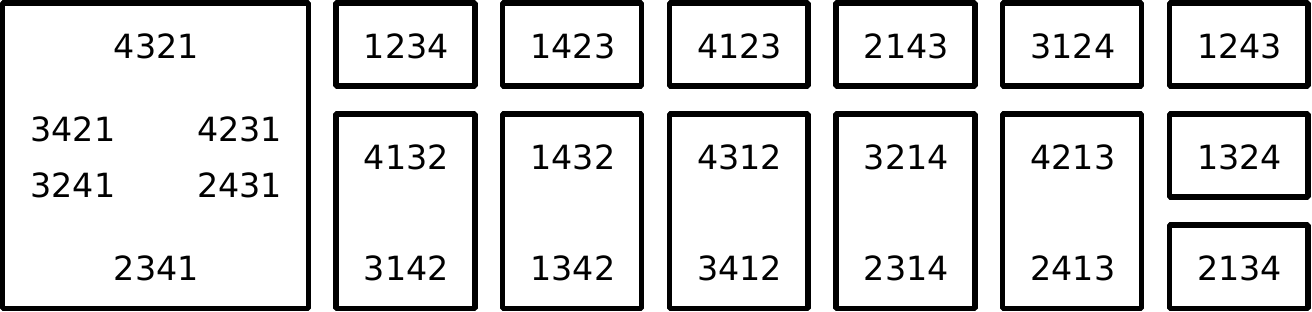}
  \caption{Fibers of the $\NDPF$ quotient for $H_0(S_4)$.}
  \label{fig.s4ndpfFibers}
\end{figure}
As a larger example, let $\sigma=[3,6,4,5,7,2,1]\in S_7$.  For Lemma~\ref{lemma:231unmatched}, we find minimal-width $[231]$-patterns, with the element corresponding to the $3$ chosen as far to the left as possible.  (The subsequence $(5, 7, 2)$ of $\sigma$ is such a minimal $[231]$-pattern.)  Then applying the transformation $[231]\rightarrow [321]$ on that instance of the pattern preserves the fiber of the $\NDPF$ quotient, and increases the length of the permutation by $1$.  By sequentially removing eight such minimal $[231]$-patterns, one obtains the long element in $S_7$, which is $[231]$-avoiding.  The fiber containing the long element also contains a $[321]$-avoiding element $[2,3,4,5,6,7,1]$, which has length $6$, and is the shortest element in its fiber.
\end{example}

We now fix bountiful width system for $[231]$- and $[321]$-patterns, which we will use for the remainder of this section.
\begin{definition}
Let $x\in S_N$, $x=[x_1, \ldots, x_N]$ in one-line notation, and consider all $[231]$-patterns $(x_p,x_q,x_r)$ in $x$.  The \textbf{width} of a $[231]$-pattern $(x_p,x_q,x_r)$ is the pair $(r-p, q-p)$.  The pattern is a \textbf{minimally chosen $[231]$-pattern} if the width is lexicographically minimal amongst all $[231]$-patterns in $x$.

On the other hand, call a $[321]$-pattern $(x_p,x_q,x_r)$ \textbf{left minimal} if for all $t$ with $p<t<q$, $x_t<x_r$, and for all $s$ with $q<s<r$, $x_s>x_q$.
\end{definition}

The following is a direct result of the proof of Lemma~\ref{lemma:231unmatched}.
\begin{corollary}
\label{cor:fiber231}
Let $x\in S_N$.  Let $(x_p,x_q,x_r)$ be a minimally chosen $[231]$-pattern in $x$.  Then the permutation 
\[
[x_1, \ldots, x_{p-1}, x_q, x_{p+1}, \ldots, x_{q-1}, x_p, x_{q+1}, \ldots, x_r, \ldots, x_N],
\]
obtained by applying the transposition $t_{p,q}$, is in the same $\NDPF$-fiber as $x$.  The result of applying this transposition is a left-minimal $[321]$-pattern.
\end{corollary}

\subsection{Involution}
\label{subsec:involution}

Let $\Psi$ be the involution on the symmetric group induced by conjugation by the longest word.  Then $\Psi$ acts on the generators by sending $s_i \to s_{N-i}$.  This descends to an isomorphism of $H_0(S_N)$ by exchanging the generators in the same way: $\pi_i \to \pi_{N-i}$.

We can thus obtain a second map from $H_0(S_N)\to \NDPF_N$ by pre-composing with $\Psi$.  This has the effect of changing the $\NDPF$ relation to a statement about unmatched \emph{descents} instead of unmatched ascents.  Then applying the $\NDPF$ relation allows one to exchange braids for unmatched descents and vice-versa, giving the following theorem.

\begin{theorem}
\label{thm:ndpfFibers312}
Each fiber of the map $\phi \circ \Psi: H_0(S_N)\to \NDPF_N$ contains a unique $[321]$-avoiding element for minimal length and a unique $[312]$-avoiding element of maximal length.
\end{theorem}

The proof is exactly the mirror of the proof in previous section.

We fix bountiful width system for $[312]$-patterns, and a second bountiful width system for $[321]$-patterns, which we will use for the remainder of this section.
\begin{definition}
Let $x\in S_N$, $x=[x_1, \ldots, x_N]$ in one-line notation, and consider all $[312]$-patterns $(x_p,x_q,x_r)$ in $x$.  The \textbf{width} of a $[312]$-pattern $(x_p,x_q,x_r)$ is the pair $(r-p, r-q)$.  The pattern is a \textbf{minimally chosen $[312]$-pattern} if the width is lexicographically minimal amongst all $[312]$-patterns in $x$.

Likewise, call a $[321]$-pattern $(x_p,x_q,x_r)$ \textbf{right minimal} if the \textbf{right width} $(p-r, r-q)$ is lexicographically minimal amongst all $[321]$-patterns in $x$.
On the other hand, call a $[321]$-pattern $(x_p,x_q,x_r)$ \textbf{right minimal} if for all $t$ with $p<t<q$, $x_t<x_q$, and for all $s$ with $q<s<r$, $x_s>x_p$.
\end{definition}

\begin{corollary}
\label{cor:fiber312}
Let $x\in S_N$.  Let $(x_p,x_q,x_r)$ be a minimally chosen $[312]$-pattern in $x$.  Then the permutation 
\[
[x_1, \ldots, x_{p-1}, x_q, x_{p+1}, \ldots, x_{q-1}, x_p, x_{q+1}, \ldots, x_r, \ldots, x_N],
\]
obtained by applying the transposition $t_{p,q}$, is in the same $\NDPF\circ \Psi$-fiber as $x$.  The result of applying this transposition is a right-minimal $[321]$-pattern.
\end{corollary}

\section{Type B $\NDPF$ and $[4321]$-Avoidance}
\label{sec:bndpfPattAvoid}

In this section, we establish a monoid morphism of $H_0(S_N)$ whose fibers each contain a unique $[4321]$-avoiding permutation.  To motivate this map, we begin with a discussion of Non-Decreasing Parking Functions of Type $B$.

The Weyl Group of Type $B$ may be identified with the \textbf{signed symmetric group} $S_N^B$, which is discussed (for example) in~\cite{Bjorner_Brenti.2005}.  Combinatorially, $S_N^B$ may be understood as a group permuting a collection of $N$ labeled coins, each of which can be flipped to heads or tails.  The size of $S_N^B$ is thus $2^NN!$.   A minimal set of generators of this group are exactly the simple transpositions $\{t_i\mid i \in \{1, \ldots, N-1\}\}$ interchanging the coins labeled $i$ and $i+1$, along with an extra generator $t_N$ which flips the last coin.  

The group $S_N^B$ can be embedded into $S_{2N}$ by identifying the $t_i$ with $s_is_{2N-i}$ for each $i \in \{1, \ldots, N-1\}$, and $t_N$ with $s_N$.

\begin{definition}
The {\bf Type B Non-Decreasing Parking Functions} $\BNDPF_N$ are the elements of the submonoid of $\NDPF_{2N}$ generated by the collection $\mu_i := \pi_i\pi_{2N-i}$ for $i$ in the set $\{1, \ldots, N\}$.
\end{definition}

Note that $\mu_N = \pi_N^2 = \pi_N $.

The number of $\BNDPF_N$ has been explicitly computed up to $N=9$, though a proof for a general enumeration has proven elusive, in the absence of a more conceptual description of the full set of functions generated thusly.  The sequence obtained (starting with the $0$-th term) is
\[
	( 1, 2, 7, 33, 183, 1118, 7281, 49626, 349999, 253507, \ldots ),
\]
which agrees with the sequence 
\[
\sum_{j=0}^N \binom{N}{j}^2 C_j 
\]
so far as it has been computed.  This appears in Sloane's On-Line Encyclopedia of Integer Sequences as sequence $A086618$~\cite{Sloane}, and was first noticed by Hivert and Thi\'ery~\cite{Hivert.Thiery.HeckeGroup.2007}.

\begin{conjecture}
\[ 
| \BNDPF_N | = \sum_{j=0}^N \binom{N}{j}^2 C_j.
\]
\end{conjecture}

Let $X$ be some object (group, monoid, algebra) defined by generators $S$ and relations $R$.  Recall that a \emph{parabolic subobject} $X_J$ is generated by a subset $J$ of the set $S$ of simple generators, retaining the same relations $R$ as the original object.  Let $\BNDPF_{N,\hat{N}}$ denote the parabolic submonoid of of $\BNDPF_N$ retaining all generators but $\mu_N$.

Consider the embedding of $\BNDPF_{N,\hat{N}}$ in $\NDPF_{2N}$.  Then a reduced word for an element of $\BNDPF_{N,\hat{N}}$ can be separated into a pairing of $\NDPF_N$ elements as follows:
\begin{eqnarray}
\mu_{i_1}\mu_{i_2}\ldots\mu_{i_k} &=& \pi_{i_1}\pi_{2N-i_1}\pi_{2N-i_2}\pi_{i_2}\ldots\pi_{i_k}\pi_{2N-i_k} \\
&=& \pi_{i_1}\pi_{i_2}\ldots\pi_{i_k}\pi_{2N-i_1}\pi_{2N-i_2}\ldots\pi_{2N-i_k}
\end{eqnarray}

In particular, one can take any element $x\in H_0(S_N)$ and associate it to the pair:
\[
\omega(x):=(\phi(x), \phi\circ\Psi(x)),
\] 
recalling that $\Psi$ is the Dynkin automorphism on $H_0(S_N)$, described in Section~\ref{subsec:involution}.

Given the results of the earlier section, one naturally asks about the fiber of $\omega$.  It is easy to do some computations and see that the situation is not quite so nice as before.  In $H_0(S_4)$ the only fiber with order greater than one contains the elements $[4321]$ and $[4231]$.  Notice what happens here: $[4231]$ contains both a $[231]$-pattern and a $[312]$-pattern, which is straightened into two $[321]$-patterns.   On the level of reduced words, two reduced words for $[4231]$ are $((3,2,1,2,3))=((1,2,3,2,1))$, one of which ends with the unmatched ascent $[2,3]$ while the other ends with the unmatched descent $[2,1]$.  Multiplying on the right by the simple transposition $s_2$ matches both of these simultaneously.

In fact, this is a perfectly general operation.  Let $x\in H_0(S_N)$.  For any minimally-chosen $[231]$-pattern in $x$, one can locate an unmatched ascent in $x$ that corresponds to the pattern.  Here the smaller element to the right remains fixed while the two ascending elements to the left are exchanged.  Then applying the $\NDPF$ relation to turn the $[231]$ into a $[321]$ preserves the fiber of $\phi$.  Likewise, one can turn a minimal $[312]$ into a $[321]$ and preserve the fiber of $\phi\circ\Psi(x)$.  Here the larger element to the left is fixed while the two ascending elements to the right are exchanged.  Hence, to preserve the fiber of $\omega$, one must find a pair of ascending elements with a large element to the left and a small element to the right: this is exactly a $[4231]$-pattern.  

One may make this more precise by defining a system of widths under which minimal $[4231]$-patterns contain a locally minimal $[231]$-pattern and a locally-minimal $[321]$-pattern.  The results of Section~\ref{sec:widthSystems} imply that this is possible.  Applying the $\NDPF$ relation, this becomes a $[4321]$. 

\begin{figure}
  \includegraphics[scale=1]{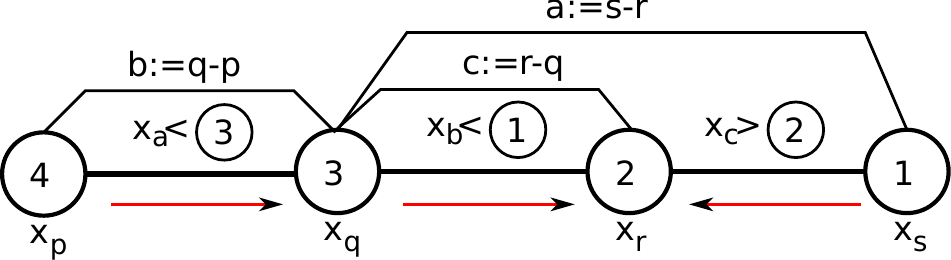}
  \caption{A diagram of a minimal $[4321]$ pattern, labeled analogously to the labeling in Figure~\ref{fig.minimal231}.}
  \label{fig.minimal4321}
\end{figure}

On the other hand, we can define a minimal $[4321]$-pattern by a tuple of widths analogous to the constructions of minimal $[231]$-patterns.  The construction of this tuple, and the constraints implied when the tuple is minimal, is depicted in Figure~\ref{fig.minimal4321}.  Such a minimal pattern may always be turned into a $[4231]$-pattern while preserving the fiber of $\omega$.

Let $x\in S_N$ and $P=(x_p, x_q, x_r, x_s)$ a $[4321]$-pattern in $x$.  For the remainder of this section, we fix the width system $(q-p, r-q, s-r)$, and use the same width system for $[4231]$-patterns.  One may check directly that this is a bountiful width system in both cases.

\begin{lemma}
Let $x$ contain a minimal $[4321]$-pattern $P=(x_p, x_q, x_r, x_s)$, and let $x' = x t_{r,s}$, where $t_{r,s}$ is the transposition exchanging $x_r$ and $x_s$.  Then $\omega(x')=\omega(x)$.
\end{lemma}
\begin{proof}
Since the width system on $[4321]$-patterns is bountiful, we can factor $x = y x_J z$, with $\len(x)=\len(y)+ \len(x_J)+\len(z)$ where 
\[x_J = s_{s-2}s_{s-1} s_s s_{s-1} s_{s-2}s_{s-1}.\]
By the discussion above, the trailing $s_{s-1}$ in $x_J$ may be removed to simultaneously yield an unmatched ascent and an unmatched descent.  Then this removal preserves the fiber of both $\phi$ and $\Psi\circ \phi$, and thus also preserves the fiber of $\omega$.
\end{proof}

Note that there need not be a unique $[4231]$-avoiding element in a given fiber of $\omega$.  The first example of this behavior occurs in $N=7$, where there is a fiber consisting of $[5274163], [5472163],$ and $[5276143]$.  In this list, the first element is $[4321]$-avoiding, and the two latter elements are $[4231]$-avoiding.  In the first element, there are $[4231]$ patterns $[5241]$ and $[7463]$ which can be respectively straightened to yield the other two elements.  Notice that either transposition moves the 4 past the bounding element of the other $[4231]$-pattern, thus obstructing the second transposition.

\begin{theorem}
\label{thm:ndpfFibers4321}
Each fiber of $\omega$ contains a unique $[4321]$-avoiding element.
\end{theorem}

\begin{proof}
Given any element of $H_0(S_N)$, we have seen that we can preserve the fiber of $\omega$ by turning locally minimal $[4321]$-patterns into $[4231]$-patterns.  Each such operation reduces the length of the element being acted upon, and thus this can only be done so many times.  Furthermore, any minimal-length element in the fiber of $\omega$ will be $[4321]$-avoiding.  We claim that this element is unique.  

First, note that one can impose a partial order on the fiber of $\omega$ with $x$ covering $y$ if $x$ is obtained from $y$ by turning a locally minimal $[4321]$-pattern into a $[4231]$-pattern.  Then the partial order is obtained by taking the transitive closure of the covering relation.  Note that if $x$ covers $y$ then $x$ is longer than $y$.  The Hasse diagram of this poset is connected, since any element of the fiber can be obtained from another by a sequence of $\NDPF$ relations respecting both the fiber of $\phi$ and $\phi\circ\Psi(x)$.  

Let $x$ be an element of $H_0(S_N)$ containing (at least) two locally minimal $[4321]$-patterns, in positions $(x_a,x_b,x_c,x_d)$ and $(x_p,x_q,x_r,x_s)$, with $a<b<c<d, p<q<r<s$.  Then one can exchange $x_b$ with $x_c$ or $x_q$ with $x_r$ and preserve the fiber of $\omega$.  Let $y$ be the element obtained from exchanging $x_b$ with $x_c$, and $z$ obtained by exchanging $x_q$ with $x_r$.  Then we claim that there exists $w$ covered by both $y$ and $z$.  (In other words, the poset structure on each fiber is a meet semilattice.)

If the tuples $(a,b,c,d)$ and $(p,q,r,s)$ are disjoint, then the claim is clearly true.  Likewise, if $a=p$ and/or $d=s$ the claim holds.  A complete but perhaps unenlightening proof of the claim can be accomplished by showing that it holds for all $\BNDPF_{N,\hat{N}}$ with $N<8$, where every possible intermingling of the tuples with every possible ordering of the entries $x_.$ occurs at least once.  It is best to perform this check with a computer, given that there are 2761 elements in $\BNDPF_{7,\hat{7}}$, with $7! = 5040$ elements in the fibers, and indeed a computer check shows that the claim holds.  The code accomplishing this is provided below.

Let's look at a couple cases, though, to get a feeling for why this should be true.  Refer to the extremal elements at the edge of the $[4321]$ pattern as the ``boundary,'' and the elements to be transposed as the ``interior.''  The main cases are the following:

Case $c=r$: Just take the smaller of $x_s$ and $x_d$ to be the common right boundary for both patterns.

Case $c=q$: The problem for $[4231]$ patterns was that one could apply a transposition that obstructed the other transposition by moving one of the interior elements past its boundary.  But here, we have $x_d<x_c$ and $x_s<x_q=x_d$, so we can use $s$ as the boundary for both patterns, and the obstruction is averted.  In this case, though, the two transpositions generate six elements in the fiber, instead of four.  We can still find a common meet, though.  $[x_ax_bx_qx_rx_s]$ becomes $[x_ax_qx_bx_rx_s]$ and $[x_ax_bx_rx_qx_s]$, which both cover $[x_ax_rx_bx_qx_s]$, for example.

Case $r=d$ or $q=d$: Again, just take $s$ as a common boundary for the two patterns.

And so on.  Many cases are symmetric to the three considered above, and every interesting case is solved by changing the boundary of one of the patterns.

Now that every pair of elements have a common meet, we are almost done.  Suppose there exist two different $[4321]$-avoiding elements $A_1$ and $A_2$ in some fiber.  Then since the fiber is connected, we can find a minimal element $x$ where a branching occurred, so that $x$ covers both $y>A_1$ and $z>A_2$, and $x$ is of minimal length.  But if both $y$ and $z$ were obtainable from $x$, then there exists a $w$ of shorter length below them both.  Now $w$ sits above some $[4321]$-avoiding element, as well.  If $w>A_1$ but not $A_2$, then in fact a branching occurred at $z$, contradicting the minimality of $x$.  The same reasoning holds if $w>A_2$ but not $A_1$.  If $w$ is above both $A_1$ and $A_2$, then in fact $y$ was comparable to $A_2$ and $z$ was comparable to $A_1$, and there was not a branching at $x$ at all.
\end{proof}

\subsection{Code for Theorem~\ref{thm:ndpfFibers4321}.}
Here we provide code for checking the claim of Theorem~\ref{thm:ndpfFibers4321} that each fiber of $\omega$ contains a unique $[4321]$-avoiding element.  The code is written for the Sage computer algebra system, which has extensive built-in functions for combinatorics of permutations, including detecting the presence of permutation patterns.  

The code below constructs a directed graph (see the function \textbf{omegaFibers}) whose connected components are fibers of $\omega$.  The vertices of this graph are permutations, and the edges correspond to straightening locally-minimal $[4231]$-patterns into $[4321]$ patterns.  A component is `bad' if it does not contain exactly one $[4321]$-avoiding permutation.
\begin{verbatim}
def width4231(p):
    """
    This function returns the width of a [4231]-instance p.
    """
    return (p[1]-p[0], p[2]-p[1], p[3]-p[2]) 

def min4231(x):
    """
    This function takes a permutation x and finds all minimal-width
    4231-patterns in x, and returns them as a list.
    """
    P=x.pattern_positions([4,2,3,1])
    if P==[]:
        return None
    minimal=[P[0]]
    for i in [1..len(P)-1]:
        if width4231(P[i])<width4231(minimal[0]):
            minimal = [ P[i] ]
        else:
            if width4231(P[i])==width4231(minimal[0]):
                minimal.append(P[i])
    return minimal
    
def localMin4231(x):
    """
    This function finds all locally-minimal 4231-patterns in a 
    permutation x, and returns them as a list.
    """
    P=x.pattern_positions([4,2,3,1])
    if P==[]:
        return None
    localMin=[]
    for p in P:
        xp=Permutation(x[ p[0]:p[3]+1 ])
        qp=[i - p[0] for i in p]
        qmin=min4231(xp)
        if qp in qmin: localMin.append(p)
    return localMin

def omegaFibers(N):
    """
    Given N, this function builds a digraph whose vertices are given by
    permutations of N, and with an edge a->b whenever b is obtained
    from a by straightening a locally minimal 4231-pattern into a 
    4321-pattern. 
    The connected components of G are the fibers of the map omega.
    """
    S=Permutations(N)
    G=DiGraph()
    G.add_vertices(S.list())
    for x in S:
        if x.has_pattern([4,2,3,1]):
            # print x, localMin4231(x)
            #add edges to G for each locally minimal 4231.
            Q=localMin4231(x)
            for q in Q:
                y=Permutation((q[1]+1,q[2]+1))*x
                G.add_edge(x,y)
    return G
    
def headCount(G):    
    """
    This function takes the diGraph G produced by the omegaFibers 
    function, and finds any connected components with more than one
    4321-pattern.  It returns a list of all such connected components.
    """
    bad=[]
    for H in G.connected_components_subgraphs():
        total=0
        for a in H:
            if not a.has_pattern([4,3,2,1]): total+=1
        if total != 1:
            #prints if any fiber has more than one 4321-av elt
            print H, total
            bad.append(H)
    print "N =", N 
    print "\tTotal connected components: \t", count
    print "\tBad connected components: \t", len(bad), '\n'
    return bad
\end{verbatim}

As explained in Theorem~\ref{thm:ndpfFibers4321}, we should check that each fiber of $\omega$ contains a unique $[4321]$-avoiding element for each $N\leq 7$.  This is accomplished by running the following commands:
\begin{verbatim}
sage: for N in [1..7]:
sage:     G=omegaFibers(N)
sage:     HH=headCount(G)
\end{verbatim}
The output of this loop is as follows:
\begin{verbatim}
N = 1
	Total connected components: 	1
	Bad connected components: 	0 

N = 2
	Total connected components: 	2
	Bad connected components: 	0 

N = 3
	Total connected components: 	6
	Bad connected components: 	0 

N = 4
	Total connected components: 	23
	Bad connected components: 	0 

N = 5
	Total connected components: 	103
	Bad connected components: 	0 

N = 6
	Total connected components: 	513
	Bad connected components: 	0 

N = 7
	Total connected components: 	2761
	Bad connected components: 	0 
\end{verbatim} 
There are no bad components, and thus the theorem holds.  

The sequence $(1, 2, 6, 23, 103, 513, 2761)$ is the beginning of the sequence counting $[4321]$-avoiding permutations.  This sequence also counts $[1234]$-avoiding permutations (reversing a $[1234]$-avoiding permutation yields a $[4321]$-avoiding permutation, and \textit{vice versa}), and is listed in that context in Sloane's On-Line Encyclopedia of Integer Sequences (sequence $A005802$)~\cite{Sloane}.

The author executed this code on a computer with a 900-mhz Intel Celeron processor (blazingly fast by 1995 standards) and 2 gigabytes of RAM.  On this machine, the $N=6$ case took 3.86 seconds of CPU time, and the $N=7$ case took just over one minute (62.06s) of CPU time.  The $N=8$ case (which is unnecessary to the proof) correctly returns 15767 connected components, none of which are bad, and took 1117.24 seconds (or 18.6 minutes) to run.

\section{Affine $\NDPF$ and Affine $[321]$-Avoidance}
\label{sec:affNdpfPattAvoid}

The affine symmetric group is the Weyl group of type $A_N^{(1)}$, whose Dynkin diagram is given by a cycle with $N$ nodes.  All subscripts on generators for type $A_N^{(1)}$ in this section will be considered $(\text{mod } N)$.  A combinatorial realization of this Weyl group is given below.

\begin{definition}
\label{def:affSn}
The \textbf{affine symmetric group} $\ASn_N$ is the set of bijections $\sigma: \ZZ \rightarrow \ZZ$ satisfying:
\begin{itemize}
    \item Skew-Periodicity: $\sigma(i+N) = \sigma(i) + N$, and
    \item Sum Rule: $\sum_{i=1}^N \sigma(i) = \binom{N+1}{2}$.
\end{itemize}
\end{definition}

We will often denote elements of $\ASn_N$ in the \textbf{window notation}, which is a one-line notation where we only write $(\sigma(1), \sigma(2), \ldots, \sigma(N))$.  Due to the skew-periodicity restriction, writing the window notation for $\sigma$ specifies $\sigma$ on all of $\ZZ$.  

The generators $s_i$ of $\ASn_N$ are indexed by the set $I=\{0, 1, \ldots, N-1\}$, and $s_i$ acts by exchanging $j$ and $j+1$ for all $j \equiv i (\text{mod } N)$.  These satisfy the relations:
\begin{itemize}
    \item Reflection: $s_i^2 = 1$,
    \item Commutation: $s_j s_i = s_i s_j$ when $|i-j| > 1$, and
    \item Braid Relations: $s_i s_{i+1} s_i = s_{i+1} s_i s_{i+1}$. \\
\end{itemize}
In these relations, all indices should be considered mod $N$.

Since the Dynkin diagram is a cycle, it admits a dihedral group's worth of automorphisms.  One can implement a ``flip'' automorphism $\Phi$ by fixing $s_0$ and sending $s_i \rightarrow s_{N-i}$ for all $i \neq 0$, extending the automorphism used in the finite case.  A ``rotation'' automorphism $\rho$ can be implemented by simply sending each generator $s_i \rightarrow s_{i+1}$.  Combinatorially, this corresponds to the following operation.  Given the window notation $(\sigma_1, \sigma_2, \ldots, \sigma_N)$, we have:
\[
\rho(\sigma) = (\sigma_N-N+1, \sigma_1 +1, \sigma_2 +1, \ldots, \sigma_{N-1}+1).
\]
This can be thought of as shifting the base window one place to the left, and then adding one to every entry.  It is clear that this operation preserves the skew periodicity and sum rules for affine permutations, and it is also easy to see that $\rho^N = 1$.

As before, we can define the Hecke algebra of $\ASn_N$, and the $0$-Hecke algebra, generated by $\pi_i$ with $\pi_i$ idempotent anti-sorting operators, exactly mirroring the case for the finite symmetric group.  As in the finite case, elements of the $0$-Hecke algebra are in bijection with affine permutations.  We can also define the $\NDPF$ quotient of $H_0(\ASn_N)$, by introducing the relation 
\[
\pi_{i+1}\pi_i \pi_{i+1} = \pi_{i+1} \pi_i.
\]  
This allows us to give combinatorial definition for the affine $\NDPF$, which we will prove to be equivalent to the quotient.

\begin{definition}
The extended affine non-decreasing parking functions are the functions $f: \ZZ\rightarrow \ZZ$ which are:
\begin{itemize}
    \item Regressive: $f(i)\leq i$,
    \item Order Preserving: $i\leq j \Rightarrow f(i)\leq f(j)$, and
    \item Skew Periodic: $f(i+N) = f(i)+N$.
\end{itemize}
Define the \textbf{shift functions} $\operatorname{sh}_t$ as the functions sending $i \rightarrow i-t$ for every $i$.

The \textbf{affine non-decreasing parking functions} $\ANDPF_N$ are obtained from the extended affine non-decreasing parking functions by removing the shift functions for all $t\neq 0$.  
\end{definition}

Notice that the definition implies that
\[
f(N)-f(1) \leq N.
\]
Furthermore, since the shift functions are not in $\ANDPF_N$, there is always some $j \in \{0, 1, \ldots, N\}$ such that $f(j)\neq f(j+1)$ unless $f$ is the identity.

We now state the main result of this section, which will be proved in pieces throughout the remainder of the chapter.
\begin{theorem}
\label{andpfMainThm}
The affine non-decreasing parking functions $\ANDPF_N$ are a $\JJ$-trivial monoid which can be obtained as a quotient of the $0$-Hecke monoid of the affine symmetric group by the relations $\pi_j\pi_{j+1}\pi_j = \pi_j \pi_{j+1}$, where the subscripts are interpreted modulo $N$.  Each fiber of this quotient contains a unique $[321]$-avoiding affine permutation.
\end{theorem}

\begin{figure}
  \includegraphics[scale=.75]{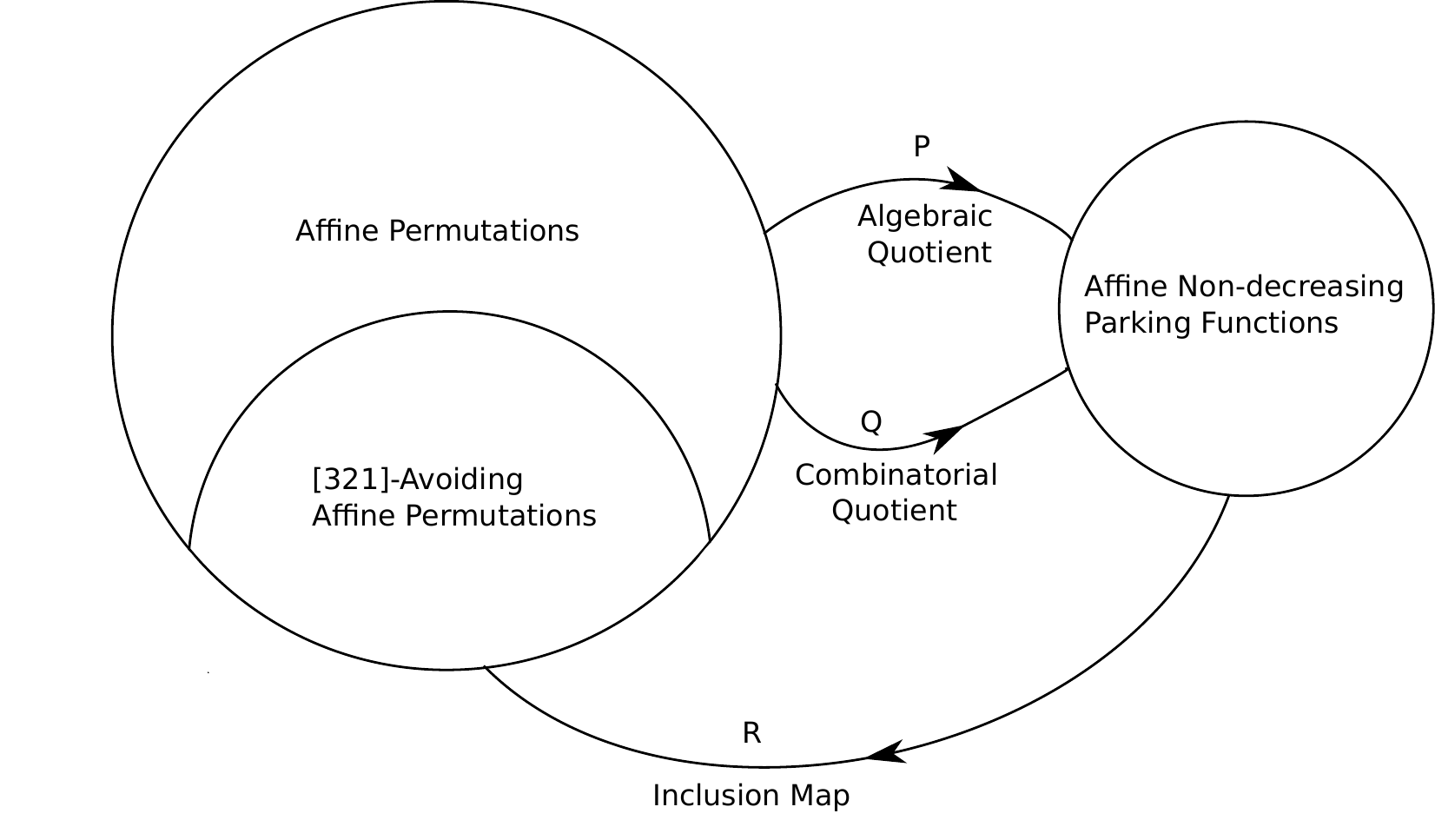}
  \caption{Maps between $H_0(\ASn_N)$ and $\ANDPF_N$.}
  \label{fig.affineNDPFMaps}
\end{figure}

\begin{proposition}
\label{andpfgens}
As a monoid, $\ANDPF_N$ is generated by the functions $f_i$ defined by:
\begin{displaymath}
   f_i(j) = \left\{
     \begin{array}{lr}
       j-1 & : j\equiv i+1 (\text{mod }N)\\
       j   & : j\not \equiv i+1 (\text{mod }N).\\
     \end{array}
   \right.
\end{displaymath} 
These functions satisfy the relations:
\begin{eqnarray*}
f_i^2 &=& f_i \\
f_if_j &=& f_jf_i \text{ when $|i-j|>1$, and} \\
f_if_{i+1}f_i = f_{i+1}f_if_{i+1} &=& f_{i+1}f_i \text{ when $|i-j|=1$,}
\end{eqnarray*}
where the indices are understood to be taken $(\text{mod }N)$.
\end{proposition}
\begin{proof}
One can easily check that these functions $f_i$ satisfy the given relations.  We then check that any $f\in \ANDPF_N$ maybe written as a composition of the $f_i$.

Let $f \in \ANDPF_N$.  If there is no $j \in \{0, \ldots, N\}$ such that $f(j)=f(j+1)$, then $f$ is a shift function, and is thus the identity.

Otherwise, we have some $j$ such that $f(j)=f(j+1)$.  We can then build $f$ using $f_i$'s by the following procedure.  Notice that, if any $g \in \ANDPF$ has $g(j)=g(j+1)$ for some $j$, we can emulate a shift function by concatenating $g$ with $f_j f_{j+1} \cdots f_{j+N-1}$, where the subscripts are understood to be taken $(\text{mod }N)$.  In other words, we have:
\[
g \operatorname{sh}_1 = g f_j f_{j+1} \cdots f_{j+N-1}.
\]

Suppose, without loss of generality, that $f(N)\neq f(N+1)$, so that $N$ and $N+1$ are in different fibers of $f$, and $N$ is maximal in its fiber.  (If the ``break'' occurs elsewhere, we simply use that break as the `top' element for the purposes of our algorithm.  Alternately, we can apply the Dynkin automorphism to $f$ until $\rho^k f(N)\neq \rho^k f(N+1)$. for some $k$.  We can use this algorithm to construct $\rho^k f$, and then apply $\rho$ $N-k$ times to obtain $f$.)  Begin with $g = 1$, and construct $g$ algorithmically as follows.

\begin{itemize}
\item Collect together the fibers.  Set $g'$ to be the shortest element in $\NDPF_N$ such that the fibers of $g'$ match the fibers of $f$ in the base window.  Let $g_0$ be the affine function obtained from a reduced word for $g'$.  This is the pointwise maximal function in $\ANDPF_N$ with fibers equal to the fibers of $f$.

\item Now that the fibers are collected, post-compose $g_0$ with $f_i$'s to move the images into place.  We begin with $g:=g_0$ and apply the following loop:\\
\begin{eqnarray*}
&&\text{while } g\neq f: \\
&&\phantom{aaaa}\text{for } i \text{ in } \{1, \ldots, N \}: \\
&&\phantom{aaaaaaaa}\text{if } g(i+1)>f(i+1) \text{ and } g^{-1}(g(i+1)-1) = \emptyset: \\
&&\phantom{aaaaaaaaaaaa}g := g.f_i.
\end{eqnarray*}

This process clearly preserves the fibers of $g_0$ (which coincide with the fibers of $f$),  and terminates only if $g=f$.  We need to show that the algorithm eventually halts.

Recall that $g_0(i)\geq f(i)$ for all $i$, and then notice that it is impossible to obtain any $g$ in the evaluation of the algorithm with $g(i)<f(i)$, so that we always have $g(i)-f(i)>0$.  With each application of a $f_j$, the sum $\sum_{i=1}^N (g(i)-f(i))$ decreases by one.

Suppose the loop becomes stuck; then for every $i$ either $f(i+1)=g(i+1)$ or $g^{-1}(g(i+1)-1) \neq \emptyset$.  If there is no $i$ with $f(i+1)=g(i+1)$, then there must be some $i$ with $g^{-1}(g(i+1)-1) = \emptyset$, since $g(N)-g(1)\leq N$ and $g \neq 1$.  Then we can find a minimal $i \in \{1, \ldots, N\}$ with $f(i+1)=g(i+1)$.  

Now, find $j$ minimal such that $f(i+j) \neq g(i+j)$, so that $f(i+j-1) = g(i+j-1)$.  In particular, notice that $i+j-1$ and $i+j$ must be in different fibers for both $f$ and $g$.  If $g^{-1}(g(i+j)-1) = \emptyset$, then the loop would apply a $f_{i+j-1}$ to $g$, but the loop is stuck, so this does not occur and we have that $f(i+j-1)=g(i+j-1)=g(i+j)-1< f(i+j)\leq g(i+j) = g(i+j-1)+1$.  This then forces $g(i+j)=f(i+j)$, contradicting the condition on $j$.

Thus, the loop must eventually terminate, with $g=f$.
\end{itemize}

We have not yet shown that these relations are all of the relations in the monoid; this must wait until we have developed more of the combinatorics of $\ANDPF_N$.  In fact, $\ANDPF_N$ is a quotient of the $0$-Hecke monoid of $\ASn_N$ by the relations $\pi_i \pi_{i+1} \pi_i = \pi_i\pi_{i+1}$ for each $i \in I$, where subscripts are understood to be taken $\text{mod } N$.  To prove this (and simultaneously prove that we have in fact written all the relations in $\ANDPF_N$), we will define three maps, $P, Q$, and $R$ (illustrated in Figure~\ref{fig.affineNDPFMaps}).  The map $P: H_0(\ASn_N)\rightarrow \ANDPF_N$ is the algebraic quotient on generators sending $\pi_i \rightarrow f_i$.  The map $Q: H_0(\ASn_N)\rightarrow \ANDPF_N$ is a combinatorial algorithm that assigns an element of $\ANDPF_N$ to any affine permutation.  In Lemma~\ref{lem.combQuotient} we show that $P=Q$.  Additionally, we have already shown that $P$ is onto (since the $f_i$ generate $\ANDPF_N$), so $Q$ is onto as well.

The third map $R: \ANDPF_N \rightarrow H_0(\ASn_N)$ assigns a $[321]$-avoiding affine permutation to an $f\in \ANDPF_N$.  In fact, $R \circ P$ is the identity on the set of $[321]$-avoiding affine permutations, and $P\circ R$ is the identity on $\ANDPF_N$.  This then implies that there are no additional relations in $\ANDPF_N$.
\end{proof}

\begin{corollary}
The map $P: H_0(\ASn_N)\rightarrow \ANDPF_N$, defined by sending $\pi_i \rightarrow f_i$ and extending multiplicatively, is a monoid morphism.
\end{corollary}
\begin{proof}
The generators $f_i$ satisfy all relations in the $0$-Hecke algebra, so $P$ is a quotient of $H_0(\ASn_N)$ by whatever additional relations exist in $\ANDPF_N$.
\end{proof}

\begin{lemma}
Any function $f \in \ANDPF_N$ is entirely determined by its set of fibers, set of images, and one valuation $f(i)$ for some $i\in \ZZ$.
\end{lemma}
\begin{proof}
This follows immediately from the fact that $f$ is regressive and order preserving.
\end{proof}

\begin{lemma}
Let $f \in \ANDPF_N$, and $F_f=\{m_j\}$ be the set of maximal elements of the fibers of $f$.  Each pair of distinct elements $m_j, m_k$ of the set $F_f \cap \{1, 2, \ldots, N\}$ has $f(m_j) \not \equiv f(m_k) (\text{mod } N)$.
\end{lemma}
\begin{proof}
Suppose not.  Then $f(m_j)-f(m_i) = kN$ for some $k \in \ZZ$, implying that 
$f(m_j)=f(m_i+kN)$.  Since $f(m_j)-f(m_i)\leq N$, we must have $k=0$.  But then $m_j$ and $m_i$ are in the same fiber, providing a contradiction.
\end{proof}

\begin{theorem}
$\ANDPF_N$ is $\JJ$-trivial.
\end{theorem}
\begin{proof}
Thi is a direct consequence of the regressiveness of functions in $\ANDPF_N$.  Let $M:= \ANDPF_N$, and $f\in M$.  Then each $g\in MfM$ has $g(i)\leq f(i)$ for all $i \in \ZZ$.  Thus, if $MgM=MfM$, we must have $f=g$.  Then the $\JJ$-equivalence classes of $M$ are trivial, so $\ANDPF_N$ is $\JJ$-trivial.
\end{proof}

Note that $\ANDPF_N$ is not aperiodic in the sense of a finite monoid.  (Aperiodicity was defined in Section~\ref{sec:bgnot}.)  Take the function $f$ where $f(i)=0$ for all $i \in \{1, \ldots, N\}$.  Then $f^k(1) = (1-k)N$, so there is no $k$ such that $f^k = f^{k+1}$.

\subsection{Combinatorial Quotient}
\label{subsec:combintorialQuotient}

\begin{figure}
  \includegraphics[scale=1]{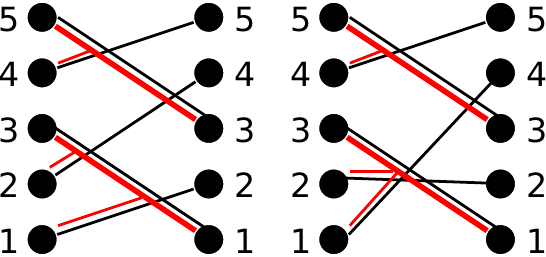}
  \caption{Example of the combinatorial quotient $Q_{cl}: H_0(S_5) \rightarrow \NDPF_5$.  The string diagram is read left-to-right, with the permutation illustrated with black strings and the image function drawn in red.  The permutation in the left diagram, then, is $x=[2,4,1,5,3]$ and $Q_{cl}(x)$ is the function $f=[1,1,1,3,3]$.  For the permutation on the right, we have $y=[4,2,1,5,3]$ and $Q_{cl}(y)=Q_{cl}(x)=[1,1,1,3,3]$.  Notice that these two permutations $x$ and $y$ are related by turning the $[321]$-pattern in $y$ into a $[231]$-pattern in $x$, preserving the fiber of $Q$.}
  \label{fig.combQuotient}
\end{figure}

A direct combinatorial map from affine permutations to $\ANDPF_N$ is now discussed.  This map directly constructs a function $f$ from an arbitrary affine permutation $x$, with the same effect as applying the algebraic $\ANDPF$ quotient to the $0$-Hecke monoid element indexed by $x$.  We first define the combinatorial quotient in the finite case and provide an example (Figure~\ref{fig.combQuotient}).

\begin{definition}
The combinatorial quotient $Q_{cl}: H_0(S_N) \rightarrow \NDPF_N$ is given by the following algorithm, which assigns a function $f$ to a permutation $x$.
\begin{enumerate}
\item Set $f(N):=x(N)$.
\item Suppose $i$ is maximal such that $f(i)$ is not yet defined.  If $x(i)<x(i+1)$, set $f(i):=f(i+1)$.  Otherwise, set $f(i):=x(i)$.
\end{enumerate}
\end{definition}

Note that the map $Q_{cl}$ is closely related to bijection of Simion and Schmidt between $[132]$-avoiding permutations and $[123]$-avoiding permutations~\cite{simion.schmidt.1985}.  (The bijection is also covered very nicely in~\citetalias{bona.permutations})  This bijection operates by marking all left-to-right minima (ie, elements smaller than all elements to their left) of a $[132]$-avoiding permutation, and then reverse-sorting all elements which are not marked. The resulting permutation is $[123]$-avoiding.  For example, the permutation
$[\textbf{5}, 6, \textbf{4}, 7, \textbf{1}, 2, 3]$ avoids the pattern $[132]$; the bold entries are the left-to-right minima.  Sorting the non-bold entries, one obtains the permutation
$[\textbf{5}, 7, \textbf{4}, 6, \textbf{1}, 3, 2]$, which avoids the permutation $[123]$.  Notice that the bold entries are still left-to-right minima after anti-sorting the other entries.

The patterns $[231]$ and $[123]$ are the respective ``reverses'' of the patterns $[132]$ and $[321]$, obtained by simply reversing the one-line notation.  It is trivial to observe that $x$ avoids $p$ if and only if the reverse of $x$ avoids the reverse of $p$.  Then the ``reverse'' of the Simion-Schmidt algorithm (which marks right-to-left minima, and sorts the other entries) gives a bijection between $[231]$- and $[321]$-avoiding permutations; in fact, this is the same bijection given by the fibers of the $\NDPF$ quotient of the $0$-Hecke monoid.

A similar combinatorial quotient may be defined from $\ASn_N\rightarrow \ANDPF_N$, generalizing the map $Q_{cl}$.  This map will assign a function $f$ to an affine permutation $x$.

Below, we will show that each fiber of the map $Q$ contains a unique $[321]$-avoiding affine permutation (Theorem~\ref{thm:affNdpfFibers321}).  However, it is too much to expect a bijection between affine $[231]$- and $[321]$-avoiding permutations.  By a result of Crites, there are infinitely many affine permutations that avoid a pattern $\sigma$ if and only if $\sigma$ contains the pattern $[321]$~\cite{Crites.2010}.  Thus, there are infinitely many $[321]$-avoiding affine permutations, but only finitely many $[231]$-avoiding affine permutations.

We first identify some $k \in \{1, 2, \ldots, N\}$ such that for every $j>k$, $x(j)>x(k)$.  

\begin{lemma}
Let $k_0 \in \{1, 2, \ldots, N\}$ have $x(k_0)\leq x(m)$ for every $m \in \{1, 2, \ldots, N\}$.  Then for every $j>k_0$, $x(j)>x(k_0)$.
\end{lemma}
\begin{proof}
 Suppose $j>k_0$ with $x(j)<x(k_0)$.  Then there exists $p\in \N$ such that $j-pN \in \{1, 2, \ldots, N\}$, so that $x(j-pN) = x(j)-pN < x(k_0)$, contradicting the minimality of $x(k_0)$.
\end{proof}

Now the affine combinatorial quotient is defined by the following algorithm.

\begin{definition}
The combinatorial quotient $Q: H_0(\ASn_N) \rightarrow \ANDPF_N$ is given by the following algorithm, which assigns a function $f$ to an affine permutation $x$.
\begin{enumerate}
\item Let $k_0 \in \{1, 2, \ldots, N\}$ have $x(k_0)\leq x(m)$ for every $m \in \{1, 2, \ldots, N\}$.  Set $f(k_0)=x(k_0)$.
\item Choose $i \in \{1, 2, \ldots, N-1\}$ minimal such that $f(k_0 - i)$ is not yet defined.  If $x(k_0-i+1)<x(k_0-i)$, set $f(k_0-i) := f(k_0-i+1)$.  Otherwise, set $f(k_0-i) := x(k_0-i)$.
\item Define $f$ on all other $i$ using skew periodicity.
\end{enumerate}
\end{definition}

\begin{lemma}
\label{lem.combQuotient}
The affine combinatorial quotient $Q$ agrees with the algebraic $\ANDPF$ quotient $P$. 
\end{lemma}
\begin{proof}
We denote the combinatorial quotient by $Q$ and the algebraic quotient by $P$.  

One can easily check that $Q(1)=P(1)=1$, and $Q(\pi_i)=P(\pi_i)=f_i$.  Since $P$ is a monoid morphism, we have that $P(x\pi_i)=P(x)P(\pi_i)= ff_i$.   We then assume that $Q(x)=P(x)=f$, and consider $Q(x\pi_i)$.  We will show that $Q(x\pi_i)= Q(x)f_i = ff_i = P(x\pi_i)$.  

If $\pi_i$ is a right descent of $x$ then $Q(x\pi_i)=Q(x)=f=P(x\pi_i)$, and we are done.

If $\pi_i$ is not a right descent of $x$, we have $x(kN+i) < x(kN+i+1)$ for all $k\in \ZZ$, and 
\[
x\pi_i(j) = 
    \left\{
     \begin{array}{lr}
       x(j)   \text{ for all } j\not \equiv i, i+1 (\text{mod } N)\\
       x(j+1) \text{ for all } j     \equiv i      (\text{mod } N)\\
       x(j-1) \text{ for all } j     \equiv    i+1 (\text{mod } N)\\
     \end{array}
   \right.
\]  
We examine the functions $Q(x\pi_i)$ and $ff_i$ on $i$ and $i+1$, since these functions are equal on $j\not \equiv i, i+1 (\text{mod } N)$, and the actions on $i$ and $i+1$ then determine the functions on all $j \equiv i, i+1 (\text{mod } N)$. 

We consider two cases, depending on whether $i$ and $i+1$ are in the same fiber of $f$.  

\begin{itemize}
\item If $i$ and $i+1$ are in the same fiber of $f$ and $i+1$ is maximal in this fiber, we must (by construction of $Q$) have $x(i+1)<x(i)$, contradicting the assumption that $\pi_i$ was not a right descent of $x$.

\item If $i$ and $i+1$ are in the same fiber of $f$ and $i+1$ is not maximal in this fiber, then there exists some (minimal) $m> i+1>i$ with $x(m)<x(i)$ and $x(m)<x(i+1)$, maximal in the fiber of $i$ and $i+1$.  Then $x(m)<x(i+1)=x\pi_i(i)$ and $x(m)<x(i)=x\pi_i(i+1)$.  Since the maximal size of a fiber of $f$ is $N$, we have that $m-i\leq N$.  Then (since $i+1$ not maximal in the fiber of $f$) $m\not \equiv i+1 (\text{mod } N)$.

If $m \equiv i (\text{mod } N)$, then $i$ is maximal in its fiber, and we must have $i$ and $i+1$ in different fibers, contrary to assumption.

If $m\not \equiv i (\text{mod } N)$, we have $x(m)=x\pi_i(m) < x\pi_i(i), x\pi_i(i+1)$, and so by the construction of $Q$, we have 
$Q(x\pi_i)(i)=Q(x\pi_i)(i+1)=Q(x\pi_i)(m)=Q(x)(m) = x(m)$.  Then in this case, $Q(x\pi_i)=f$.

On the other hand, $ff_i(i)=f(i)=f(m)=ff_i(m)$, and $ff_i(i+1)=f(i)=f(m)=ff_i(m)$, so $ff_i=f$.

\item If $i$ and $i+1$ are in different fibers of $f$, then we have $i$ maximal in its fiber, and take $m$ (possibly equal to $i+1$) to be the maximal element of the fiber in which $i+1$ sits.  We note that if $m\equiv i+1 (\text{mod } N)$, then we must have $i$ and $i+1$ in the same fiber, reducing to the previous case.

Otherwise, applying the construction of $Q$, we find that $Q(x\pi_i)(i+1) = x(i)$, and that $Q(x\pi_i)(i) = x(i)$; thus $i+1$ is removed from its fiber and merged into the fiber with $i$.  The resulting function is equal to $ff_i$.
\end{itemize}
This exhausts all cases, completing the proof.
\end{proof}

\begin{corollary}
The finite type combinatorial quotient agrees with the $\NDPF_N$ quotient of $H_0(S_N)$ obtained by introducing the relations $\pi_i \pi_{i+1}\pi_i = \pi_{i+1} \pi_i$ for $i\in \{1, \ldots, N-2\}$. 
\end{corollary}
\begin{proof}
This follows immediately from Lemma~\ref{lem.combQuotient} by parabolic restriction to the finite case.  In the finite case, the index set is $\{1, 2, \ldots, N-1\}$, so we must have $i\in \{1, \ldots, N-2\}$.
\end{proof}

\subsection{Affine $[321]$-Avoidance}

An affine permutation $x$ avoids a pattern $\sigma \in S_k$ if there is no subsequence of $x$ in the same relative order as $\sigma$.  This ostensibly means that an infinite check is necessary, however one may show that only a finite number of comparisons is necessary to determine if $x$ contains a $[321]$-pattern.  The following lemma is equivalent to~\cite[Lemma 2.6]{Green.2002}.

\begin{lemma}
Let $x$ contain at least one $[321]$-pattern, with $x_i> x_j> x_k$ and $i<j<k$.  Then $x$ contains a $[321]$-pattern $x_{i'}> x_j> x_{k'}$ such that $i\leq i'<j<k'\leq k$, $j-i'< N$, and $k'-j < N$.
\end{lemma}
\begin{proof}
We have $x_j>x_k>x_{k-aN} = x_k - aN$ for $a \in \N$, so if $k-j>N$, we can find a $[321]$ pattern replacing $x_k$ with $x_{k-aN}$.  A similar argument allows us to replace $i$ with $i+bN$ for the maximal $b\in \N$ such that $j-(i+bN) < N$.
\end{proof}

\begin{figure}
  \includegraphics[scale=1]{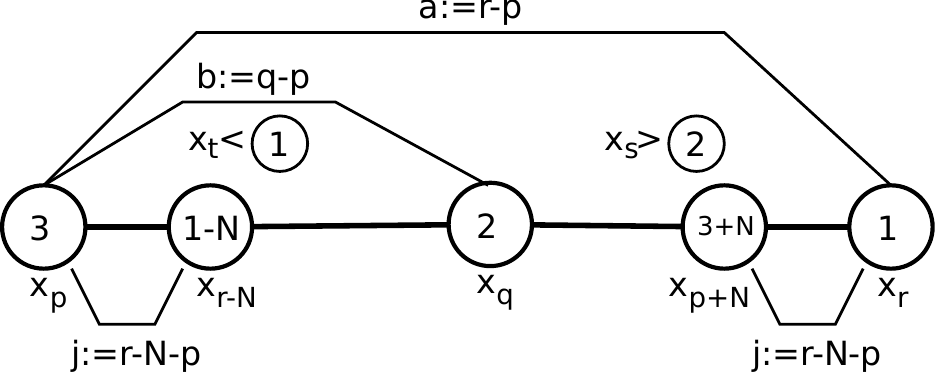}
  \caption{Diagram of a bountiful width system for the pattern $[321]$ for affine permutations.  The pattern occurs at positions $(x_p, x_q, x_r)$, with width system given by $(r-p, q-p)$.  In the case where $r-p>N$, there is an `overlap' of $j=r-N-p$.  Bountifulness of the width system ensures that the elements in the overlap may be moved moved out of the interior of the pattern instance by a sequence of simple transpositions, each decreasing the length of the permutation by one, just as in the non-affine case.}
  \label{fig.minimalAffine321}
\end{figure}

As noted by Green, one can then check whether an affine permutation contains a $[321]$-pattern using at most $\binom{N}{3}$ comparisons.  Green also showed that any affine permutation containing a $[321]$-pattern contains a braid; we can actually replicate this result using a width system on the affine permutation, as depicted in Figure~\ref{fig.minimalAffine321}.  The Lemma ensures that the width of a minimal $[321]$-pattern under this width system has a total width of at most $2N-2$.  One must consider the case when the total width of a minimal $[321]$-instance is greater than $N$, but nothing untoward occurs in this case: the width system is bountiful and allows a factorization of $x$ over $[321]$.

We now prove the main result of this section.

\begin{theorem}
\label{thm:affNdpfFibers321}
Each fiber of the $\ANDPF_N$ quotient of $\ASn_N$ contains a unique $[321]$-avoiding affine permutation.
\end{theorem}
\begin{proof}
We first establish that each fiber contains a $[321]$-avoiding affine permutation, and then show that this permutation is unique.

Recall the algebraic quotient map $P: H_0(\ASn_N) \rightarrow \ANDPF_N$, which introduces the relation $\pi_i \pi_{i+1}\pi_i = \pi_{i+1} \pi_i$.  

Choose an arbitrary affine permutation $x$; we show that the fiber $Q^{-1}\circ Q(x)$ contains a $[321]$-avoiding permutation.  If $x$ is itself $[321]$-avoiding, we are already done.  So assume $x$ contains a $[321]$-pattern.  As shown by Green~\cite{Green.2002}, an affine permutation $x$ contains a $[321]$-pattern if and only if $x$ has a reduced word containing a braid; thus, $x=y \pi_i \pi_{i+1} \pi_i z$ for some permutations $y$ and $z$ with $\len(x)=\len(y)+3+\len(z)$.  Applying the $\ANDPF_N$ relations, we may set $x'= y \pi_{i+1} \pi_i z$, and have $Q(x)=Q(x')$, with $\len(x')=\len(x)-1$.  If $x'$ contains a $[321]$, we apply this trick again, reducing the length by one.  Since $x$ is of finite length, this process must eventually terminate; the permutation at which the process terminates must then be $[321]$-avoiding.  Then the fiber $Q^{-1}\circ Q(x)$ contains a $[321]$-avoiding permutation.

We now show that each fiber contains a unique $[321]$-avoiding affine permutation, using the combinatorial quotient map.

Let $x$ be $[321]$-avoiding, and let $Q(x)=f$ an affine non-decreasing parking function; we use information from $f$ to reconstruct $x$.  Let $\{m_i\}$ be the set of elements of $\ZZ$ that are maximal in their fibers under $f$.  By the construction of the combinatorial quotient map, we have $x(m_i)=f(m_i)$ for every $i$.  Since $f$ is in $\ANDPF_N$, we have $x(m_i)<x(m_{i'})$ whenever $i<i'$; thus $\{x(m_i)\}$ is a strictly increasing sequence.

Let $\{m_{i,j}\} = f^{-1}\circ f(m_i) \setminus \{m_i\}$, with $m_{i,j}<m_{i,j+1}$ for every $j$.  Notice that if $i<i'$ and $j<j'$ then $m_{i,j}<m_{i',j'}$.  

We claim that if $i<i'$ and $j<j'$, then $x(m_{i,j})<x(m_{i',j'})$.  If not, then we have 
\[
x(m_{i'}) < x(m_{i',j'}) < x(m_{i,j}), \text{ with }  m_{i,j} < m_{i',j'} < m_{i'},
\]
in which case $x$ contains a $[321]$-pattern, contrary to assumption.  Thus, the sequence $\{x(m_{i,j})\}$ with $i$ and $j$ arbitrary is a strictly increasing sequence.

Now $\{f(m_i) = x(m_{i})\}$ and $\{x(m_{i,j})\}$ are two increasing sequences.  Since $x$ is a bijection, and every $z\in \ZZ$ is either an $m_i$ or an $m_{i,j}$, $x$ is determined by the choice of $x(m_{1,1})$.  A valid choice for $x(m_{1,1})$ exists, since every $f$ arises as the image of some affine permutation under $Q$, and every fiber contains some $[321]$-avoiding element.  

One can show that the choice of $x(m_{1,1})$ is uniquely determined by the following argument.  Suppose two valid possibilities exist for $x(m_{1,1})$, giving rise to two different $[321]$-avoiding affine permutations $x$ and $x'$.  Suppose without loss of generality that $1\leq m_{1,1}\leq N$, and that $x(m_{1,1})<x'(m_{1,1})$.  Then: 
\begin{eqnarray*}
\binom{N+2}{2} &=&  \sum_{k=1}^{N} x(k) \\
 &=& \sum ( x(m_i) + \sum x(m_{i,j}) ) \text{ where $m_i$, $m_{i,j} \in \{1, \ldots, N\}$ }\\
 &<& \sum ( x'(m_i) + \sum x'(m_{i,j}) ) \text{ where $m_i$, $m_{i,j} \in \{1, \ldots, N\}$ }\\
 &=& \sum_{k=1}^{N} x'(k) \\
 &=& \binom{N+2}{2},
\end{eqnarray*}
providing a contradiction.  Hence $x(m_{1,1})$ is uniquely determined, and thus each fiber of $Q$ contains a unique $[321]$-avoiding permutation.
\end{proof}

%% file: chapter-crystals.tex
We give a brief introduction to crystal bases, and discuss two results.  

First, we give a complementary result to the work of Bandlow, Schilling, and Thi\'ery for two-tensors of finite-dimensional affine crystals.  Shimozono proved that in any crystal of rectangular shape of Type $A_n$, there exists a unique promotion operator which implements the Dynkin diagram automorphism in type $A_n^{(1)}$.  The promotion operator allows one to define an affine structure on the crystal; one may thus obtain the corresponding Kirillov-Reshetikhin crystal~\cite{Shimozono.2002}.
Bandlow, Schilling and Thi\'ery extended Shimozono's result, and showed that given a tensor product of two crystals of rectangular shape in type $A_n$, with $n\geq 2$, there exists a unique connected promotion operator, which in turn defines an affine structure isomorphic to the structure of the tensor product of two Kirrilov-Reshetikhin crystals~\cite{BST.2010}.  When $n=1$, there are, in general, many possible promotion operators on a tensor product of crystals which define many affine structures on the tensor product.  This occurs because the promotion operator at $n=1$ satisfies $\pr=\pr^{-1}$, thus providing less information than is available when $n\geq 2$.  However, we show (Theorem~\ref{thm:twotensorA1} that of these many possible promotion operators, only two give affine structures arising from representations of the quantum affine algebra of type $A_1^{(1)}$.

In Section~\ref{sec:stembridge}, we provide a computer implementation of Stembridge's local axioms for crystals arising from highest weight representations.

\section{ Background and Notation }
\label{sec:crystalBackground}

We first fix a number of definitions that will be useful in the subsequent sections.  Helpful books for this background are Hong and Kang's Introduction to Quantum Groups and Crystal Bases~\cite{hongKang} and Klimyk and Schm\"{u}dgen's Quantum Groups and Their Representations~\cite{klimykSchmudgen}.  The notation in this chapter is chosen to agree with Hong and Kang.

\begin{definition}
A \textbf{generalized Cartan matrix} is a square matrix $A=(a_{ij})$ such that:
\begin{itemize}
\item the diagonal entries $a_{ii}=2$,
\item the off diagonal entries $a_{ij}\leq0$,
\item $a_{ij}=0$ if and only if $a_{ji}=0$, 
\item and there exist a diagonal matrix $D$ and symmetric matrix $S$ such that $A=DS$.
\end{itemize}
\end{definition}

The generalized Cartan matrix encodes information which determines a number of important mathematical objects.  In particular, there exists a classification of quantum groups (defined below), which arise as deformations of the Universal enveloping algebra of a Kac-Moody algebra, and their associated Weyl groups by their Cartan matrices.  The matrices in this classification are labeled by a capital letter between $A$ and $G$, subscripted with the rank of the matrix $A$, possibly with a superscript.  This label is called the \textbf{Cartan type}.  (For example, $B_6$, $A_{365}^{(1)}$, and $D_4^{(3)}$ all Cartan types.)  When $A$ is of full rank, then $A$ is said to be of \textbf{finite type}.  When $A$ contains only the entries $2, 0,$ and $-1$, then $A$ is said to be \textbf{simply-laced}.  The \textbf{index set} is, for finite types, the set $I:=\{1, 2, \ldots, \rank(A) \}$, and for affine type, the set $I:=\{0, 1, 2, \ldots, \rank(A) \}$.

The information in the Cartan matrix may also be encoded in the \textbf{Dynkin diagram}, which is a graph with one node for each row of the Cartan matrix, with nodes labeled by the index set $I$.  Nodes $i$ and $j$ are connected by an edge whenever $a_{ij}\neq 0$; in simply laced types, these entries $a_{ij}$ are always $-1$, and the corresponding nodes in the Dynkin diagram are connected by a single undirected edge.  

Of particular interest in this chapter will be the cases of type $A_N$ and $A_N^{(1)}$.
For type $A_N$, the Dynkin diagram is the chain with $N$ vertices, and the Cartan matrix is:
\[
A = 
\begin{pmatrix}
 2 & -1 &  0 &  0 & \cdots & 0 & 0 \\
-1 &  2 & -1 &  0 & \cdots & 0 & 0 \\
 0 & -1 &  2 & -1 & \cdots & 0 & 0 \\
 \vdots & \vdots & \vdots & \vdots & \ddots & \vdots & \vdots \\
 0 &  0 &  0 &  0 & \cdots & 2 & -1 \\
 0 &  0 &  0 &  0 & \cdots & -1 & 2 \\
\end{pmatrix}.
\]
For the affine type $A_N^{(1)}$ with $N\geq 2$, the Dynkin diagram is a cycle with $N+1$ vertices, and Cartan matrix:
\[
A' = 
\begin{pmatrix}
 2 & -1 &  0 &  0 & \cdots & 0 & -1 \\
-1 &  2 & -1 &  0 & \cdots & 0 & 0 \\
 0 & -1 &  2 & -1 & \cdots & 0 & 0 \\
 \vdots & \vdots & \vdots & \vdots & \ddots & \vdots & \vdots \\
 0 &  0 &  0 &  0 & \cdots & 2 & -1 \\
 -1 &  0 &  0 &  0 & \cdots & -1 & 2 \\
\end{pmatrix}.
\]
Given a Dynkin diagram, one may define the associated \textbf{Weyl group} $W$.  This is a group generated by reflections $r_i$ for $i\in I$.  When $i$ and $j$ are not connected by an edge in the Dynkin diagram, then $W$ has the relation $r_ir_j=r_jr_i$.  When $i$ and $j$ are connected by a single edge (ie, $a_{ij}=-1$), then we have the relation $r_ir_jr_i=r_jr_ir_j$.  For affine types, the restriction to the finite type yields the \textbf{classical Weyl group} $W_{cl}$.  For type $A_N$, the Weyl group is isomorphic to the permutation group $S_{N+1}$. For type $A_N^{(1)}$ with $N\geq 2$ the Weyl group is the affine permutation group $\ASn_{N+1}$.  

Type $A_1^{(1)}$ is often treated as a special case, as the behavior of the associated Weyl group and quantum Kac-Moody algebra is quite distinct from that of $A_N^{(1)}$ with $N\geq 2$.  In this case $N=1$, the Weyl group $\ASn_{2}$ is isomorphic to the infinite dihedral group.  In the case where $N=1$, the Dynkin diagram for affine $A_1^{(1)}$ is given by two nodes connected by an edge labeled $\infty$ (so this is not a simply laced type), and the Cartan matrix is:
\[
A' = 
\begin{pmatrix}
 2 & -2 \\
-2 &  2 \\
\end{pmatrix}.
\]

Given the Cartan matrix, we define the \textbf{dual weight lattice} $\pv$ to be a free Abelian group of rank $2|I|-\rank A$, with basis $\{h_i \mid i \in I\} \cup \{d_s \mid s=1, \ldots, |I|-\rank A \}$.  (The basis elements $d_s$ are the \textbf{grading element(s)}; for affine types, there is one grading element, and for finite type there are none.)  Since we are only concerned with finite and affine types, we will henceforth assume that there is at most one grading element, labeled $d$.  From $\pv$, we construct the \textbf{Cartan subalgebra} $\mathfrak{h}=\CC \otimes \pv$.  The \textbf{weight lattice} $P$ for finite types is defined as the $\ZZ$-span of the set $\{\Lambda_i \in \mathfrak{h}* \mid i\in I, \Lambda_i(h_j)=\delta_{ij}, \Lambda_i(d_s)=0 \}$; the $\Lambda_i$ are the \textbf{fundamental weights}.  For affine type, there is one additional fundamental weight, the \textbf{null root} $\delta$.  In type $A_N^{(1)}$, the null root is given by $\delta = \alpha_0 + \alpha_1 + \cdots + \alpha_N$.

We then set the \textbf{simple coroots} $\Pi^\vee = \{h_i \mid i\in I\}$, and the \textbf{simple roots} 
\[
\Pi = \{\alpha_i \mid i \in I, \alpha_i \in \mathfrak{h}*, \alpha_i(h_j) = a_{ij}, \alpha_i(d)= \delta_{i,0} \},
\] 
where $a_{ij}$ is given by the Cartan matrix.

We may define the \textbf{simple reflections} $s_i$ on $\mathfrak{h}*$ by
\[
r_i(\lambda) = \lambda - \lambda(h_i)\alpha_i.
\]
The simple reflections generate a Coxeter group called the \textbf{Weyl group}.  For type $A_N$, the Weyl group is that symmetric group $S_{N+1}$, and for type $A_N^{(1)}$, the Weyl group is the affine symmetric group $\ASn_{N+1}$.

\begin{definition}
Let $A$ be a generalized Cartan Matrix such that there exists a diagonal matrix $D_{ii}=s_i$ with $DA$ symmetric.  Let $\Pi=\{\alpha_i\}$ be a collection of simple roots, $\Pi^\vee$ the simple coroots, $P$ be the weight lattice, and $\pv$ the dual weight lattice.  Then the \textbf{quantum Kac-Moody algebra} $U_q$
associated with Cartan datum
$(A, \pv, P, \Pi^\vee, \Pi)$ is the associative algebra
over $\QQ(q)$ with the unit $1$
generated by the symbols $e_i$, $f_i$ $(i \in I)$ and $q^h$
$(h\in \pv)$ subject to the following defining relations\,:
\begin{align}
\begin{aligned}
\ & q^0 = 1, \ \
q^h q^{h'} = q^{h+ h'} \quad (h, h'\in \pv),\\
\ & q^h e_i q^{-h} = q^{\alpha_i(h)} e_i, \quad
q^h f_i q^{-h} = q^{-\alpha_i(h)} f_i, \\
\ & e_i f_j - f_j e_i
= \delta_{ij} \frac{K_i - K_i^{-1}}{q_i - q_i^{-1}},
\quad \mbox{where} \quad K_i = q^{s_i h_i}, \\
\ & \sum_{r=0}^{1-a_{ij}} (-1)^r
{\begin{bmatrix} 1-a_{ij} \\ r \end{bmatrix}}_i
e_i^{1-a_{ij}-r} e_j e_i^{r} = 0 \quad \mbox{if}
\ \ a_{ii}=2, \ i \ne j, \\
\ & \sum_{r=0}^{1-a_{ij}} (-1)^r
{\begin{bmatrix} 1-a_{ij} \\ r \end{bmatrix}}_i
f_i^{1-a_{ij}-r} f_j f_i^{r} = 0 \quad \mbox{if}
\ \ a_{ii}=2, \ i \ne j, \\
\ & e_ie_j-e_je_i=0, \quad
f_if_j - f_j f_i = 0
\quad \mbox{if} \ \ a_{ij}=0.
\end{aligned}
\end{align}
Here, $K_i=q^{s_ih_i}$, where $s_i$ are the diagonal entries of the diagonal matrix $D$ in the expression $A=DS$.  In particular, when $A$ is simply-laced $s_i=1$ for every $i$.
\end{definition}

At $q=1$, this definition specializes to the enveloping algebra of a Lie algebra $\mathfrak{g}$; the quantum Kac-Moody algebra which specializes to $U(\mathfrak{g})$ at $q=1$ is denoted $U_q(\mathfrak{g})$.  For type $A_N$, this specialization is the enveloping algebra $U(\mathfrak{sl}_{N+1})$.  For affine types, we denote the Lie algebra by $\hat{\mathfrak{g}}$, and the associated Lie algebra of finite type by $\mathfrak{g})$.

On the topic of $q$-deformations, we define the quantum integer for $n\in \N$:
\[
[n]_q = \frac{q^n-q^{-n}}{q - q^{-1}},
\]
and define the $q$-factorial and $q$-binomial coefficients as:
\[
[n]_q! = \prod_{i=1}^n [i]_q, \phantom{lalala} \binom{n}{k}_q = \frac{[n]_q!}{[k]_q! [n-k]_q!}
\]
for $n, k\in \N$.
At $q=1$, these specialize to the usual integers, factorials, and binomial coefficients.  We note that the $q$-binomial coefficients satisfy the recurrence~\cite{hongKang}[Section 3.1]:
\[
\binom{n+m+1}{n}_q = q^{-n}\binom{n+m+}{n}_q + q^{m}\binom{n+m}{n-1}_q.
\]

The algebra $U_q(\mathfrak{g})$ carries a Hopf algebra structure, with comultiplication $\Delta$, counit $\epsilon$ and antipode $S$ defined as follows:
\begin{eqnarray*}
\Delta K_i = K_i \otimes K_i, & & 
    \Delta K_i^{-1} = K_i^{-1} \otimes K_i^{-1},\\
\Delta f_i = f_i \otimes 1 + K_i\otimes f_i, & &
    \Delta e_i = e_i \otimes K_i^{-1} + 1 \otimes e_i,\\
\epsilon K_i = 1, & &
    \epsilon e_i =\epsilon f_i = 0, \\
S(K_i) &=& K_i^{-1}, \\
S(e_i) = -e_i K_i^{-1}, & &
    S(f_i) =  K_i f_i. \\
\end{eqnarray*}

We may define a slightly different beast by tweaking the weight and coweight lattices.  Set $\bar{\pv}$ and $\bar{P}$ to be the $\ZZ$-spans of $\{h_i \mid i \in I\}$ and $\{\Lambda_i \mid i\in I \}$ respectively, obtained by throwing out the grading element and the null root.  The Cartan datum $(A, \bar{\pv}, \bar{P}, \Pi^\vee, \Pi)$ is called the \textbf{classical Cartan datum}, and the quantum Kac-Moody algebra generated by this set is the \textbf{quantum affine algebra} $U_q'(\mathfrak{\hat{g}})$.  This may be regarded as the subalgebra of $U_q(\mathfrak{\hat{g}})$ generated by $e_i, f_i$ and $K_i^{\pm 1}$ for $i \in I$. 

Now we describe an important class of modules for $U_q(\mathfrak{\hat{g}})$.  In what follows, we use the terms ``modules'' and ``representation'' interchangeably.  

\begin{definition}[Weight Module]
A $U_q(\mathfrak{\hat{g}})$ module $V$ is a \textbf{weight module} if it admits a \textbf{weight space decomposition}:
$V = \bigoplus_{\mu\in P}V_\mu$, where 
\[
V_\mu:=\{v\in V \mid q^hv=q^{\mu(h)}v \hspace{5pt} \forall h \in P^\vee \}.
\]
A vector $v\in V_\mu$ for some $\mu$ is called a \textbf{weight vector}, in which case $\mu$ is the \textbf{weight} of $v$,denoted $\wt(v)=\mu$.

A weight module $V$ is a \textbf{highest weight module} if there exists a vector $v_\lambda\in V$ of weight $\lambda$ such that $v_\lambda$ generates $M$ as a $U_q(\mathfrak{\hat{g}})$ module, and $e_iv_\lambda=0$ for all $i\in I.$
\end{definition}

We define a partial order on the weight lattice $P$ by $\lambda \geq \mu$ if and only if $\lambda-\mu \in \bigoplus_{i\in I} \ZZ_{\geq 0} \alpha_i$.
For a weight $\lambda \in P$, set $D(\lambda)$ to be the set of $\mu \leq \lambda$.

\begin{definition}
A $U_q(\mathfrak{\hat{g}})$-module $V$ is in \textbf{category $\mathcal{O}^q_{int}$} if:
\begin{itemize}
\item $V$ is a weight module, with $V=\bigoplus_{\mu\in P}V_\mu$ and $\dim V_\mu<\infty$ for every $\mu$,
\item There exist a finite number of weights $\lambda_1, \ldots, \lambda_k$ such that for every weight vector $v\in V$, $\wt(v) \in D(\lambda_1) \cup \cdots \cup D(\lambda_k)$,
\item All $e_i$ and $f_i$ are locally nilpotent on $V$.
\end{itemize}
Modules in $\mathcal{O}^q_{int}$ are also called \textbf{integrable} modules.
\end{definition}

We now construct the Kashiwara operators, $\tilde{e}_i$ and $\tilde{f}_i$, following the presentation of~\cite{hongKang}.  The divided operators $f^{(n)}=\frac{f^n}{[n]_q!}$ are important for the crystal basis theory and the definition of the Kashiwara operators.  

\begin{lemma}
Let $V=\bigoplus_{\mu\in P}V_\mu$ be a weight module for $U_q(\mathfrak{\hat{g}})$.  Then, for each $i\in I$, every weight vector $u\in V_\mu$ may be uniquely expressed as
\[
u = u_0 + f_i u_1 + f_i^{(2)}u_2 + \cdots + f_i^{(k)}u_k,
\]
where $k\in \ZZ_{\geq 0}$ and $u_k \in V_{\mu + k\alpha_i} \cup \ker e_i$.
\end{lemma}
The proof of this lemma may be found in~\cite[Chapter 4]{hongKang}.

\begin{definition}
The \textbf{Kashiwara operators} $\tilde{e}_i$ and $\tilde{f}_i$ on $V$ are defined on a weight vector $u$ by:
\[
\tilde{e}_i u = \sum_{k=1}^N f_i^{(k-1)} u_k, \phantom{blech.} \tilde{f}_i u = \sum_{k=0}^N f_i^{(k+1)} u_k.
\]
\end{definition}
In particular, $\wt( \tilde{e_i}u) = \wt(u)+\alpha_i$, and $\wt( \tilde{f_i}u) = \wt(u)-\alpha_i$.

In the construction of the crystal basis, we first construct the 
crystal lattice.  Set $A_0$ to be the ring of rational functions in $q$ that evaluate at $q=0:$ 
\[
 A_0=\{ \frac{f}{g} \mid f,g\in \CC[q], g(0) \neq 0 \}.
\]

\begin{definition} Let $M$ an integrable $U_q'(\mathfrak{\hat{g}})$ module.  A free $A_0$-submodule $\mathcal{L}$ of $M$ is a \textbf{crystal lattice} if:
\begin{itemize}
\item $\mathcal{L}$ generates $M$ as a vector space over $F(q)$
\item $\mathcal{L}$ decomposes as a direct sum of weight spaces, compatible with those of $M$.
\item $\tilde{e}_i\mathcal{L} \subset \mathcal{L}, \tilde{f}_i\mathcal{L} \subset \mathcal{L}$ for all $i \in I.$
\end{itemize}
\end{definition}

We can realize the crystal basis as a $\mathbb{Q}$-basis of $\mathcal{L}\setminus q \mathcal{L}$, which is often thought of crudely as a limit as $q$ goes to $0$. 

\begin{definition}
\label{def:crystalBasis}
A \textbf{crystal basis} for $M$ is a pair $(L,B)$ such that:
\begin{enumerate}
\item $\LL$ is a crystal lattice for $M.$
\item $B$ is a $\mathbb{Q}$-basis of $\LL\setminus q\LL \cong \mathbb{Q} \otimes_{A_0} \LL$
\item The elements of $B$ are weight vectors.
\item $\tilde{e_i}B\subset B\bigcup {0}, \tilde{f_i}B\subset B\bigcup {0}$ for all $i\in I.$
\item For all $i\in I, b,b'\in B,$ we have $\tilde{f_i}b=b'$ iff $b=\tilde{e_i}b'$
\end{enumerate}
\end{definition}

In particular, we see that the Kashiwara operators must preserve the crystal basis as a set; this will have consequences below.

As an abstraction of crystal bases, we have crystal graphs, which we define as follows.
\begin{definition}
\label{def:crystalGraph}
A \textbf{crystal graph} associated to the Cartan datum $(A, \pv, P, \Pi^\vee, \Pi)$ is a set $\mathcal{B}$ and maps $\wt: \mathcal{B} \rightarrow P $, $\tilde{e_i}, \tilde{f_i}: \mathcal{B} \rightarrow \mathcal{B} \cup \{ \emptyset \}$, and $\epsilon_i, \phi_i: \mathcal{B} \rightarrow \ZZ_{\geq 0}$ (with $i\in I$), satisfying the following properties for all $i\in I, b, b' \in \mathcal{B}$:
\begin{enumerate}
\item $\phi_i(b) = \epsilon_i(b) + \langle h_i, \wt(b) \rangle $,
\item $\wt(\tilde{e_i} b) = \wt(b)+\alpha_i $, $\wt(\tilde{f_i} b) = \wt(b)-\alpha_i $ when $\tilde{e_i} b, \tilde{f_i} b \neq \emptyset$,
\item $\epsilon_i(\tilde{e_i} b) = \epsilon_i(b)-1 $, $\phi_i(\tilde{e_i} b) = \phi_i(b)+1 $ when $\tilde{e_i} b \neq \emptyset$,
\item $\epsilon_i(\tilde{f_i} b) = \epsilon_i(b)+1 $, $\phi_i(\tilde{f_i} b) = \phi_i(b)-1 $ when $\tilde{f_i} b \neq \emptyset$,
\item (Semiregular) $\epsilon_i(b) = \operatorname{max}\{k\geq 0 \mid \tilde{e_i}^k b\in \mathcal{B} \} $,
    $\phi_i(b) = \operatorname{max}\{k\geq 0 \mid \tilde{f_i}^k b\in \mathcal{B} \} $,
\item $\tilde{f_i} b'=b$ if and only if $b'= \tilde{e_i} b$.
\end{enumerate}
We say that $\mathcal{B}$ is a $U_q(\mathfrak{g})$-crystal, where $U_q(\mathfrak{g})$ is the quantum Kac-Moody algebra associated to the Cartan datum $(A, \Pi, \Pi^\vee, P, \pv)$.  We set $\mathcal{B}_\lambda = \{b \in \mathcal{B} \mid \wt(b)=\lambda \}$.
\end{definition}

In effect, a crystal graph is an edge-colored directed graph with the map $\wt$ to the weight lattice.  By the last property, each vertex in the crystal graph has at most in-degree and out-degree $1$ in any particular color.  Consider the sub-graph with the same vertex set but with only the $i$-colored edges.  This graph is a collection of directed chains, called $i$-strings.  

The crystal graph admits a \textbf{Weyl group action} by $W$, the Weyl group associated to the Cartan matrix $A$.  The Weyl group is generated by reflections $r_i$ for $i \in I$; the action of $r_i$ on an element $b$ is to ``flip'' $b$ in its $i$-string.  More precisely, if the $i$-string of $b$ has $n$ vertices, and 
$\tilde{e_i}^kb=\emptyset$, then $r_i(b)$ is the vertex $b'$ in the same $i$-string satisfying $\tilde{f_i}^kb'=\emptyset$.

Many different constructions of crystal graphs exist, explicitly designed to give graphs which arise from crystal bases for actual $U_q(\mathfrak{g})$ representations.  One such construction is given by \textbf{crystals of tableaux}, which we will now partially describe.  

A \textbf{Ferrer's Diagram} of a partition $\lambda=(\lambda_1, \lambda_2, \ldots, \lambda_k)$ is a collection of ``boxes shoved in a corner,'' with $k$ rows and $\lambda_i$ boxes in row $i$, and all rows ``left-justified.''  Thus, the total number of boxes is $n$ if $\lambda$ is a partition of $n$.  Also, since $\lambda$ is a decreasing sequence, the boxes appear to be ``bottom-justified.''  (Up to one's choice of bottom, anyway.  The French, in a rare display of practicality, place the bottom at the bottom of the page.  The English, known for being wily and inscrutable, place the bottom at the top.  The Russians, being clever, rotate the French convention counter-clockwise by $45$ degrees, so that the ``corner'' is at the bottom of the page and all of the justification happens by gravity.  Being neither wily nor clever, the author will stick with the French convention.)  

A \textbf{Young tableaux} $T$ is an assignment of a positive integer to each box in a Ferrer's diagram $\lambda$, which is called the \textbf{shape} of the tableaux.  $T$ is \textbf{semi-standard} if:
\begin{itemize}
\item $T$ has strictly increasing columns, read bottom-to-top, and 
\item $T$ has weakly-increasing rows, read left-to-right.
\end{itemize}

\begin{figure}
\begin{tikzpicture}[scale=.5]
\coordinate (base point) at (0,12);
\node (2211) at (1,13) [rectangle, draw, minimum size=1cm] {};
\foreach \entry / \r / \c in { 2/1/0, 2/1/1, 1/0/0, 1/0/1  }
{ \draw (base point) ++(\c,\r) rectangle +(1,1) node [midway] {\entry}; }

\coordinate (base point) at (0,9);
\node (2311) at (1,10) [rectangle, draw, minimum size=1cm] {};
\foreach \entry / \r / \c in { 2/1/0, 3/1/1, 1/0/0, 1/0/1  }
{ \draw (base point) ++(\c,\r) rectangle +(1,1) node [midway] {\entry}; }

\coordinate (base point) at (-3,6);
\node (2312) at (-2,7) [rectangle, draw, minimum size=1cm] {};
\foreach \entry / \r / \c in { 2/1/0, 3/1/1, 1/0/0, 2/0/1  }
{ \draw (base point) ++(\c,\r) rectangle +(1,1) node [midway] {\entry}; }

\coordinate (base point) at (3,6);
\node (3311) at (4,7) [rectangle, draw, minimum size=1cm] {};
\foreach \entry / \r / \c in { 3/1/0, 3/1/1, 1/0/0, 1/0/1  }
{ \draw (base point) ++(\c,\r) rectangle +(1,1) node [midway] {\entry}; }

\coordinate (base point) at (0,3);
\node (3312) at (1,4) [rectangle, draw, minimum size=1cm] {};
\foreach \entry / \r / \c in { 3/1/0, 3/1/1, 1/0/0, 2/0/1  }
{ \draw (base point) ++(\c,\r) rectangle +(1,1) node [midway] {\entry}; }

\coordinate (base point) at (0,0);
\node (3322) at (1,1) [rectangle, draw, minimum size=1cm] {};
\foreach \entry / \r / \c in { 3/1/0, 3/1/1, 2/0/0, 2/0/1  }
{ \draw (base point) ++(\c,\r) rectangle +(1,1) node [midway] {\entry}; }

	\draw [->>] (2211.south) to (2311.north) [red] node [auto,midway] {2};
	\draw [->>] (2311.east) to [out=0, in=90] (3311.north) [red] node [auto,midway] {2};
	\draw [->>] (2312.south) to [out=270, in=180] (3312.west) [red] node [auto,midway] {2};
	\draw [->>] (2311.west) to [out=180, in=90] (2312.north) [blue] node [auto,midway] {1};
	\draw [->>] (3311.south) to [out=270, in=0] (3312.east) [blue] node [auto,midway] {1};
	\draw [->>] (3312.south) to (3322.north) [blue] node [auto,midway] {1};

\end{tikzpicture}
\caption{Example of a crystal graph of type $A_2$.  The vertices are given by semi-standard Young tableaux, and the edges are determined by a combinatorial algorithm.  The weight of a tableaux $T$ is the vector whose $i$th entry is the number of $i$'s appearing in $T$.}
  \label{fig.crystalTableaux}
\end{figure}
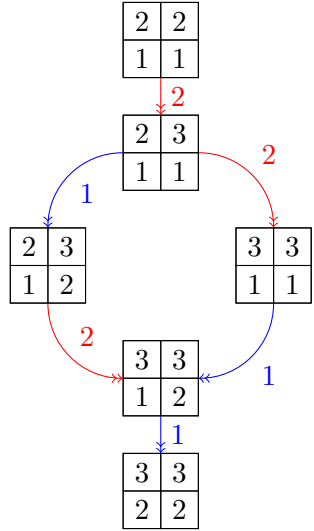

The \textbf{crystal graph of type $A_N$ of shape $\lambda$}, where $\lambda$ is a partition, is an edge-colored digraph whose vertex set is the set of all semi-standard Young tableaux of shape $\lambda$ whose entries are in the set $\{1, 2, \ldots, N+1\}$.  The edges are determined by a combinatorial algorithm, described in~\cite{hongKang}.  The weight of a tableaux $T$ is simply the vector whose $i$th coordinate is the number of $i$'s appearing in $T$.  An example is shown in Figure~\ref{fig.crystalTableaux}, which depicts the crystal of tableaux of highest weight $\lambda=(2,2)$ in type $A_2$.

We say that a crystal is \emph{rectangular} if it is a highest-weight crystal with highest weight $r \Lambda_s$ for $r\in \N$ and $s\in I$.  This corresponds to a crystal of tableaux of shape $\lambda$ where $\lambda$ is a rectangle of height $s$ and width $r$.

\subsection{Tensor Products}
Modules for the quantum group $U_q(\mathfrak{\hat{g}})$ admit natural tensor products, thanks to the Hopf algebra structure on $U_q(\mathfrak{\hat{g}})$.  Tensor products of modules with crystal bases also admit crystal bases, structured according to the Tensor Product Rule, described in~\cite{hongKang}[Theorem 4.4.1].

\begin{theorem}[Tensor Product Rule.]
Let $U_1$ and $U_2$ be $U_q(\mathfrak{\hat{g}})$ modules, with crystal lattices $\LL_1$ and $\LL_2$  and crystal bases $B_1$ and $B_2$.  Set $\LL= \LL_1 \otimes_{A_0} \LL_2$ and $B= B_1 \times B_2$.  Then $\LL$ and $B$ are a crystal lattice and crystal basis for $U_1\otimes U_2$, where the action of the Kashiwara operators is given by:
\begin{eqnarray*}
   \tilde{e_i}(b_1\otimes b_2) &=& \left\{
     \begin{array}{lr}
       (\tilde{e_i}b_1)\otimes b_2    & : \text{ if } \phi_i(b_1)\geq \epsilon_i(b_2), \\
        b_1\otimes (\tilde{e_i}b_2)   & : \text{ if } \phi_i(b_1)<    \epsilon_i(b_2),\\
     \end{array}
   \right\}, \\
   \tilde{f_i}(b_1\otimes b_2) &=& \left\{
     \begin{array}{lr}
       (\tilde{f_i}b_1)\otimes b_2    & : \text{ if } \phi_i(b_1)>    \epsilon_i(b_2) \\
        b_1\otimes (\tilde{f_i}b_2)   & : \text{ if } \phi_i(b_1)\leq \epsilon_i(b_2).\\
     \end{array}
   \right\}
\end{eqnarray*}
Therefore we have:
\begin{eqnarray*}
\wt(b_1\otimes b_2) &=& \wt(b_1) + \wt(b_2), \\
\epsilon_i(b_1\otimes b_2) &=& \text{max}( \epsilon_i(b_1),  \epsilon_i(b_2)-\langle h_i, \wt(b_1) \rangle ), \\
\phi_i(b_1\otimes b_2) &=& \text{max}( \phi_i(b_2),  \epsilon_i(b_1)-\langle h_i, \wt(b_2) \rangle ).
\end{eqnarray*}
We write $b_1\otimes b_2$ for the element $b_1 \times b_2$ and understand that $b_1\otimes \emptyset = \emptyset \otimes b_2 = \emptyset$.
\end{theorem}

The tensor product rule descends naturally to crystal graphs.

\begin{theorem}[Tensor Product Rule.]
Let $B_1$ and $B_2$ be crystal graphs.  Then $B= B_1 \otimes B_2$ is a crystal graph with vertex set $B_1 \times B_2$ whose crystal structure is defined by:
\begin{eqnarray*}
    \wt(b_1\otimes b_2) &=& \wt(b_1) + \wt(b_2), \\
    \epsilon_i(b_1\otimes b_2) &=& \text{max}( \epsilon_i(b_1),  \epsilon_i(b_2)-\langle h_i, \wt(b_1) \rangle ), \\
    \phi_i(b_1\otimes b_2) &=& \text{max}( \phi_i(b_2),  \epsilon_i(b_1)-\langle h_i, \wt(b_2) \rangle ), \\
   \tilde{e_i}(b_1\otimes b_2) &=& \left\{
     \begin{array}{lr}
       (\tilde{e_i}b_1)\otimes b_2    & : \text{ if } \phi_i(b_1)\geq \epsilon_i(b_2), \\
        b_1\otimes (\tilde{e_i}b_2)   & : \text{ if } \phi_i(b_1)<    \epsilon_i(b_2),\\
     \end{array}
   \right\}, \\
   \tilde{f_i}(b_1\otimes b_2) &=& \left\{
     \begin{array}{lr}
       (\tilde{f_i}b_1)\otimes b_2    & : \text{ if } \phi_i(b_1)>    \epsilon_i(b_2) \\
        b_1\otimes (\tilde{f_i}b_2)   & : \text{ if } \phi_i(b_1)\leq \epsilon_i(b_2).\\
     \end{array}
   \right\}
\end{eqnarray*}
We write $b_1\otimes b_2$ for the element $b_1 \times b_2$ and understand that $b_1\otimes \emptyset = \emptyset \otimes b_2 = \emptyset$.
\end{theorem}

Using the tensor product rule, one may show that the tensor product of two crystals $B(k)$ and $B(j)$ of type $A_2$ is 
\[
B(k)\otimes B(j) \cong \bigoplus_{t=0}^{\lfloor \frac{k+j}{2} \rfloor}  B(t).
\]

\subsection{Promotion Operators}

Let $B$ be a classical crystal of type $A_N$.  The Dynkin diagram of type $A_N$ is a chain, but by adjoining a new node and attaching edges to make a cycle, we may obtain the Dynkin diagram of type $A_N^{(1)}$.  A promotion operator is a map $B\rightarrow B$ which implements the affine Dynkin diagram automorphism on $B$, allowing one to place $0$ arrows in $B$ and obtain a crystal of type $A_N^{(1)}$.

\begin{definition}
We define a \textbf{promotion operator} to be a map $\pr: B\rightarrow B$ satisfying the properties:
\begin{itemize}
\item If $\wt(b) = (w_1, \ldots, w_{N+1})$, then $\wt(\pr(b)) = (w_{N+1}, w_1, \ldots, w_N)$,
\item $\pr^{N+1} = \operatorname{id}$, and 
\item For all $i \in \{1, \ldots, N \}$, we have: 
\[
    \pr \circ f_i = f_{i+1}\circ \pr , \text{ and } \pr \circ e_i =e_{i+1} \circ \pr.
\]
\end{itemize}
Given a promotion operator on a classical crystal $B$, one may define an affine structure on $B$ by placing the $0$-arrows according to:
\[
    \pr \circ f_n = f_{0}\circ \pr , \text{ and } \pr \circ e_n =e_{0} \circ \pr.
\]
A promotion operator is \textbf{connected} if the resulting affine crystal is connected.
\end{definition}

Recall that Bandlow, Schilling, and Thi\'ery studied two-tensors of crystals of rectangular shape in type $A_N$ with $N\geq 2$.  Their result showed that in such a crystal admits a unique connected promotion operator~\cite{BST.2010}.

In the case where $N=1$, a promotion operator satisfies $\pr^2=\operatorname{id}$, so that $\pr = \pr^{-1}$.  In practice, this degeneracy means that the promotion operator at $N=1$ provides less information than the case where $N\geq 2$.  In fact, more than one connected promotion operators may be found in a two tensor of rectangular $A_1$ crystals.  Below, we show that there is a unique connected promotion operator on such a two-tensor which gives an affine structure actually coming from a representation of $U_q'(\slnh_2)$.

\subsection{Evaluation Modules}
From a result of Chari and Pressley~\cite{chariPressley.95}[Theorem 4.3], it is known that every finite dimensional irreducible weight module of $U_q(\hat{sl_2})$ is isomorphic to a tensor product of evaluation representations of $U'_q(\hat{sl_2})$.

The irreducible representation $V(j)$ is the irreducible $j+1$-dimensional representation of $U_q(sl_2)$ generated by $u$, with basis $\{ u, f_1u, \ldots , f_1^ju \}$.  Then $V(j)$ carries a $U'_q(\hat{sl_2})$ structure defined by the relations:
\begin{eqnarray*}
	K_0f_1^xu & = & q^{-j+2x}f_1^xu, \\
	e_0f_1^xu & = & f_1^{x+1}u, \\
	f_0f_1^xu & = & [x]_q[j-x+1]_qf_1^{x-1}u.
\end{eqnarray*}

The evaluation module $V(j)_a$ modifies this structure by the introduction of a deformation by a constant $a$ as follows:
\begin{eqnarray*}
	K_0f_1^xu & = & q^{-j+2x}f_1^xu \\
	e_0f_1^xu & = & af_1^{x+1}u \\
	f_0f_1^xu & = & a^{-1}[x]_q[j-x+1]_qf_1^{x-1}u
\end{eqnarray*}
The classical structure remains the same as in $V(j)$.

\subsection{Simple Crystals and Extremal Vectors}
\begin{definition}
Call an element $b$ of a crystal $B$ for $U_q(g)$ with associated Weyl group $W$ \textbf{extremal} if:
\begin{enumerate} 
\item Either $e_i(b)=\emptyset$ or $f_i(b)=\emptyset$ for every $i$, and 
\item For every $w\in W$, $wX$ satisfies condition $(1)$.  
\end{enumerate}
A crystal is \textbf{simple} if there exists some $b$ in $B$ such that any extremal vector of $B$ is contained in $W_{cl}b$.  
\end{definition}

In particular, simplicity implies that for any finite-dimensional irreducible representation $M$ of $U'_q(\hat{sl_2})$, the crystal for $M$ can have at most two extremal vectors.  In $V(j)_a\otimes V(k)_b$ these vectors are already spoken for: since $f_0$ and $e_1$ increase weights and $e_0$ and $f_1$ decrease weights, the vectors $u\otimes v$ and $f_1^{(j)}u\otimes f_1^{(k)}v$ must be the only extremal vectors.

\section{ Classifying Finite Dimensional Crystals for $U'_q(\hat{sl_2})$ }

Let $V(j)$, $V(k)$ be two irreducible representations of $U_q(sl_2)$, generated by highest weight vectors $u$ and $v$ with highest weights $j$ and $k$, respectively. Then the the coproduct on $U_q(sl_2)$ determines the structure of the tensor product $U\otimes V$.  Recall that the coproduct of $f_1$ is $\Delta(f_1) = f_1\otimes 1 + K_1\otimes f_1$.  Define $p_{x,y}(q)$ such that $f_1^n(u\otimes v) = \sum_{x+y=n} p_{x,y}(q)f_1^x(u)\otimes f_1^y(v)$.  So $p_{x,y}(q)$ are the structure constants for the action of $U_q(sl_2)$ on $U\otimes V$.  

To get some sense of $p_{x,y}(q)$, consider the integer point $(x, y)$ in the plane, and notice that $p_{x,y}$ is a $q$-counting of all possible ways to get from $(0,0)$ to $(x, y)$ incrementing by either $(1, 0)$ or $(0, 1)$.  There are $\binom{x+y}{x}$ such paths, weighted by $q^n$, where $n$ is determined by the particular path taken.  For example, one can get from $u\otimes v$ to to $f_1^1(u)\otimes f_1^1(v)$ in two ways: by going left then up or by going up then left.  

We now consider the weight of a path heuristically.  (Below, in Lemma~\ref{lem:structConst} an explicit recurrence will be detailed.)  Each increase in $y$ applies $K_1$ to $f_1^\cdot(u)$, thus changing the power of $q$ associated to the path, while increasing $x$ does not change the exponent.  Keeping track solely of the action of $f_1$'s and $K_1$'s on $u$, to each path we associate a word in the letters $f$ and $K$.  Any path to $f_1^x(u)\otimes f_1^y(v)$ will have a word with $x$ $f_1$'s and $y$ $K_1$'s.  There are $\binom{x+y}{x}$ such paths, each contributing some power of $q$.  So we can tell at the outset that the evaluation of $p$ at $q=1$ will be $\binom{x+y}{x}$.  Also notice that $p_{x,0} = 1$ for all $x$, and $p_{0,y} = q^{jy}$ for all y, since the weight of the words $f^x$ and $K_1^y$ are those values, respectively.  In light of all this, the following lemma should not be too surprising.

\begin{lemma} 
\label{lem:structConst}
$p_{x,y}(q) = q^{jy-xy} {x+y\choose x}_q$, where the symbol on the right is the quantum binomial coefficient given by $[x+y]_q!/([x]_q![y]_q!)$. 
\end{lemma}

\begin{proof}  There is a simple recurrence on $p_{x,y}(q)$, analogous to the recurrence on the binomial coefficients.  Any word in $x$ $f_1$'s and $y$ $K_1$'s begins with either an $f_1$ or a $K_1$.  Strip away the first letter to get a path contributing to either $p_{x-1,y}(q)$ or $p_{x,y-1}(q)$.  Adding an $f$ to a path to $p_{x-1,y}(q)$ does not change the weight of the path.  Recalling that $K_1f_1=q^{-2}f_1K_1$, we can see that adding a $K_1$ to a path to $p_{x,y-1}(q)$ will change the weight by $q^{-2x}$, since the relation must be applied $x$ times.  Thus, $p_{x,y}(q) = p_{x-1,y}(q) + q^{-2x}p_{x,y-1}(q)$.

There is a similar recurrence on the quantum binomial coefficients.  In particular, ${x+y+1\choose x}_q = q^{-x}{x+y\choose x}_q + q^y{x+y\choose x-1}_q$.  

Notice that the Lemma holds in the boundary cases $x=0$ or $y=0$.  For the induction, set $p_{x-1,y}(q) = q^{jy-(x-1)y} {x+y-1\choose x-1}_q$ and $p_{x,y-1}(q) = q^{j(y-1)-x(y-1)} {x+y-1\choose x-1}_q$.  Plug these values into the recurrence on $p_{x,y}$ to complete the proof.
\end{proof}

We have the following corollary.
\begin{corollary} 
\label{cor:structConst}
$f_1^{(n)}(u\otimes v) = \sum_{x+y=n} q^{jy-xy} f_1^{(x)}(u)\otimes f_1^{(y)}(v)$.  Furthermore, at $n=j+k$, this reduces to 
\[
    f_1^{(j+k)}(u\otimes v) = f_1^{(j)}u \otimes f_1^{(k)}v.
\]
\end{corollary}
\begin{proof}
The first computation is immediate from the lemma and the definition of the divided difference operators.

At $n=j+k$, there is only one additive term, since $f_1^{j+1}u=f^{k+1}v=0$.  The second identity then follows from evaluation of the first identity at $x=j, y=k$.
\end{proof}

Now we wish to show that the combinatorial promotion operator (described in~\cite{Shimozono.2002, BST.2010}) yields the only possible connected affine structure on $V(j)\otimes V(k)$ arising from a weight representation of $U'_q(\hat{sl_2})$.  From the Chari and Pressley's theorem~\cite{chariPressley.95}[Theorem 4.3], we can identify any such affine structure as a tensor of evaluation modules, so we examine $V(j)_a\otimes V(k)_b$.  Chari and Pressley further showed that these evaluation modules are almost always finite-dimensional; we consider the question of when these modules have a crystal basis.

The crystal for the evaluation module $V(j)_a\otimes V(k)_b$ at $a=b=1$ is the tensor product of the Kirillov-Reshetikhin crystals $B(j)$ and $B(k)$, as explained in~\cite[Chapter 10]{hongKang}.  $B(j)\otimes B(k)$ carries an affine structure corresponding to the canonical promotion operator.  We also know from the Littlewood-Richardson rule that the underlying $U_q(sl_2)$ structure is a direct sum of $V(j+k-2i)$ for $i$ in $0$ to 
$\lfloor \frac{j+k}{2} \rfloor$.  In particular, there is one vector of weight $\pm(j+k)$, two of weight $\pm(j+k-2)$, and so on.

\begin{example}
Consider for a moment the case of only one tensor factor, the evaluation module $V(j)_a$ with highest weight vector $u$.  The underlying classical module $V(j)$ has a crystal lattice $L$ and a crystal basis $B$ which are unique up to scalar multiplication; thus, this basis must be equal (as a set) to the crystal basis for $V(j)_a$.   For the classical crystal basis we take: 
\[
B = \{f_1^{(x)}u \mid x \in \{0, 1, \ldots, j\}\}.
\]
Other bases my be obtained by multiplying every element of $B$ by the same constant multiple; multiplying by such a scalar will not change the crystal structure, and does not affect our further arguments.

Then, using the definition of the evaluation module, we have:
\[
	e_0^{(y)} u = a^y f_1^{(y)}u.
\]
We note that $a$ must be in $A_0$ in order to preserve the lattice; otherwise, $e_0(u)$ would lie outside the crystal lattice violating the definition of the crystal lattice.

Consider the choice $a=1+q$.  Then $e_0(u) = a f_1(u) = f_1(u) + q f_1(u) = f_1(u)$, where the last equality is as elements of $\LL \setminus q\LL$.  Then this choice of $a$ yields is consistent with the classical crystal basis.  However, the choice $a=2$ would give $e_0(u) = 2 f_1(u)$, which is not in the classical crystal basis.

Then one can see that any choice of $a$ such that $V(j)_a$ has a crystal basis must have $a(0)=1$.
\end{example}

\begin{lemma}
Let $B$ be an affine crystal structure on $V(j)\otimes V(k)$ arising from a promotion operator $\pr$.  Then $\tilde{e_0}(u\otimes v)$ and $\tilde{f_0}\tilde{f_1}^{(j+k)}(u\otimes v)$ are non-empty in $B$.
\end{lemma}
\begin{proof}
Consider that $\tilde{f_0}^p \pr(u\otimes v)= \pr \tilde{f_1}^p (u\otimes v) \neq \emptyset$ for $0\leq p \leq j+k$, and $\tilde{e_0} \pr(u\otimes v)= \pr \tilde{e_1} (u\otimes v) = \emptyset$, so that $\pr(u\otimes v)$ is a $0$-highest weight vector, sitting in a $0$-string of length $j+k$.  By weight considerations, we can then conclude that $\pr(u\otimes v)=\tilde{f_1}^{(j+k)}(u\otimes v)$.  Since $\pr^2=id$, we also have $u\otimes v=\pr \tilde{f_1}^{(j+k)}(u\otimes v)$.

As a result, we have:
\begin{eqnarray*}
\tilde{e_0}(u\otimes v) &=& \tilde{e_0}\pr^2((u\otimes v)) \\
        &=& \pr \tilde{e_1} \pr (u\otimes v)) \\
        &=& \pr \tilde{e_1} \tilde{f_1}^{(j+k)}(u\otimes v) \\
        &=& \pr \tilde{f_1}^{(j+k-1)}(u\otimes v) \neq \emptyset. \\
\end{eqnarray*}
A similar calculation shows that $\tilde{f_0}\tilde{f_1}^{(j+k)}(u\otimes v)\neq \emptyset$.
\end{proof}

\begin{corollary}
We have:
\begin{eqnarray*}
\tilde{e_0}(u\otimes v) &=& \tilde{f_1}u \otimes v \text{ or } u \otimes \tilde{f_1}v, \\
\tilde{f_0}(\tilde{f_1}^{(j+k)}(u\otimes v)) &=& 
    \tilde{f_1}^{j-1}u \otimes \tilde{f_1}^k v \text{ or } \\
    &   & \tilde{f_1}^{j}u \otimes \tilde{f_1}^{k-1} v.
\end{eqnarray*}
\end{corollary}
\begin{proof}
This follows directly from the existence of these elements and weight considerations in the classical crystal $V(j)\otimes V(k)$.
\end{proof}

We now obtain restrictions on what choices of $a$ and $b$ may give rise to a crystal structure with promotion operator.

\begin{proposition}
\label{prop:crysStructure}
Let $B$ be an affine crystal structure on $V(j)_a\otimes V(k)_b$ arising from a promotion operator $\pr$.  Then one of two cases hold:
\begin{itemize}
\item Either $a(0)=b(0)=1$, which recovers the tensor product $V(j)\otimes V(k)$ which gives rise to the usual tensor of finite affine crystals,

\item Or we have $b = q^j b'$ and $a=q^{-k} a'$, with $a', b' \in A_0$ and $a'(0)=b'(0)=1$.
\end{itemize}
\end{proposition}
\begin{proof}

We now consider the ``top'' and ``bottom'' of the crystal individually, and obtain restrictions on $a$ and $b$.

First, we consider the action of $e_0$ at the top of the crystal, applying it to the classical highest weight element of extremal weight.  Since the crystal basis $B \cup \{ \emptyset \}$ is preserved as a set by the actions of $\tilde{e_0}$ and $\tilde{f_0}$, we have:
\begin{eqnarray*}
 \tilde{e_0}(u\otimes v) &=& (e_0 u) \otimes K_0^{-1}v + u \otimes (e_0 v) \\
                         &=& aq^{k} f_1(u) \otimes v +   b u \otimes f_1(v) \\
                         &=&         f_1(u) \otimes v + q\LL \\
                         & & \text{ or }        u \otimes f_1(v) + q\LL
\end{eqnarray*}
Then there are two possibilities for the evaluation of $a$ and $b$ at $q=0$.  The first possibility is that $b(0)=1$ and $a(0)$ is chosen in such a way that $aq^{k} f_1(u) \otimes v \in q\LL$.  (Note that $a(0)=1$ satisfies this possibility.)  The second possibility is that or $a=q^{-k} a'$, with $a'(0)=1$,  and $b$ is chosen such that $b u \otimes f_1(v) \in q\LL$.

Next, consider the action of $f_0$ at the bottom of the crystal:
\begin{eqnarray*}
f_0^{(1)}(f_1^{(j)}u\otimes f_1^{(k)}v) 
    &=&  f_0 f_1^{(j)}u\otimes f_1^{(k)}v  + K_0 f_1^{(j)}u\otimes f_0 f_1^{(k)}v \\
    &=&  a^{-1} f_1^{(j-1)}u\otimes f_1^{(k)}v + b^{-1} q^{j} f_1^{(j)}u\otimes f_1^{(k-1)}v \\
    &=&               f_1^{(j-1)}u\otimes f_1^{(k)}v + q\LL \\
    & & \text{ or }                                                 f_1^{(j)}u\otimes f_1^{(k-1)}v + q\LL.
\end{eqnarray*}
Again, there are two possibilities.  In the first possibility, we have $a(0)=1$ and $b(0)$ chosen such that 
$b^{-1} q^{j} f_1^{(j)}u\otimes f_1^{(k-1)}v \in q\LL$.  (For this, the choice $b(0)=1$ suffices.)  In the second possibility, we take $b = q^j b'$ where $b'(0)=1$, and $a$ chosen such that 
$a^{-1} f_1^{(j-1)}u\otimes f_1^{(k)}v \in q\LL$.  

One can see that there are only two consistent choices of possibilities at the top and bottom of the crystal yield the cases outlined in the statement of the proposition.  

\end{proof}

These choices of $a$ and $b$ completely determine the crystal structure; thus, there are at most two possible affine crystal structures on $V(j)_a \times V(k)_b$.  Re-applying the arguments in~\ref{lem:structConst}, we may obtain the following result, helpful in describing the $0$-string through the extremal weight vector $u\otimes v$.

\begin{lemma}
\label{lem:e0StructConsts}
The structure constants of $e_0^{(n)}$ on the vector $u\otimes v$ are given by:
\[
e_0^{(n)}(u\otimes v) = \sum_{x+y=n} q^{kx-xy}a^x b^y f_1^{(x)}(u)\otimes f_1^{(y)}(v).
\]
Furthermore, at $n=j+k$, this reduces to:
\[
    e_0^{(j+k)}(u\otimes v) = a^j b^k f_1^{(j)}u \otimes f_1^{(k)}v.
\]
\end{lemma}
\begin{proof}
Notice that each application of $e_0$ to the left side of the tensor contributes an $a$ and each application on the right contributes a $b$, explaining the $a^x b^y$ in the formula.  Otherwise, the result follows from an inductive argument on quantum binomial coefficients exactly analogous to that in Lemma~\ref{lem:structConst}.
\end{proof}

In the case where $a(0)=b(0)=1$, the leading term of $e_0^{(n)}(u\otimes v)$ is that with $x=0$ (if $n<k$) or $y=k$ (if $n\geq k$).  These correspond in the crystal limit exactly to the highest/lowest weight vectors in each $1$-string in the tensor product, showing that the resulting crystal is connected.

In the case where $a=q^{-k}$ and $b=q^j$, we obtain:
\[
   e_0^{(n)}(u\otimes v) = \sum_{x+y=n} q^{jy-xy} f_1^{(x)}(u)\otimes f_1^{(y)}(v) = f_1^{(n)}(u\otimes v).
\]
This means that the $0$-string and $1$-string in the crystal limit coincide, and the crystal is disconnected.  (We may also use $a=q^{-k}a'$ and $b=q^{-k}b'$ with $a'(0)=b'(0)=1$ to obtain the same result.)

\begin{theorem}
\label{thm:twotensorA1}
Let $B$ be a finite-dimensional affine crystal arising from the evaluation representation $V(j)_a\otimes V(k)_b$, with promotion operator $\pr$.  Then $B$ is one of exactly two crystals, one of which is the tensor of the Kirillov-Reshetikhin crystals $B^{1,j}\otimes B^{1,k}$, and the other of which is disconnected.  The first case may be obtained from the evaluation module at $a=b=1$, and the disconnected case may be obtained from $a=q^{-k}$ and $b=q^j$.
\end{theorem}
\begin{proof}
Since $B$ admits a promotion operator, it satisfies the conditions of Proposition~\ref{prop:crysStructure}.  The constants $a$ and $b$ must fulfill one of the two possibilities of that Proposition.  

In the first case (where $a(0)=b(0)=1$), the resulting crystal must be that of the tensor of the Kirillov-Reshetikhin crystals $B^{1,j}\otimes B^{1,k}$.  

In the second case (where $a=q^{-k}a'$ and $b=q^{-k}b'$ with $a'(0)=b'(0)=1$), Lemma~\ref{lem:e0StructConsts} may be used to describe the $0$-string of $u\otimes v$, and thus show that the resulting crystal is disconnected, as discussed above.
\end{proof}

\section{Implementation of Stembridge Local Axioms}
\label{sec:stembridge}

Computers are playing an increasingly important role in the study of algebraic combinatorics, allowing researchers to grapple with concrete examples that are far outside the realm of what may be computed by hand.  The Sage computer algebra system is a free, open-source mathematics system developed in recent years by a network of volunteers spanning the globe.  Sage ties together numerous pre-existing pieces of open-source mathematics software, using the Python programming language as a kind of glue between them.  Sage is more than the sum of these pieces, though, with numerous areas of mathematics that have been coded specifically for Sage.

In this section, we will discuss an implementation of Stembridge's local characterization of simply-laced crystals~\cite{stembridge2003} in the Sage system.  These local axioms give a finite check at each elements of a crystal basis which are sufficient to determine whether the crystal arises from a representation.  It provides functions for computations required to check the Stembridge characterization, as well as a function that carries out the full check at the level of a crystal element and for the entire crystal.

This implementation extends an existing implementation of general crystal bases in Sage, coded mostly by Prof. Anne Schilling.  

\subsection{Organization of the Crystal Code}

Sage is structured by an enhanced object-oriented framework, mimicking mathematical categories.  Object orientation is very useful for managing large projects like Sage, allowing one to replicate code across many varieties of objects and provide an organizing scheme for the project.  

In a traditional object-oriented programming language, one has \textbf{classes} which describe classes of objects with common characteristics.  Classes act as a specification for the objects which inhabit the class, providing a number of properties which every object in the class must have.  The definition of a class also often gives default values for these properties, which may be over-ridden by particular objects in the class.  Properties may be of many types, such as numbers, strings, or even functions.  Every class should also have a function for generating objects in the class, unless the class describes purely fictional objects!

As an example, we can imagine a class called \textbf{BICYCLES}, which describes my favourite mode of transportation.  The class specification for \textbf{BICYCLES} might include constants such as \textit{wheel size} and \textit{crank length}, arrays such as \textit{chain-ring sizes} and \textit{cog sizes}, two more constants specifying the current chain ring size and current cog size, and a function determining the gear-inches determined by the wheel size, crank length, and current chain ring and cog sizes.  (Gear inches determines how far the bike moves given one revolution of the pedals.)  The class specification may give default values for any or all of these values (e.g. wheel size of $700c$, for a fast road bike), which may be over-ridden by a particular instance of a bicycle (e.g. if a particular bicycle has $26''$ wheels, as most mountain bikes do).

The function for determining gear inches might also be over-ridden, which could be convenient if one has elliptical chain rings.

One may nest classes, allowing specification of objects in varying degrees of generality and also allowing \textbf{inheritance} of properties from the containing class to the contained class.  For example, one can imagine a class called \textbf{VEHICLES} which contains the classes \textbf{CARS} (describing the H\=inay\=ana) and \textbf{BICYCLES} (describing the Mah\=ay\=ana).  Code written for \textbf{VEHICLES} is inherited by \textbf{BICYCLES}, but may also be over-ridden.  For example, this could happen if the code written for \textbf{VEHICLES} works in full generality but is very slow to evaluate, while things may be very fast to calculate in the special case of \textbf{BICYCLES}.

Modern mathematics carries some inherent object orientation in the form of Category Theory.  A mathematical category consists of objects and morphisms (or `arrows') between objects, satisfying certain conditions; for example, in the category $\mathfrak{Groups}$, the objects are groups and the morphisms are group homomorphisms.   In the category $\mathfrak{Sets}$, the objects are sets and the morphisms from a set $S$ to a set $T$ are functions $f: S\rightarrow T$.  Categories may be included one into another; the category $\mathfrak{Groups}$ includes into the category $\mathfrak{Sets}$, since every group may be regarded as a set, and every group homomorphism is also a set function.  One furthermore has maps between categories, called \textbf{functors}, which map objects to objects and arrows to arrows.  The inclusion of $\mathfrak{Groups}$ into $\mathfrak{Sets}$ is an example of a functor.  (This is known as a ``forgetful functor,'' which simply forgets the extra structure of groups and group homomorphisms.)

Thus, a category is much like a class, but with additional information specifying ways of getting from one object to another, and also ways for getting from one category to another.

Crystals in Sage are defined as a category.  Objects in this category are crystals, which themselves contain crystal elements, the vertices of the crystal.  The crystal itself is called the Parent, and is endowed with various Parent methods.  These methods include a method which returns the Cartan type of the crystal, for example.  Elements of the crystal (called ``Elements'') have a number of Element methods specified, which include the crystal operators $\tilde{e_i}$ and $\tilde{f_i}$, the map to the weight lattice, and the maps $\phi_i$ and $\epsilon_i$ describing the position of the element in its $i$-string.

Since an abstract crystal graph is defined by axioms, many important functions are defined as ``abstract methods'' and left to a particular realization of a crystal basis to define.  For example, here is the code for the $e_i$ element method in the crystal code:
\begin{verbatim}
    @abstract_method
    def e(self, i):
        r"""
        Returns `e_i(x)` if it exists or ``None`` otherwise.

        This method should be implemented by the element class of
        the crystal.

        EXAMPLES::

            sage: C = Crystals().example(5)
            sage: x = C[2]; x
            3
            sage: x.e(1), x.e(2), x.e(3)
            (None, 2, None)
        """
\end{verbatim}
This does not actually do anything; the portion enclosed by the triple quotations are a documentation string describing the expected input and output of the function, and also some examples of the function in action.

This, on the other hand, is the code for the $\epsilon$ function:
\begin{verbatim}
    def epsilon(self, i):
        r"""
        EXAMPLES::

            sage: C = CrystalOfLetters(['A',5])
            sage: C(1).epsilon(1)
            0
            sage: C(2).epsilon(1)
            1
        """
        assert i in self.index_set()
        x = self
        eps = 0
        while True:
            x = x.e(i)
            if x is None:
                break
            eps = eps+1
        return eps
\end{verbatim}
Again, we have a documentation string which provides examples of the function in action.  This time we also have some actual code (which may be over-ridden by a particular realization of the crystal graph), which finds the maximum number of times on may apply $e_i$ to the given crystal element before killing it, which is the definition of the $\epsilon_i$ function for a semiregular crystal.

\subsection{Stembridge Local Axioms}

Suppose a weight representation $V$ of a quantum Kac-Moody algebra $U_q(\mathfrak{g})$ has a crystal basis.  Then the crystal basis will satisfy the axioms of a crystal graph, and may be regarded as such.  There is  still a question, though, as to when a crystal graph $\mathcal{B}$ actually arises from a weight representation of $U_q(\mathfrak{g})$.  Stembridge gave an answer to this question in the case where $\mathcal{B}$ is simply laced via a local characterization of the crystal $\mathcal{B}$~\cite{stembridge2003}.  In particular, Stembridge provided a list of axioms which, if fulfilled by a highest-weight crystal graph $\mathcal{B}$, imply that $\mathcal{B}$ may be obtained from a highest-weight representation of $U_q(\mathfrak{g})$.

The first two axioms are as follows:
\begin{itemize}
\item (P1) The $i$-strings in the crystal are of finite length, and cycle-free.
\item (P2) For any $x, y \in \mathcal{B}$, there is at most one $i$-colored edge $x\rightarrow y$.
\end{itemize}
These are direct consequences of Definition~\ref{def:crystalGraph}.

To write the additional axioms, we first define some additional functions on the crystal.
\begin{definition}
We define four functions $\mathcal{B}\rightarrow \ZZ$:
\begin{enumerate}
\item The \textbf{delta-depth operator} $\Delta \epsilon (x,i,j) = - \epsilon_j \tilde{e_i} x + \epsilon_j x$,
\item The \textbf{delta-rise operator} $\Delta \phi (x,i,j) =  \phi_j \tilde{e_i} x - \phi_j x$,
\item The \textbf{del-depth operator} $\Delta \epsilon (x,i,j) = \epsilon_j \tilde{f_i} x - \epsilon_j x$,
\item The \textbf{del-rise operator} $\Delta \phi (x,i,j) = - \phi_j \tilde{f_i} x + \phi_j x$,
\end{enumerate}
\end{definition}

Also from Definition~\ref{def:crystalGraph}, we have that for any $x \in \mathcal{B}$:
\[
\phi_i x = \epsilon_i x + \langle h_i, \wt(x) \rangle.
\]
As a result, we have:
\begin{eqnarray*}
\Delta \epsilon (x,i,j) + \Delta \phi (x,i,j) & = & 
    - \epsilon_j \tilde{e_i} x + \epsilon_j x + \phi_j \tilde{e_i} x - \phi_j x \\
     & = & \langle h_j, \wt(\tilde{e_i} x) \rangle - \langle h_j, \wt(x) \rangle \\
     & = & \langle h_j, \wt(\tilde{e_i} x) - \wt(x) \rangle \\
     & = & \langle h_j, \alpha_i \rangle \\
     & = & A_{ij}. \\
\end{eqnarray*}
This is the third Stembridge axiom:
\begin{itemize}
\item (P3) $\Delta \epsilon (x,i,j) + \Delta \phi (x,i,j) = A_{ij}$.
\end{itemize}

The remaining axioms do not follow directly from the definition of the crystal graph.  The fourth is:
\begin{itemize}
\item (P4) $\Delta \epsilon (x,i,j)\leq 0$ and $\Delta \phi (x,i,j) \leq 0$, when $\tilde{e_i} x\neq \emptyset$.
\end{itemize}
We define the \textbf{Stembridge triple} for a crystal element $x$ to be the tuple $( A_{ij}, \Delta \epsilon (x,i,j), \Delta \phi (x,i,j))$, defined when $\tilde{e_j} x\neq \emptyset$.  In a simply-laced crystal, we have $A_{ij} \in \{0, -1\}$ when $i\neq j$.  As a result, the Stembridge triple must be one of:
\[
(0, 0, 0), \phantom{aaa} (-1, 0, -1), \phantom{aaa} (-1, -1, 0).
\]

Axioms (P5) and (P6) deal with these cases separately.  Assume that $\tilde{e_i} x\neq \emptyset$ and $\tilde{e_j} x\neq \emptyset$.
\begin{itemize}
\item (P5) If $\Delta \epsilon (x,i,j)=0$, then:
\[
\tilde{e_j}\tilde{e_i}x=\tilde{e_i}\tilde{e_j}x.
\]
\item (P6) If $\Delta \epsilon (x,i,j)=\Delta \epsilon (x,j,i)=-1$, then:
\[
\tilde{e_i}\tilde{e_j}^2\tilde{e_i}x=\tilde{e_j}\tilde{e_i}^2\tilde{e_j}x.
\]
\end{itemize}

There is also a collection of dual axioms, which add no additional information.

\subsection{The Code.}

The following functions are defined on elements of a crystal graph; the ``self'' in the code refers to the particular crystal element.  Each function contains a documentation string, which is automatically compiled into the documentation for Sage; this documentation describes the purpose of the function, as well as its expected inputs and outputs.  In the documentation string are a number of examples.

\begin{verbatim}
def stembridgeDelta_depth(self,i,j):
    r"""
    The `i`-depth of a crystal node `x` is ``-x.epsilon(i)``.
    This function returns the difference in the `j`-depth of `x` and 
    ``x.e(i)``, where `i` and `j` are in the index set of the underlying
    crystal.  This function is useful for checking the Stembridge local 
    axioms for crystal bases.

    EXAMPLES::

        sage: T = CrystalOfTableaux(['A',2], shape=[2,1])
        sage: t=T(rows=[[1,2],[2]])
        sage: t.stembridgeDelta_depth(1,2)
        0
        sage: s=T(rows=[[2,3],[3]])
        sage: s.stembridgeDelta_depth(1,2)
        -1
    """
    if self.e(i) is None: return 0
    return -self.e(i).epsilon(j) + self.epsilon(j)

def stembridgeDelta_rise(self,i,j):
    r"""
    The `i`-rise of a crystal node `x` is ``x.phi(i)``.

    This function returns the difference in the `j`-rise of `x` and 
    ``x.e(i)``, where `i` and `j` are in the index set of the underlying 
    crystal.  This function is useful for checking the Stembridge local 
    axioms for crystal bases.

    EXAMPLES::

        sage: T = CrystalOfTableaux(['A',2], shape=[2,1])
        sage: t=T(rows=[[1,2],[2]])
        sage: t.stembridgeDelta_rise(1,2)
        -1
        sage: s=T(rows=[[2,3],[3]])
        sage: s.stembridgeDelta_rise(1,2)
        0
    """
    if self.e(i) is None: return 0
    return self.e(i).phi(j) - self.phi(j)

def stembridgeDel_depth(self,i,j):
    r"""
    The `i`-depth of a crystal node `x` is ``-x.epsilon(i)``.
    This function returns the difference in the `j`-depth of `x` and 
    ``x.f(i)``, where `i` and `j` are in the index set of the underlying 
    crystal.  This function is useful for checking the Stembridge local 
    axioms for crystal bases.

    EXAMPLES::
    
        sage: T = CrystalOfTableaux(['A',2], shape=[2,1])
        sage: t=T(rows=[[1,1],[2]])
        sage: t.stembridgeDel_depth(1,2)
        0
        sage: s=T(rows=[[1,3],[3]])
        sage: s.stembridgeDel_depth(1,2)
        -1
    """
    if self.f(i) is None: return 0
    return -self.epsilon(j) + self.f(i).epsilon(j)

def stembridgeDel_rise(self,i,j):
    r"""
    The `i`-rise of a crystal node `x` is ``x.phi(i)``.
    This function returns the difference in the `j`-rise of `x` and 
    ``x.f(i)``, where `i` and `j` are in the index set of the underlying 
    crystal.  This function is useful for checking the Stembridge local 
    axioms for crystal bases.

    EXAMPLES::

        sage: T = CrystalOfTableaux(['A',2], shape=[2,1])
        sage: t=T(rows=[[1,1],[2]])
        sage: t.stembridgeDel_rise(1,2)
        -1
        sage: s=T(rows=[[1,3],[3]])
        sage: s.stembridgeDel_rise(1,2)
        0
    """
    if self.f(i) is None: return 0
    return self.phi(j)-self.f(i).phi(j)

def stembridgeTriple(self,i,j):
    r"""
    Let `A` be the Cartan matrix of the crystal, `x` a crystal element, 
    and let `i` and `j` be in the index set of the crystal.
    Further, set
    ``b=stembridgeDelta_depth(x,i,j)``, and
    ``c=stembridgeDelta_rise(x,i,j))``.
    If ``x.e(i)`` is non-empty, this function returns the triple 
    `( A_{ij}, b, c )`; otherwise it returns ``None``.  By the Stembridge
    local characterization of crystal bases, one should have `A_{ij}=b+c`.

    EXAMPLES::

        sage: T = CrystalOfTableaux(['A',2], shape=[2,1])
        sage: t=T(rows=[[1,1],[2]])
        sage: t.stembridgeTriple(1,2)
        sage: s=T(rows=[[1,2],[2]])
        sage: s.stembridgeTriple(1,2)
        (-1, 0, -1)

        sage: T = CrystalOfTableaux(['B',2], shape=[2,1])
        sage: t=T(rows=[[1,2],[2]])
        sage: t.stembridgeTriple(1,2)
        (-2, 0, -2)
        sage: s=T(rows=[[-1,-1],[0]])
        sage: s.stembridgeTriple(1,2)
        (-2, -2, 0)
        sage: u=T(rows=[[0,2],[1]])  
        sage: u.stembridgeTriple(1,2)
        (-2, -1, -1)
    """
    if self.e(i) is None: return None
    A=self.cartan_type().cartan_matrix()
    b=self.stembridgeDelta_depth(i,j)
    c=self.stembridgeDelta_rise(i,j)
    dd=self.cartan_type().dynkin_diagram()
    a=dd[j,i]
    return (a, b, c)

def _test_stembridge_local_axioms(self, index_set=None, verbose=False, \\
    **options):
    r"""
    This implements tests for the Stembridge local characterization on the
    element of a crystal ``self``.  The current implementation only uses 
    the axioms for simply-laced types.  Crystals of other types should 
    still pass the test, but in non-simply-laced types, passing is not a 
    guarantee that the crystal arises from a representation.
    
    One can specify an index set smaller than the full index set of the 
    crystal, using the option ``index_set``.

    Running with ``verbose=True`` will print warnings when a test fails.

    REFERENCES::

        .. [S2003] John R. Stembridge, A Local Characterization of 
           Simply-Laced Crystals,
           Transactions of the American Mathematical Society, Vol. 355, 
           No. 12 (Dec., 2003), pp. 4807-4823

    EXAMPLES::

        sage: T = CrystalOfTableaux(['A',2], shape=[2,1])
        sage: t=T(rows=[[1,1],[2]])
        sage: t._test_stembridge_local_axioms()
        True
        sage: t._test_stembridge_local_axioms(index_set=[1,3])
        True
        sage: t._test_stembridge_local_axioms(verbose=True)
        True
    """
    tester = self._tester(**options)
    goodness=True
    A=self.cartan_type().cartan_matrix()
    if index_set is None: index_set=self.index_set()

    for (i,j) in Subsets(index_set, 2):
        if self.e(i) is not None and self.e(j) is not None:
            triple=self.stembridgeTriple(i,j)
            #Test axioms P3 and P4.
            if not triple[0]==triple[1]+triple[2] or \\
                triple[1]>0 or triple[2]>0:
                if verbose:
                    print 'Warning: Failed axiom P3 or P4 at vector ', \\
                    self, 'i,j=', i, j, 'Stembridge triple:', \\
                    self.stembridgeTriple(i,j)
                    goodness=False
                else:
                    tester.fail()
            if self.stembridgeDelta_depth(i,j)==0:
                #check E_i E_j(x)= E_j E_i(x)
                if self.e(i).e(j)!=self.e(j).e(i) or \\
                self.e(i).e(j).stembridgeDel_rise(j, i)!=0:
                    if complete:
                        print 'Warning: Failed axiom P5 at: vector ', \\
                        self, 'i,j=', i, j, 'Stembridge triple:', \\
                        stembridgeTriple(x,i,j)
                        goodness=False
                    else:
                        tester.fail()
            if self.stembridgeDelta_depth(i,j)==-1 and \\
            self.stembridgeDelta_depth(j,i)==-1:
                #check E_i E_j^2 E_i (x)= E_j E_i^2 E_j (x)
                y1=self.e(j).e(i).e(i).e(j)
                y2=self.e(j).e(i).e(i).e(j)
                a=y1.stembridgeDel_rise(j, i)
                b=y2.stembridgeDel_rise(i, j)
                if y1!=y2 or a!=-1 or b!=-1:
                    if verbose:
                        print 'Warning: Failed axiom P6 at: vector ', x,\\
                        'i,j=', i, j, 'Stembridge triple:', \\
                        stembridgeTriple(x,i,j)
                        goodness=False
                    else:
                        tester.fail()
    tester.assertTrue(goodness)
    return goodness
\end{verbatim}

The following function is defined on any finite crystal graph, and simply checks that the Stembridge local axioms hold on each element of the crystal.

\begin{verbatim}
def _test_stembridge_local_axioms(self, index_set=None, verbose=False, \\
complete=False, **options):
    r"""
    This implements tests for the Stembridge local characterization on the
    finite crystal ``self``.  The current implementation only uses the 
    rules for simply-laced types.  Crystals of other types should still 
    pass the test, but expansion of this test to non-simply laced type
    would be desirable.
    
    One can specify an index set smaller than the full index set of the 
    crystal, using the option ``index_set``.
    
    Running with ``verbose=True`` will print each node for which a local
    axiom test applies.
    
    Running with ``complete=True`` will continue to run the test past the 
    first failure of the local axioms.  This is probably only useful in 
    conjunction with the verbose option, to see all places where the local
    axioms fail.

    EXAMPLES::

        sage: T = CrystalOfTableaux(['A',3], shape=[2,1])
        sage: T._test_stembridge_local_axioms()
        True
        sage: T._test_stembridge_local_axioms(verbose=True)
        True
        sage: T._test_stembridge_local_axioms(index_set=[1,3])
        True
    """
    tester = self._tester(**options)
    goodness=True

    for x in self:
        goodness=x._test_stembridge_local_axioms(index_set, verbose)
        if goodness==False and not complete:
            tester.fail()
    tester.assertTrue(goodness)
    return goodness
\end{verbatim}